\documentclass[11pt,twoside]{article}
\usepackage{fullpage}
\usepackage{epsf}
\usepackage{fancyhdr}
\usepackage{graphics}
\usepackage{graphicx}
\usepackage{psfrag}
\usepackage{microtype}
\usepackage{subfigure}
\usepackage{algorithmic}
\usepackage{color,xcolor}
\usepackage[linesnumbered,ruled]{algorithm2e}% http://ctan.org/pkg/algorithm2e
\DontPrintSemicolon
\usepackage{color}
\usepackage{tabularx}
\usepackage{amsthm}
\usepackage{amsfonts}
\usepackage{amsmath}
\usepackage{amssymb,bbm}
\usepackage{stackengine}
\usepackage{footnote}
\makesavenoteenv{tabular}
\makesavenoteenv{table}

\newcommand{\bali}{$$\begin{aligned}}
\newcommand{\eali}{\end{aligned}$$}

\renewcommand{\tilde}{\widetilde}

\newcommand{\aone}{\ensuremath{\alpha_1}}
\newcommand{\atwo}{\ensuremath{\alpha_2}}
\newcommand{\Tarb}{{T^*}}
\newcommand{\Quant}{\ensuremath{Q}}
\newcommand{\sglip}{\ell_\Xi}

\newcommand{\remnoise}{\zeta}

\newcommand{\filtration}{\mathcal{F}}

\newcommand{\holderexponent}{\gamma}

\newcommand{\holderconst}{L_\holderexponent}
\newcommand{\holderconstprime}{L_1}

\newcommand{\quadvar}[2]{[ #1 ]_{#2}}
\newcommand{\crossvar}[3]{[ #1, #2 ]_{#3}}

\newcommand{\testMat}{G}

\newcommand{\LSN}{LSN\xspace}
\newcommand{\ISC}{\textbf{ISC}\xspace}

\newcommand{\norm}[1]{\vecnorm{#1}{2}}

\newcommand{\epochsrestart}{B}
\newcommand{\shortepoch}{\niters^\flat}

\newcommand{\stpszmax}{{\stpsz_{max}}}
\newcommand{\kappamax}{{\omega_{max}}}

\newcommand{\sigmastar}{{\sigma_*}}
\newcommand{\sigmastarsq}{\ensuremath{\sigma_*^2}}

\newcommand{\niters}{T}

\usepackage{thmtools, thm-restate}

\newcommand{\stpsz}{\eta}

\newcommand{\highorder}{\mathcal{H}}
\newcommand{\highordertilde}{\widetilde{\highorder}}

\newcommand{\sigstar}{\sigma_*}

\newcommand{\real}{\ensuremath{\mathbb{R}}}

\newcommand{\thetastar}{\ensuremath{{\theta^*}}}
\newcommand{\thetahat}{\ensuremath{\widehat{\theta}}}
\newcommand{\abss}[1]{\left| #1 \right |}

\newcommand{\simiid}{\overset{\mathrm{i.i.d.}}{\sim}}

\newcommand{\sphere}{\ensuremath{\mathbb{S}}}

\newcommand{\mydefn}{\ensuremath{:=}}

\newcommand{\noise}{\ensuremath{\varepsilon}}

\newcommand{\smoothness}{L} % smoothness parameter

\newcommand{\defn}{:=}

\newcommand{\matsnorm}[2]{|\!|\!| #1 | \! | \!|_{{#2}}}
\newcommand{\vecnorm}[2]{\left\| #1\right\|_{#2}}
\newcommand{\opnorm}[1]{\ensuremath{\matsnorm{#1}{\tiny{\mbox{op}}}}}
\newcommand{\inprod}[2]{\ensuremath{\langle #1 , \, #2 \rangle}}

\newcommand{\Exs}{\ensuremath{{\mathbb{E}}}}
\newcommand{\Prob}{\ensuremath{{\mathbb{P}}}}

\newtheoremstyle{named}{}{}{\itshape}{}{\bfseries}{.}{.5em}{\thmnote{#3's }#1}
\theoremstyle{named}

\theoremstyle{plain}

\newtheorem{theorem}{Theorem}

\newtheorem{lemma}{Lemma}

\newlength{\widebarargwidth}
\newlength{\widebarargheight}
\newlength{\widebarargdepth}

\makeatletter
\long\def\@makecaption#1#2{
        \vskip 0.8ex
        \setbox\@tempboxa\hbox{\small {\bf #1:} #2}
        \parindent 1.5em  %% How can we use the global value of this???
        \dimen0=\hsize
        \advance\dimen0 by -3em
        \ifdim \wd\@tempboxa >\dimen0
                \hbox to \hsize{
                        \parindent 0em
                        \hfil
                        \parbox{\dimen0}{\def\baselinestretch{0.96}\small
                                {\bf #1.} #2
                                }
                        \hfil}
        \else \hbox to \hsize{\hfil \box\@tempboxa \hfil}
        \fi
        }
\makeatother

\long\def\comment#1{}
\definecolor{battleshipgrey}{rgb}{0.52, 0.52, 0.51}
\definecolor{darkgray}{rgb}{0.66, 0.66, 0.66}
\definecolor{darkgreen}{rgb}{0.0, 0.2, 0.13}
\definecolor{darkspringgreen}{rgb}{0.09, 0.45, 0.27}
\definecolor{dukeblue}{rgb}{0.0, 0.0, 0.61}
\definecolor{olivedrab7}{rgb}{0.24, 0.2, 0.12}
\definecolor{darkblue}{rgb}{0.0, 0.0, 0.55}
\definecolor{darkscarlet}{rgb}{0.34, 0.01, 0.1}
\definecolor{candyapplered}{rgb}{1.0, 0.03, 0.0}
\definecolor{ao(english)}{rgb}{0.0, 0.5, 0.0}
\definecolor{applegreen}{rgb}{0.55, 0.71, 0.0}

\usepackage{url}% for url's in bib
\usepackage[colorlinks,linkcolor=magenta,citecolor=blue, pagebackref=true,backref=true]{hyperref}
\renewcommand*{\backrefalt}[4]{%
    \ifcase #1 \footnotesize{(Not cited.)}%
    \or        \footnotesize{(Cited on page~#2.)}%
    \else      \footnotesize{(Cited on pages~#2.)}%
    \fi}

\usepackage{nicefrac}
\usepackage{comment}

\usepackage{chngpage}

 \usepackage{tabularx}%

\usepackage{enumitem}
\usepackage{booktabs}
\usepackage{caption}

\usepackage{bm,bbm}
\usepackage{mathtools}

\newcommand{\strongconvex}{\mu}
\newtheorem{assumption}{Assumption}

\theoremstyle{definition}

\newtheorem{theorem2}{Proposition}
\newtheorem{corollaryprime}[theorem]{Corollary}%opening
\newcommand{\lammin}{\ensuremath{\lambda_{\operatorname{min}}}}

\newcommand{\SigStar}{\ensuremath{\Sigma^*}}

\newcommand{\NoiseAplain}{\ensuremath{\Xi}}
\newcommand{\NoiseAt}{\ensuremath{\NoiseAplain_t}}

\makeatletter
\long\def\@makecaption#1#2{
        \vskip 0.8ex
        \setbox\@tempboxa\hbox{\small {\bf #1:} #2}
        \parindent 1.5em  %% How can we use the global value of this???
        \dimen0=\hsize
        \advance\dimen0 by -3em
        \ifdim \wd\@tempboxa >\dimen0
                \hbox to \hsize{
                        \parindent 0em
                        \hfil 
                        \parbox{\dimen0}{\def\baselinestretch{0.96}\small
                                {\bf #1.} #2
                                } 
                        \hfil}
        \else \hbox to \hsize{\hfil \box\@tempboxa \hfil}
        \fi
        }
\makeatother

\newcommand{\normb}[1]{\vecnorm{#1}{2}^2}

\newcommand{\usedim}{\ensuremath{d}}

\newcommand{\Sphere}[1]{\ensuremath{\mathbb{S}}}

\newcommand{\cF}{{\mathcal{F}}}

\newcommand{\burnin}{{T_0}}

\def\ROOTSGD{\textsf{ROOT-SGD}\xspace}

\newcommand{\hessianstar}{H^*}
\newcommand{\Lmax}{{L_{max}}}
\newcommand{\Lmaxsq}{{L_{max}^2}}

\newcommand{\mprob}{\ensuremath{\mathbb{P}}}

\newcommand{\sglipfour}{\sglip}
\newcommand{\sglipfourbracket}{\sglipfour}

\newcommand{\XiDomain}{\Xi}

\usepackage{mathrsfs}
\newcommand{\naF}[1]{\nabla F (#1)}
\newcommand{\Epoch}{B}
\newcommand{\inneriters}{\shortepoch}
\newcommand{\numobs}{n}
\newcommand{\nsamples}{{n}}%COMBINE LATER
\newcommand{\sglipsq}{\sglip^2}
\newcommand{\binprod}[2]{\ensuremath{ \big \langle #1 , \, #2 \big \rangle}}
\newcommand{\Tburnin}{{\burnin}}
\newcommand{\nbF}[1]{\nabla^2 F(#1)}
\newcommand{\sgliptild}{\ell_\Xi}
\newcommand{\sigstartild}{\tilde{\sigstar}}


\newcounter{parentnumber}

\newcommand{\kcon}{\ensuremath{c_0}}
\newcommand{\Contwo}{\ensuremath{c_1}}
\newcommand{\aconprime}{a'}
\newcommand{\para}[1]{\left( #1 \right)}
\newcommand{\normd}[1]{\vecnorm{#1}{2}^4}
\newcommand{\Twindow}{\tilde{T}^*}
\newcommand{\ISCspace}{\ISC}%SHOULD REMOVE
\newcommand{\noiset}{\noise_t}
\newcommand{\normc}[1]{\vecnorm{#1}{2}^3}
\newcommand{\nafxit}[1]{\nabla f(#1; \xi_t)}

\newcommand{\LSNspace}{\LSN}

\newcommand{\asympconst}{a}

\def\pb{}
\begin{document}

\begin{center}
{\bf{\LARGE{\ROOTSGD: Sharp Nonasymptotics and Asymptotic Efficiency in a Single Algorithm}}}

\vspace*{.2in}
{\large{
\begin{tabular}{cccc}
Chris Junchi Li$^{\diamond, \star}$
&
Wenlong Mou$^{ \diamond, \star}$
&
Martin J.~Wainwright$^{\diamond, \dagger}$ 
&
Michael I.~Jordan$^{\diamond, \dagger}$
\end{tabular}
}}

\vspace*{.2in}

\begin{tabular}{c}
Department of Electrical Engineering and Computer Sciences$^\diamond$
\\
Department of Statistics$^\dagger$
\\
UC Berkeley
\\
\end{tabular}

\vspace*{.2in}

\today

\vspace*{.2in}

\let\thefootnote\relax\footnotetext{$\,^\star$
CJL and WM contributed equally to this work.}

\begin{abstract}
We study the problem of solving strongly convex and smooth unconstrained optimization problems using stochastic first-order algorithms. We devise a novel algorithm, referred to as \emph{Recursive One-Over-T SGD} (\ROOTSGD), based on an easily implementable, recursive averaging of past stochastic gradients. We prove that it simultaneously achieves state-of-the-art performance in both a finite-sample, nonasymptotic sense and an asymptotic sense. On the nonasymptotic side, we prove risk bounds on the last iterate of \ROOTSGD with leading-order terms that match the optimal statistical risk with a unity pre-factor, along with a higher-order term that scales at the sharp rate of $O(n^{-3/2})$ under the Lipschitz condition on the Hessian matrix. On the asymptotic side, we show that when a mild, one-point Hessian continuity condition is imposed, the rescaled last iterate of (multi-epoch) \ROOTSGD converges asymptotically to a Gaussian limit with the Cram\'{e}r-Rao optimal asymptotic covariance, for a broad range of step-size choices.
\end{abstract}

\end{center}


\vspace*{.2in}

\pb\section{Introduction}\label{sec_intro}
Let $f: \real^d \times \XiDomain \rightarrow \real$ be differentiable as a function of its first argument, and consider the following unconstrained minimization problem:
\begin{align}\label{stoch_opt}
\min_{\theta \in \real^d} F(\theta)
, \qquad \mbox{where $F(\theta) \mydefn \Exs \big[ f (\theta; \xi) \big]$,}
\end{align}
and where the expectation is taken over a random vector $\xi \in \XiDomain$ with distribution $\mprob$. 
Our goal is to approximately solve this minimization problem based on samples $(\xi_i)_{i = 1,2,  \cdots} \simiid \mprob$, and moreover to do so in a way that is computationally efficient and statistically optimal. 
When the samples arrive as an online stream, it is desirable to compute the approximate solution in a single pass, without storing the data, and this paper focuses on this online setting.

Stochastic optimization problems of this type underpin a variety of methods in large-scale machine learning and statistical inference. 
One of the simplest methods is \emph{stochastic gradient descent} (SGD), which recursively updates a parameter vector $\theta_t$ by taking a step in the direction of a single stochastic gradient, with a (possibly) time-varying step-size $\stpsz_t$ \cite{robbins1951stochastic}. 
This simple strategy has been surprisingly successful in modern large-scale statistical machine learning problems~\cite{nemirovski2009robust,bottou2018optimization,nguyen2019new}; however, it can be substantially improved, both in theory and in practice, by algorithms that make use of more than a single stochastic gradient.
Such algorithms belong to the general family of \emph{stochastic first-order methods}.  
Various procedures have been studied, involving different weightings of past stochastic gradients, and also a range of analysis techniques. 
The diversity of approaches is reflected by the wide range of terminology, including \emph{momentum}, \emph{averaging}, \emph{acceleration}, and \emph{variance reduction}. 
All of these ideas center around two main underlying goals---that of proceeding quickly to a minimum, and that of arriving at a final state that achieves the optimal statistical efficiency and also provides a calibrated assessment of the uncertainty associated with the solution.

More concretely, the former goal requires the algorithm to achieve a fast rate of convergence and low sample complexity, ideally matching that of the noiseless case and the information-theoretic limit. 
For example, gradient descent takes $O\big( \frac{\smoothness}{\strongconvex}\big)$ number of iterations to optimize a $\smoothness$-smooth and $\strongconvex$-strongly convex function. 
It is therefore desirable that the sample-size requirement for a stochastic optimization algorithm scales linearly with $O\big(\frac{\smoothness}{\strongconvex})$, with additional terms characterizing the effect of random noise on optimality. 
On the other hand, the latter goal imposes a more fine-grained requirement on the estimator produced by the algorithm. 
Roughly speaking, we need the estimator to share the same \emph{optimal statistical properties} typically possessed by the empirical risk minimizer (were it be computed exactly in the batch setting).  
The notion of \emph{statistical efficiency}, in both its asymptotic and nonasymptotic forms, allows for a fine-grained study of these issues.

Let $\thetastar$ denote the minimizer of $F$, and define the matrices
\begin{align*}
\hessianstar \mydefn \nabla^2 F (\thetastar)
, \quad \mbox{and} \quad
\SigStar \mydefn \Exs \left[\nabla f (\thetastar; \xi) \nabla f(\thetastar; \xi)^\top \right]
.
\end{align*}
Under certain regularity assumptions, given a collection of $n$ samples $(\xi_i)_{i \in [n]} \simiid \mprob$, classical statistical theory guarantees that the minimizer $\thetahat_\numobs^{\mathrm{ERM}} \mydefn \arg\min_{\theta \in \real^\usedim} \sum_{i = 1}^\numobs f (\theta; \xi_i) $ of the associated empirical risk has the following asymptotic behavior:
\begin{align}
\sqrt{\numobs} \left( \thetahat_\numobs^{\mathrm{ERM}} - \thetastar \right) 
\xrightarrow{d} 
\mathcal{N} \left(0, (\hessianstar)^{-1} \SigStar (\hessianstar)^{-1} \right)
.\label{eq:erm-optimal-asymptotic}
\end{align}
Furthermore, the asymptotic distribution~\eqref{eq:erm-optimal-asymptotic} is known to be locally asymptotic minimax, i.e.~given a bowl-shaped loss function, the asymptotic risk of \emph{any} estimator is lower bounded by the expectation under such a Gaussian distribution, in a suitably defined sequence of local neighborhoods.  
See, for example,~\cite{duchi2021asymptotic} for a precise statement.

Unfortunately, the goals of rapid finite-sample convergence and optimal asymptotic behavior are in tension, and the literature has not yet arrived at a single algorithmic framework that achieves both goals simultaneously.  
Consider in particular two seminal lines of research:

\begin{enumerate}[leftmargin=0mm,label=(\arabic*)]
%\bcar
\item
The Polyak-Ruppert-Juditsky (PRJ) procedure~\cite{polyak1992acceleration,polyak1990new,ruppert1988efficient} incorporates slowly diminishing step-sizes into SGD, thereby achieving asymptotic normality with an optimal covariance matrix (and unity pre-factor).  
This meets the goal of calibrated uncertainty.  
However, the PRJ procedure is \emph{not} optimal from a nonasymptotic point of view: rather, it suffers from large high-order nonasymptotic terms and fails to achieve the optimal sample complexity in general~\cite{bach2011non}.
\item 
On the other hand, variance-reduced stochastic optimization methods have been designed to achieve reduced sample complexity that is the sum of a \emph{statistical error} and an \emph{optimization error}~\cite{roux2012stochastic, shalev2013stochastic, johnson2013accelerating, lei2017less,  defazio2014saga}.  
These methods yield control on the optimization error, with sharp nonasymptotic rates of convergence, but the guarantees for the statistical error term yield an asymptotic behavior involving constant pre-factors that are strictly greater than unity, and is hence sub-optimal.
% \ecar
\end{enumerate}

\paragraph{An open question:}  
Given this state of affairs, we are naturally led to the following question: can a single stochastic optimization algorithm simultaneously achieve optimal asymptotic and nonasymptotic guarantees? 
In particular, we would like such guarantees to enjoy the fine-grained statistical properties satisfied by the empirical risk minimizer, for a commensurate set of assumptions on the function $F$ and the observations $f(\cdot; \xi)$ and including the same rate of decay of high-order terms.

In this paper, we resolve this open question, in particular by proposing and analyzing a novel algorithm called \emph{Recursive One-Over-T Stochastic Gradient Descent} (\ROOTSGD).  
It is very easy to describe and implement, and we prove that it is optimal in both asymptotic and nonasymptotic senses:

\begin{enumerate}[leftmargin=0mm,label=(\arabic*)]
%\bcar
\item 
On the nonasymptotic side, under suitable smoothness assumptions, we show that the estimator $\thetahat_{\numobs}^{\mathrm{ROOT}}$ produced by the last iterate of the (multi-epoch) \ROOTSGD satisfies a bound of the following form:
\begin{align}
\Exs \Vert \thetahat_{\numobs}^{\mathrm{ROOT}} - \thetastar \Vert_2^2 
\leq 
\frac{\mathrm{Tr} \left( (\hessianstar)^{-1} \SigStar (\hessianstar)^{-1} \right)}{\numobs}
+ 
O\left( \frac{1}{\numobs^{3/2}} \right)
.\label{eq:matching-cramer-rao-informal}
\end{align}
Note that the leading-order term of the bound~\eqref{eq:matching-cramer-rao-informal} is exactly the squared norm of the Gaussian random vector in the local asymptotic minimax limit or Cram\'{e}r-Rao lower bound, with unity pre-factor. 
Moreover, our bound is entirely nonasymptotic, valid for all finite $\numobs$.  
We also prove that high-order term $O \big( \numobs^{-3/2} \big)$ is unavoidable under a natural setup, and it improves upon existing $O (n^{-7/6})$ and $O(n^{-5/4})$ rates for the PRJ procedure~\cite{bach2011non,xu2011towards,gadat2017optimal}. We also derive
similar bounds for the objective gap $F\big(\thetahat_\numobs^{\mathrm{ROOT}}\big) - F (\thetastar)$ and the gradient norm $\Vert \nabla F\big(\thetahat_\numobs^{\mathrm{ROOT}}\big) \Vert^2$.

\item 
Furthermore, the nonasymptotic bound~\eqref{eq:matching-cramer-rao-informal} holds true under a mild sample-size requirement.  
Indeed, given a $\smoothness$-smooth and $\strongconvex$-strongly convex population-level function $F$, and assuming that the noise $\noise (\cdot; \xi) \mydefn \nabla f (\cdot; \xi) - \nabla F (\cdot)$ satisfies a stochastic Lipschitz condition with parameter $\sglip$, the finite-sample bounds are viable as long as $\numobs \gtrsim \frac{\smoothness}{\strongconvex} + \frac{\sglip^2}{\strongconvex^2}$. The first term $O \big( \frac{\smoothness}{\strongconvex} \big)$ matches the iteration complexity of gradient descent, and the $O \big( \frac{\sglip^2}{\strongconvex^2} \big)$ term is the sample complexity needed for distinguishing a $\strongconvex$-strongly convex function from a constant function. 
The high-order terms in Eq.~\eqref{eq:matching-cramer-rao-informal} also depend on the parameters $(\strongconvex, \smoothness, \sglip)$ in a similar way. 
This exhibits the fast nonasymptotic convergence of our algorithm, matching state-of-the-art variance reduction algorithms.

\item 
We also establish asymptotic guarantees for the \ROOTSGD algorithm.
Assuming insetad a mild one-point Hessian continuity condition at the minimizer, for a broad range of step-size choices the last iterate $\thetahat_\numobs^{\mathrm{ROOT}}$ converges in distribution to the optimal Gaussian law~\eqref{eq:erm-optimal-asymptotic} whenever the Hessian matrix $\nabla^2 F$ is continuous at $\thetastar$, a much weaker condition manifesting the difference between \ROOTSGD and the Polyak-Ruppert averaging procedure~\cite{polyak1992acceleration}.  
%      \ecar
\end{enumerate}
\noindent
Notably, both the MSE bound of the form~\eqref{eq:matching-cramer-rao-informal} and the asymptotic normality are fine-grained guarantees that are satisfied by the empirical risk minimizer, under comparable assumption posed on the continuity of Hessian matrix.
To the best of our knowledge, such guarantees have not been available heretofore in the literature on stochastic optimization. 
The \ROOTSGD algorithm proposed in this paper achieves these guarantees not only simultaneously, but also with sharp nonasymptotic sample complexity.

The rest of the paper is organized as follows.  
We present the \ROOTSGD algorithm in \S\ref{sec_ROOTSGD}, and delineate the asymptotic normality and nonasymptotic upper bounds in \S\ref{sec_algoresult_hessiansmooth}.  
We present our conclusions in \S\ref{sec_summaryfuture}.  
Full proofs and discussions are provided in the appendix.

\paragraph{Notations.} 
Given a pair of vectors $u, v \in \real^d$, we write $\binprod{u}{v}$ for the inner product, and $\vecnorm{v}{2}$ for the Euclidean norm.
For a matrix $M$, the $\ell_2$-operator norm is defined as $\opnorm{M}\mydefn \sup_{\vecnorm{v}{2} = 1} \vecnorm{M v}{2}$.  
For scalars $a, b \in \real$, we adopt the shorthand notation $a \land b \mydefn \min(a,b)$ and $a\lor b\mydefn \max(a,b)$.  
Throughout the paper, we use the $\sigma$-fields $\filtration_t \mydefn \sigma (\xi_1, \xi_2, \cdots, \xi_t)$ for any $t \geq 0$.  
Unless indicated otherwise, $C$ denotes some positive, universal constant whose value may change at each appearance.  
For two sequences $\left\{a_{n}\right\}$ and $\left\{b_{n}\right\}$ of positive scalars, we denote $a_{n} \gtrsim b_{n}$ (resp.~$a_{n} \lesssim b_{n}$) if $a_{n} \geq C b_{n}$ (resp.~$a_{n} \leq C b_{n}$) for all $n$, and $a_{n} \asymp b_{n}$ if $a_{n} \gtrsim b_{n}$ and $a_{n} \lesssim b_{n}$ hold simultaneously.  
We also write $a_{n}=O\left(b_{n}\right), a_{n}=\Theta\left(b_{n}\right), a_{n}=\Omega\left(b_{n}\right)$ as $a_{n} \lesssim b_{n}, a_{n} \asymp b_{n}, a_{n} \gtrsim b_{n}$, respectively.

We finally introduce some martingale-related notations.
Given vector-valued martingales $(X_t)_{t \geq \burnin}, (Y_t)_{t \geq \burnin}$ adapted to the filtration $(\filtration_t)_{t \geq \burnin}$, we use the following notation for cross variation for $t \geq \burnin$:
\begin{align*}
\crossvar{X}{Y}{t} \mydefn \sum_{s = \burnin + 1}^{t} \binprod{ X_t - X_{t - 1}}{Y_t - Y_{t - 1} }.
\end{align*}
We also define $\quadvar{X}{t} \mydefn \crossvar{X}{X}{t}$ to be the quadratic variation of the process $(X_t)_{t \geq \burnin}$.

\pb\section{Constructing the \ROOTSGD algorithm}\label{sec_ROOTSGD} 
In this section, we introduce the \ROOTSGD algorithm that is the focus of our study.
We first motivate the algorithm from an averaging and variance reduction perspective.  
We then describe the burn-in and restarting mechanism, which contributes to the superior theoretical guarantees in the overall algorithm.

\pb\subsection{Motivation and gradient estimator}\label{sec:propose_rootsgd} 
%The algorithm is closely related to existing first-order algorithms, but differs critically in its choice of step-sizes.  
Our choice of step-size emerges from an overarching statistical perspective---rather than viewing the problem as one of correcting SGD via particular mechanisms such as averaging, variance reduction or momentum, we instead view the problem as one of utilizing all previous online data samples, $\xi_1,\dots,\xi_t\sim P$, to form an \emph{estimate} $\mathsf{Estimator}_t$ of $\nabla F(\theta_{t-1})$ at each round $t$.  
We then perform a gradient step based on this estimator---that is, we compute $\theta_t = \theta_{t-1} - \stpsz_t\cdot \mathsf{Estimator}_t$.

Concretely, our point of departure is the following \emph{idealized} estimate of the error in the current gradient:
\begin{align}\label{onlinemean}
\mathsf{Estimator}_t - \nabla F(\theta_{t-1}) 
= 
\frac{1}{t} \sum_{s=1}^t (\nabla f(\theta_{s-1};\xi_s) - \nabla F(\theta_{s-1})) 
.
\end{align}
Treating the terms $\nabla f(\theta_{s-1};\xi_s) - \nabla F(\theta_{s-1}), s=1,\dots,t$ as martingale differences, and assuming that the conditional variances of these terms are identical almost surely, it is straightforward to verify that the choice of equal weights $\frac{1}{t}$ minimizes the variance of the estimator over all such convex combinations.  
This simple but very specific choice of weights is central to our algorithm, which we refer to as \emph{Recursive One-Over-T SGD} (\ROOTSGD).

The recursive aspect of the algorithm arises as follows.  
We set $\mathsf{Estimator}_1 = \nabla f(\theta_0; \xi_1)$ and express \eqref{onlinemean} as follows:
\begin{align*}
\mathsf{Estimator}_t - \nabla F(\theta_{t-1}) 
= 
\frac{1}{t} (\nabla f(\theta_{t-1};\xi_t) - \nabla F(\theta_{t-1})) + \frac{t-1}{t} (\mathsf{Estimator}_{t-1} - \nabla F(\theta_{t-2}))
.
\end{align*}
Rearranging gives
\begin{align*}
\mathsf{Estimator}_t 
= 
\frac{1}{t} \nabla f(\theta_{t-1};\xi_t) + \frac{t-1}{t} \left(\nabla F(\theta_{t-1}) - \nabla F(\theta_{t-2})\right) + \frac{t-1}{t} \mathsf{Estimator}_{t-1}
.
\end{align*}
We now note that we do \emph{not} generally have access to the bracketed term $\nabla F(\theta_{t-1}) - \nabla F(\theta_{t-2})$, and replace the term by an unbiased estimator, $\nabla f(\theta_{t-1};\xi_t) - \nabla f(\theta_{t-2};\xi_t)$, based on the current sample $\xi_t$. 
Intuitively, the replacement should not affect much as long as the stochastic function admits some smoothness condition.
Letting $v_t$ denote $\mathsf{Estimator}_t$ with this replacement, we obtain the following recursive update:
\begin{align}
v_t 
&= 
\frac{1}{t} \nabla f(\theta_{t-1};\xi_t) + \frac{t-1}{t} \left(\nabla f(\theta_{t-1};\xi_t) - \nabla f(\theta_{t-2};\xi_t) \right) + \frac{t-1}{t} v_{t-1} \notag
\\&= 
\underbrace{\nabla f (\theta_{t-1}; \xi_t)}_{\text{stochastic gradient}} + \underbrace{
\frac{t-1}{t} \left(v_{t-1} - \nabla f (\theta_{t-2}; \xi_t)\right)
}_{\text{correction term}}
, 
\label{vupdate_new}
\end{align}
consisting of both a stochastic gradient and a correction term.

Finally, performing a gradient step based on our estimator yields the \ROOTSGD algorithm:
\begin{subequations}\label{algoROOTSGD}
\begin{align}
v_t		&=	\nabla f (\theta_{t-1}; \xi_t) + \frac{t-1}{t} \left(v_{t-1} - \nabla f (\theta_{t-2}; \xi_t)\right)
\label{algoROOTSGD-v}
\\
\theta_t	&=	\theta_{t-1} - \stpsz_t v_t
\label{algoROOTSGD-theta}
,
\end{align}
\end{subequations}
where $\{\stpsz_t\}_{t \geq 1}$ is a suitably chosen sequence of positive step-sizes.  
Note that $v_t$ defined in Eq.~\eqref{vupdate_new} is a recursive estimate of $\nabla F(\theta_{t-1})$ that is \emph{unconditionally} unbiased in the sense that $\Exs[ v_t ] = \Exs[ \nabla F(\theta_{t-1}) ]$. 
So the $\theta$-update is an approximate gradient-descent step that moves
along the negative direction $-v_t$.%
\footnote{Unlike many classical treatments of stochastic approximation, we structure the subscripts so they match up with those of the filtration corresponding to the stochastic processes.}%

We initialize $\theta_0 \in \real^d$, and, to avoid ambiguity, we define the update~\eqref{algoROOTSGD} at $t=1$ to use only $v_1 = \nabla f(\theta_0; \xi_1)$. 
Overall, given the initialization $(\theta_0,v_0,\theta_{-1}) = (\theta_0, 0, \text{arbitrary})$,
at each step $t \geq 1$ we take as input $\xi_t\sim P$, and perform an update of $(\theta_t, v_t, \theta_{t-1})$. 
This update depends only on $(\theta_{t-1}, v_{t-1}, \theta_{t-2})$ and $\xi_t$, and is first-order and Markovian.

\pb\subsection{Two-time-scale structure and burn-in period}\label{sec_hidden}
For the purposes of both intuition and the proof itself, it is useful to observe that the iterates~\eqref{algoROOTSGD} evolve in a two-time-scale manner.  
Define the process $z_t \mydefn v_t - \nabla F(\theta_{t - 1})$ for $t = 1, 2, \cdots$, which characterizes the \emph{tracking error} of $v_t$ as an estimator for the gradient.  
For each $\theta \in \real^d$ and $\xi\sim P$ we define the noise term $\noise(\theta; \xi) = \nabla_\theta f(\theta; \xi) - \nabla F(\theta)$, and use the shorthand notation $\noise_s(\cdot) = \noise(\cdot; \xi_s)$.  
Some algebra yields the decomposition
\begin{subequations}
\begin{align}\label{eq:zt-first-appear}
t \cdot z_t 
= 
\sum_{s = 1}^{t} \noise_s(\theta_{s - 1}) + \sum_{s = 1}^{t} (s - 1) \big( \noise_s(\theta_{s - 1}) - \noise_s(\theta_{s - 2}) \big)
, \quad \mbox{valid for $t = 1, 2, \ldots.$}
\end{align}
In this way, we see that the process $(t \cdot z_t)_{t \geq 1}$ is a martingale difference sequence adapted to the natural filtration $(\filtration_t)_{t \geq 0}$.  Indeed, the quantity $z_t$ plays the role of averaging the noise as well as performing a weighted averaging of consecutive differences collected along the path.  
On the other hand, the process $(t \cdot v_t)_{t \geq 1}$ moves rapidly driven by the strong convexity of the function $F$:
\begin{align}\label{eq:vt-first-appear}
t v_t 
= 
(t - 1) \big\{ v_{t - 1} + \nabla F (\theta_{t - 1})- \nabla F (\theta_{t - 2}) \big\} 
+ 
\nabla F(\theta_{t - 1}) 
+
\noise_t (\theta_{t - 1}) 
+ 
(t - 1) \big( \noise_t (\theta_{t - 1}) - \noise_t (\theta_{t - 2}) \big)
.
\end{align}
\end{subequations}
Given an appropriate step-size $\stpsz_t$, the first term on the RHS of Eq.~\eqref{eq:vt-first-appear} exhibits a contractive behavior.
Consequently, the process $(t v_t)_{t \geq 1}$ plays the role of a \emph{fast process}, driving the motion of iterates $(\theta_t)_{t \geq 0}$, and the noise-collecting process $z_t$ is a \emph{slow process}, collecting the noise along the path and contributing to the asymptotic efficiency of $\theta_t$. 
Note that the fast process moves with a step-size $\stpsz_t$, making $\stpsz_t \strongconvex$ progress when $F$ is $\strongconvex$-strongly convex, while the slow process works with a step-size $\frac{1}{t}$. 
In order to make the iterates stable, we need the fast process to be fast in a relative sense, requiring that $\stpsz_t \strongconvex \geq \frac{1}{t}$. 
This motivates a burn-in period in the algorithm, namely, in the first $\burnin$ iterations, we run the recursion~\eqref{algoROOTSGD} with step-size zero and simply average the noise at $\theta_0$; we then start the algorithm with an appropriate choice of step-size. 
Concretely, given some initial vector $\theta_0 \in \real^d$, we set $\theta_t = \theta_0$ for all $t = 1, \ldots, \burnin-1$, and compute
\begin{align}
v_t	=	\frac{1}{t} \sum_{s = 1}^t \nabla f(\theta_0; \xi_s) 
,\qquad
\mbox{for all $t = 1, \ldots, \burnin$.}
\end{align}

As suggested by our discussion, an algorithm with step-size $\stpsz_t = \stpsz$ will need a burn-in period of length $\burnin \asymp \frac{1}{\stpsz\strongconvex}$ for a $\strongconvex$-strongly convex function $F$. 
%
%
%\iffalse
Equivalently, we can view the step-sizes in the update for $\theta_t$ as being scheduled as follows:
\begin{align}\label{step-sizet}
\stpsz_t =
\begin{cases}
\stpsz,	& \text{for } t\ge \Tburnin
, \\ 
0,	& \text{for } t=1,\dots,\Tburnin-1
,
\end{cases}
\end{align}
briefed as $\stpsz_t = \stpsz\cdot 1[t\ge \burnin]$, and, accordingly, the update rule from Eqs.~\eqref{algoROOTSGD-v} and \eqref{algoROOTSGD-theta} splits into two phases:
\begin{subequations}
\begin{align*}
v_t
& =
\begin{cases}
\nabla f (\theta_{t-1}; \xi_t) + \frac{t-1}{t} \left(v_{t-1} - \nabla f (\theta_{t-2}; \xi_t)\right)	& \text{for } t \ge \Tburnin + 1 
, \\
\frac{1}{t} \displaystyle\sum_{s = 1}^t \nabla f (\theta_0; \xi_s)						& \text{for } t = 1, \dots, \Tburnin 
,
\end{cases}
\\ 
\theta_t 
&=
\begin{cases}
\theta_{t-1} - \stpsz v_t		&\text{for } t \ge \Tburnin
 , \\ 
\theta_{0}				&\text{for } t = 1, \dots, \Tburnin-1 .
\end{cases} 
\end{align*}
\end{subequations}
Such an algorithmic design has the length of the burn-in period for our algorithm is identical to the number of processed samples, so it features that the iteration number is identical to the sample complexity.
%\fi
The \ROOTSGD scheme is presented formally as Algorithm~\ref{algo_singleepoch}; for the remainder of this paper, when referring to \ROOTSGD, we mean Algorithm~\ref{algo_singleepoch} unless specified otherwise.

\begin{algorithm}[!tb]
\caption{\ROOTSGD}%
\begin{algorithmic}[1]
\STATE
\textbf{Input:}
initialization $\theta_0$; step-size sequence $(\stpsz_t)_{t\ge 1}$
\FOR{$t = 1, 2, \dots, T$} 
\STATE
$ v_t = \nabla f(\theta_{t-1}; \xi_t) + \frac{t-1}{t} \left(v_{t-1} - \nabla f(\theta_{t-2}; \xi_t)\right) $ 
\STATE
$ \theta_t = \theta_{t-1} - \stpsz_t v_t $ 
\ENDFOR 
\STATE\textbf{Output:}
$\theta_{\niters}$
\end{algorithmic}
\label{algo_singleepoch}
\end{algorithm}

\pb\section{Main results}\label{sec_algoresult_hessiansmooth}
In this section, we present our main nonasymptotic and asymptotic results.  
We first establish a preliminary nonasymptotic result in \S\ref{sec_algoresult_nonasy}.  
With augmented smoothness and moment assumptions, we then introduce in \S\ref{sec_nonasymp_assumptions} sharp nonasymptotic upper bounds with unit pre-factor on the term characterizing the optimal statistical risk.  
Finally, in \S\ref{sec_algoresult_asy_hessiansmooth}, we establish the asymptotic efficiency of \ROOTSGD.

\pb\subsection{Preliminary nonasymptotic results}\label{sec_algoresult_nonasy}
We begin by presenting preliminary nonasymptotic results for \ROOTSGD.
Before formally presenting the result, we detail our assumptions for the stochastic function $f(\cdot;\xi)$ and the expectation $F$.

First, we impose strong convexity and smoothness assumptions on the objective function:

\begin{assumption}[Strong convexity and smoothness]\label{assu_StrcvxSmooth}
The population objective objective function $F$ is twice continuously differentiable, $\strongconvex$-strongly-convex and $\smoothness$-smooth for some $0<\strongconvex\le \smoothness<\infty$:
\begin{align*}
\vecnorm{\nabla F(\theta) - \nabla F(\theta')}{2}
\le
\smoothness\vecnorm{\theta - \theta'}{2}
,\qquad \mbox{and} \quad
\binprod{\nabla F(\theta) - \nabla F(\theta')}{ \theta - \theta'}
\ge
\strongconvex\vecnorm{\theta - \theta'}{2}^2
,
\end{align*}
for all pairs $\theta, \theta' \in \real^d$.
\end{assumption}

Second, we assume sufficient regularity for the covariance matrix at the global minimizer $\thetastar$:

\begin{assumption}[Finite variance at optimality]\label{assu_noisethetastar}
At any minimizer $\thetastar$ of $F$, the stochastic gradient $\nabla f(\thetastar; \xi)$ has a positive definite covariance matrix, $
\SigStar	\defn		\Exs\left[ \nabla f(\thetastar; \xi) (\nabla f(\thetastar; \xi))^\top \right]
$, with its trace $
\sigstar^2	\defn		\Exs\vecnorm{\nabla f(\thetastar; \xi)}{2}^{2}
$ assumed to be finite.
\end{assumption}
\noindent
Note that we only assume a finite variance on the stochastic gradient at the global minimizer $\thetastar$.  
This is significantly weaker than the standard assumption of a globally bounded noise variance.  
See \cite{nguyen2019new} and \cite{lei2020adaptivity} for a detailed discussion of this assumption on the noise.

Third, we impose a mean-squared Lipschitz condition on the stochastic noise:

\begin{assumption}[Lipschitz stochastic noise]\label{assu_smoothnoise}
The noise function $\theta \mapsto \noise(\theta; \xi)$ in the associated stochastic gradients satisfies the bound
\begin{align}\label{smoothnoise}
\Exs\vecnorm{\noise(\theta; \xi) - \noise(\theta'; \xi)}{2}^2
\le
\sglipsq \vecnorm{\theta - \theta'}{2}^2 
,\qquad 
\text{for all pairs $\theta; \theta' \in \real^d$.}
\end{align}
\end{assumption}

\noindent
We note that in making Assumption~\ref{assu_smoothnoise}, we separate the stochastic smoothness (in the $L^2$ sense) of the noise, $\noise (\theta; \xi) = \nabla f (\theta; \xi) - \nabla F (\theta)$, from the smoothness of the population-level objective. The magnitude of $\sglip$ and $\smoothness$ are not comparable in general. This flexibility permits, for example, mini-batch algorithms where the population-level Lipschitz constant $\smoothness$ remains fixed but the parameter $\sglip$ decreases with batch size.
Such a separation has been adopted in nonconvex stochastic optimization literature \cite{arjevani2020second}.%
\footnote{Observe that Assumptions~\ref{assu_StrcvxSmooth} and~\ref{assu_smoothnoise} imply a mean-squared Lipschitz condition on the stochastic gradient function:
\begin{align*}
\Exs\vecnorm{ \nabla f(\theta; \xi) - \nabla f(\theta'; \xi)}{2}^2 
&= 
\vecnorm{\nabla F(\theta) - \nabla F(\theta')}{2}^2 + \Exs\vecnorm{\noise(\theta; \xi) - \noise(\theta'; \xi)}{2}^2
%\\& \le 
\le
\left(L^2 + \sglip^2 \right) \vecnorm{\theta - \theta'}{2}^2
,
\end{align*}
where the final step uses the $L$-Lipschitz condition on the population function $F$.
}%

Finally, we remark that all of these assumptions are standard in the stochastic
optimization and statistical literature; and specific instantations of these assumptions are satisfied by a broad class of statistical models and estimators.  
We should note, however, that the strong convexity and smoothness (Assumption~\ref{assu_StrcvxSmooth}) is a global condition stronger than those typically used in the asymptotic analysis of M-estimators in the statistical literature. 
These conditions are needed for the fast convergence of the algorithm as an optimization algorithm, making it possible to establish nonasymptotic bounds.
Assumptions~\ref{assu_noisethetastar} and~\ref{assu_smoothnoise} are standard for proving asymptotic normality of M-estimators and Z-estimators (e.g.~\cite{VAN[derVaart]}, Theorem 5.21).
%They are satisfied by a broad class of statistical models and estimat
In contrast to some prior work (e.g.~\cite{ghadimi2012optimal,ghadimi2013optimal}), we \emph{do not} assume uniform upper bounds on the variance of the stochastic gradient noise; this assumption fails to hold for various statistical models of interest, and theoretical results that dispense with it are of practical interest.

With the aforementioned assumptions in place, we provide our first preliminary nonasymptotic result for single-epoch \ROOTSGD, as follows:

\begin{theorem}[Preliminary nonasymptotic results, single-epoch \ROOTSGD]\label{theo_finitebdd_single}
Under Assumptions~\ref{assu_StrcvxSmooth},~\ref{assu_noisethetastar},~\ref{assu_smoothnoise}, suppose that we run Algorithm~\ref{algo_singleepoch} with
burn-in period $\Tburnin$ and step-size $\stpsz$ such that
\begin{align}\label{eq:step-size-and-burnin}
\Tburnin \mydefn \frac{24}{\stpsz\strongconvex}
\qquad \mbox{and}\quad 
\stpsz \in (0, \stpszmax]
, \qquad 
\mbox{where}\quad
\stpszmax \defn \frac{1}{4L} \land \frac{\strongconvex}{8\sglip^2}
.
\end{align}
Then, for any iteration $\niters \ge 1$, the iterate $\theta_{\niters}$ satisfies the bound
\begin{align}\label{finitebdd_single}
\Exs\| \nabla F(\theta_{\niters}) \|_2^2 
&\le 
\frac{28 \; \sigstar^2}{\niters} 
+ 
\frac{2700 \; \vecnorm{\nabla F(\theta_0)}{2}^2}{\stpsz^2 \strongconvex^2 \niters^2} 
.
\end{align}
\end{theorem}
\noindent
We provide a complete analysis of Theorem \ref{theo_finitebdd_single} in \S\ref{sec_proof,theo_finitebdd_single}. 
In order to interpret the result, we make few remarks in order.

\begin{enumerate}[label=(\roman*)]
\item
Note when stating the upper bound~\eqref{finitebdd_single} we adopt the convergence metrics in expected squared gradient norm.
The guarantee in~\eqref{finitebdd_single} consists of the sum of two terms which differs in magnitude as $\niters\to\infty$.
The leading-order first term is contributed by the \emph{optimal statistical risk}, and is determined by the noise variance $\sigstar^2$ at the minimizer.  
The higher-order second term, on the other hand, exhibits an $O(\frac{1}{\niters^2})$-dependency on the initial condition, which is suboptimal and can be improved to an exponential dependency by properly restarting the algorithm.
When the step-size $\stpsz$ is fixed, a comparison of the two summand terms \eqref{finitebdd_single} yields that the optimal asymptotic risk $\frac{\sigstar^2}{\niters}$ for the squared gradient holds up to an absolute constant whenever $\niters \gtrsim \frac{1}{\stpsz\mu} \lor \frac{\vecnorm{\nabla F(\theta_0)}{2}^2}{\stpsz^2 \strongconvex^2 \sigstar^2}$.

\item
Suppose that we use the maximal step-size $\stpszmax$ permitted by condition~\eqref{eq:step-size-and-burnin}.
By converting the convergence rate bound~\eqref{finitebdd_single} into a sample complexity bound, we then find that it suffices to take
\iffalse
\footnote{
Indeed, we choose $\niters$ in Eq.~\eqref{finitebdd_single} to be sufficiently large such that it satisfies the inequalities $
T\ge	\Tburnin	=	\lceil \frac{24}{\eta\mu }\rceil
$, as well as $
\frac{2700 \vecnorm{\nabla F(\theta_0)}{2}^2}{\eta^2 \strongconvex^2 \niters^2}	\le	\frac{\varepsilon^2}{2} 
$ and $
\frac{28\sigstar^2}{\niters}			\le	\frac{\varepsilon^2}{2}
$.
Here and on, we assume without loss of generality that $\varepsilon^2 \le \vecnorm{\nabla F(\theta_0)}{2}^2$.
It is then straightforward to see that \eqref{complexity_finitebdd_single} serves as a tight sample complexity upper bound.
}
\fi
\begin{align}
C_{\ref{theo_finitebdd_single}}(\varepsilon)
	&=
\max\left\{
\frac{74}{\eta_{\max}\mu}\cdot\frac{\vecnorm{\nabla F(\theta_0)}{2}}{\varepsilon}
, 
\frac{56\sigstar^2}{\varepsilon^2} 
\right\}
% \notag	\\&
	\asymp
\max\left\{
\left(\frac{\smoothness}{\strongconvex} + \frac{\sglip^2}{\strongconvex^2} \right) \cdot \frac{\vecnorm{\nabla F(\theta_0)}{2}}{\varepsilon}
,
\frac{\sigstar^2}{\varepsilon^2}
\right\}
 \label{complexity_finitebdd_single}
\end{align}
samples in order to obtain an estimate of $\thetastar$ with gradient norm bounded as $O(\varepsilon)$.
When the asymptotics holds as $\varepsilon$ tends to zero with other problem-dependent constants being bounded away from zero, the leading-order term $\asymp\frac{\sigstar^2}{\varepsilon^2}$ in $C_{\ref{theo_finitebdd_single}}(\varepsilon)$ matches the optimal statistical risk up to a constant pre-factor.
To the best of our knowledge, such a bound on sample complexity is achieved for the first time by a stochastic first-order algorithm in the setting where only first-order smoothness condition holds, i.e.~no continuity condition on the Hessians are posed.
The only prior results which reported a near-optimal statistical risk under comparable settings in the leading-order stochastic optimization are due to~\cite{nguyen2021inexact} and \cite{allen2018make}, achieving the optimal risk up to a polylogarithmic factor.
In Appendix~\S\ref{sec_discussion}, we compare our non-asymptotic results with these existing works in detail.

\end{enumerate}

\pb\subsection{Improved nonasymptotic upper bounds}\label{sec_nonasymp_assumptions} 
The convergence rate bound of Theorem~\ref{theo_finitebdd_single} matches the optimal risk by a constant pre-factor $c$---to be precise, $c=28$ in the provided analysis.  
In addition to this non-optimal pre-factor, this result does not match the efficiency of M-estimators in its higher-order dependency.
So as to overcome these limitations, we now show how to apply Theorem~\ref{theo_finitebdd_single} as the building block to seek to obtain a sharp fine-grained convergence rate via two-time-scale characterization, under additional smoothness and moment assumptions.

First, we need the following \emph{Lipschitz continuity condition} for the Hessian at the optimum. 
We denote $\hessianstar \mydefn \nabla^2 F (\thetastar)$ throughout.
\begin{assumption}[Lipschitz continuous Hessians]\label{assu_holderhessian}
There exists a Lipschitz constant $\holderconst>0$ such that
\begin{align}\label{holder_continuity}
\opnorm{\nbF{\theta} - \nbF{\thetastar}}
\le
\holderconst \norm{\theta - \thetastar }
.
\end{align}
\end{assumption}

We also need fourth-moment analogue of Assumptions~\ref{assu_noisethetastar} and~\ref{assu_smoothnoise}, stated as follows. 
Note that these conditions are also exploited in prior work on nonasymptotic analyses of PRJ averaging procedure~\cite{bach2011non,xu2011towards,gadat2017optimal} and Streaming SVRG~\cite{frostig2015competing}.

\begin{assumption}[Finite fourth moment at minimizers]\label{assu_noise_bdd_holder}
Let Assumption~\ref{assu_noisethetastar} hold, and let $\sigstartild^2 \defn\sqrt{\Exs\vecnorm{\nabla f(\thetastar; \xi)}{2}^4}$ be finite.
\end{assumption}
\noindent
Observe that $\sigstar \leq \sigstartild$ by H\"{o}lder's inequality. 
This distinction is important, as $\sigstar^2$ corresponds to the optimal statistical risk (measured in gradient norm), while $\sigstartild^2$ does not.

For the higher-order moments of stochastic gradients, we introduce the following

\begin{assumption}[Lipschitz stochastic noise in fourth moment]\label{assu_noise_smooth_holder}
The noise function $\theta \mapsto \noise(\theta; \xi)$ in the associated stochastic gradients satisfies the bound
\begin{align}\label{noise_smooth_holder}
\sqrt{\Exs\vecnorm{\noise (\theta_1; \xi) - \noise (\theta_2; \xi)}{2}^4}
\le
\sglip^2 \vecnorm{\theta_1 - \theta_2}{2}^2
,\qquad
\text{for all pairs $\theta_1, \theta_2 \in \real^d$.}
\end{align}
\end{assumption}
\noindent
Note that we slightly abuse the notation and denote $\sglip$ by both moment Lipschitz constants in Assumptions~\ref{assu_smoothnoise} and \ref{assu_noise_smooth_holder}. 
In the presentations for the rest of this subsection, the notation $\sglip$ should be understood as the parameter in Assumption~\ref{assu_noise_smooth_holder}, which is strictly stronger than Assumption~\ref{assu_smoothnoise}.

Formally, we present a multi-epoch version of the \ROOTSGD algorithm in Algorithm~\ref{algo_multiepoch}. 
The algorithm runs $\epochsrestart$ short epochs and one long epoch. 
The goal of each short epoch is to ``halve’’ the dependency on the initial condition $\vecnorm{\nabla F(\theta_0)}{2}$, and it suffices to take $\shortepoch = c \burnin$ for some universal constant $c > 1$. 
We further impose the mild condition that the quantity $\vecnorm{\nabla F(\theta_0)}{2} / \sigstar $ scales as a polynomial function of $\numobs$.%
\footnote{This assumption is used only to simplify the presentation. If it does not hold true, the $\log \numobs$ terms in the bounds will be replaced by $\log \numobs + \log \big(1 +  \vecnorm{\nabla F(\theta_0)}{2} / \sigstar \big)$.\label{footnote:log-terms}}
In the following Theorem~\ref{theo_rootHolder_multi}, we present the gradient norm bounds satisfied by the multi-epoch \ROOTSGD algorithm:
\begin{algorithm}[!tb]
\caption{\ROOTSGD, multi-epoch version}
\begin{algorithmic}[1]
\STATE
\textbf{Input:}
initialization $\theta_0$; fixed step-size $\stpsz$; burn-in time $\burnin$; short epochs length $\inneriters\ge \burnin$; short epochs number $\Epoch$
\STATE
Set initialization for first epoch $\theta_0^{(1)} = \theta_0$
\FOR{$b = 1, 2, \cdots, \Epoch$}
\STATE
Run \ROOTSGD (Algorithm \ref{algo_singleepoch}) for $\inneriters$ iterates with burn-in time $\burnin$ (i.e.~step-size sequence $(\stpsz_t)_{t\ge 1}$ defined as in Eq.~\eqref{step-sizet})
\STATE 
Set the initialization $\theta_0^{(b + 1)} \mydefn \theta_{\inneriters}^{(b)}$ for the next epoch
\ENDFOR
\STATE
Run \ROOTSGD (Algorithm \ref{algo_singleepoch}) for $\niters	\mydefn	\nsamples - \inneriters\Epoch$ iterates with burn-in time $\burnin$
\STATE
\textbf{Output:}
The final iterate estimator $\theta^{\text{\rm final}}_\nsamples	\mydefn	\theta_\niters^{(\Epoch + 1)}$
\end{algorithmic}
\label{algo_multiepoch}
\end{algorithm}

\begin{theorem}[Improved nonasymptotic upper bound, multi-epoch \ROOTSGD]\label{theo_rootHolder_multi}
Under Assumptions \ref{assu_StrcvxSmooth}, \ref{assu_holderhessian}, \ref{assu_noise_bdd_holder}, \ref{assu_noise_smooth_holder}, suppose that we run Algorithm~\ref{algo_multiepoch} with the number of short epochs $\Epoch = \left\lceil \frac12\log\left(\frac{e\normb{\naF{\theta_0}}}{\stpsz\strongconvex\sigstar^2}\right) \right\rceil $, the burn-in time $ \burnin = \frac{24}{\stpsz\strongconvex} $, and the small epoch length $ \inneriters = \frac{7340}{\stpsz\strongconvex}$.  
Then for any step-size $\stpsz \in \big( 0, \frac{1}{56\smoothness} \wedge \frac{\strongconvex}{64\sgliptild^2} \big]$ and $\nsamples \ge \inneriters \Epoch + 1$, it returns an estimate $\theta^{\text{\rm final}}_\nsamples$ such that
\begin{align}\label{rootHolder_multi_raw}
\Exs\normb{\naF{\theta^{\text{\rm final}}_\nsamples}} -
\frac{\sigstar^2}{\niters} \le C \left\{
\frac{\sglip^2\stpsz}{\strongconvex} 
+ \frac{\log \niters}{\stpsz\strongconvex \niters} 
+ \frac{\sglip^2\log \niters}{\strongconvex^2\niters} 
\right\} \frac{\sigstar^2}{\niters} 
+ 
\frac{C \holderconst \sigstartild^3}{\stpsz^{1/2}\strongconvex^{5/2} \niters^{2}}
\end{align}
where $\niters := \nsamples - \inneriters \Epoch$, and $C$ is a universal constant.
\end{theorem}
\noindent See \S\ref{sec_proof,theo_rootHolder_multi} for the proof of this theorem.
%We observe that \eqref{rootHolder_multi_raw} essentially keeps the terms in the first line of \eqref{rootHolder_single} via restarting, where the logarithmic dependency on initialization is incorporated in the burn-in time.

In order to interpret this result, let us take $\numobs \geq 2 \inneriters \Epoch$ so that we have $\tfrac{1}{T} \leq \frac{1}{\numobs} + \frac{2 \inneriters \Epoch}{\numobs^2}$. 
When the number of online samples $\nsamples$ is given a priori and the (constant) step-size is optimized as $\stpsz = \frac{c}{\sglip\sqrt{\nsamples}} \wedge \frac{1}{4\smoothness}$ where $c=0.49$, some algebra reduces the bound to
\begin{align}\label{rootHolder_multi_accu}
\Exs\normb{\naF{\theta^{\text{\rm final}}_\nsamples}} - \frac{\sigstar^2}{\nsamples} 
\lesssim
\left\{
\frac{\sglip}{\strongconvex\sqrt{\nsamples}}
+
\frac{\smoothness}{\strongconvex \nsamples} 
\right\}
\frac{\sigstar^2 \log \numobs}{\nsamples} 
+
\underbrace{
\left\{ 
\frac{\sgliptild}{\strongconvex \sqrt{\numobs}}
+ 
\frac{\smoothness}{\strongconvex \numobs} 
\right\}^{1/2} \frac{\holderconst}{\strongconvex^2} \bigg( 
\frac{\sigstartild^2}{\numobs} 
\bigg)^{3/2}
}_{=:\highordertilde_\numobs }
,
\end{align}
where $\highordertilde_\numobs$ is the linearization-induced term.
Given the sufficiently large sample size $\numobs$ satisfying the requirement $\numobs \gtrsim \frac{\smoothness}{\strongconvex} + \frac{\sgliptild^2}{\strongconvex^2}$, the pre-factors in the second term of~\eqref{rootHolder_multi_accu}, as well as the linearization error term $\tilde{\mathcal{H}}_\numobs$, start to diminish when the sample size $\numobs$ grows. The gradient norm bound~\eqref{rootHolder_multi_accu} consists of three terms. We discuss each of them as follows:

\begin{enumerate}[label=(\roman*)]
\item
The leading-order term $\frac{\sigstar^2}{\numobs}$ is exactly the asymptotic risk of the optimal limiting Gaussian random vector in the local asymptotic minimax theorem, measured with squared gradient norm. 
Note that this term depends only on the noise at the optimum $\thetastar$, instead of some uniform upper bounds.

\item 
The first term on the right-hand side consists of two parts that decay at different rates. 
If the sample size satisfies $\frac{\numobs}{\log \numobs} \gtrsim \frac{\smoothness}{\strongconvex} + \frac{\sgliptild^2}{\strongconvex^2}$ (which is a slightly stronger condition than the sample size needed for the theorem to hold true), this term is always smaller than the $\frac{\sigstar^2}{\numobs}$ term in the left hand side. 
For a large sample size $\numobs$, the dominating high-order term decays at the rate $O (\numobs^{-3/2} \log \numobs)$.

\item
The remaining high-order term $\highordertilde_\numobs$ in the bound~\eqref{rootHolder_multi_accu} scales as $O (\numobs^{-7/4})$, a faster rate of decay than the previous term. This term is induced by a linearization argument in our proof, and therefore depends on the Lipschitz constant $\smoothness_1$ of the Hessian matrix.

\end{enumerate}
Additionally, we remark that the sample size requirement $\numobs \gtrsim \frac{\smoothness}{\strongconvex} + \frac{\sgliptild^2}{\strongconvex^2}$ in Theorem~\ref{theo_rootHolder_multi} is natural: 
on the one hand, under the noise assumption~\ref{assu_smoothnoise}, it requires $\Theta \big( \frac{\sglip^2}{\strongconvex^2} \big)$ samples to distinguish the function $x \mapsto \frac{\strongconvex}{2} \vecnorm{x}{2}^2$ from a constant function $0$; 
on the other hand, the $\frac{\smoothness}{\strongconvex}$ term is consistent with the optimization essence of the problem --- as the \ROOTSGD algorithm reduces to gradient descent in the noiseless case, we need to pay for the complexity of gradient descent to achieve any meaningful guarantees.

%\item
Besides the expected gradient norm squared metric, we also establish guarantees in alternative metrics including the estimation error $\vecnorm{\theta_\numobs - \thetastar}{2}$ and the objective gap $F (\theta_\numobs) - F (\thetastar)$ with expectation taken.
In order to state the theorem, we define the following linearization-induced error terms that appears in the bound

\begin{subequations}
\begin{align}
\highordertilde_\numobs^{\mathrm{(MSE)}} 
&\mydefn 
\frac{c}{\lammin (\hessianstar)^2} \cdot \left\{ \frac{\holderconstprime}{  \strongconvex^2} \cdot \bigg( \frac{\sigstartild^2}{\numobs} \bigg)^{3/2} + \frac{\holderconstprime^2}{  \strongconvex^4} \cdot \bigg( \frac{\sigstartild^2}{\numobs} \bigg)^2 \right\}
,& \mbox{and}\\
\highordertilde_\numobs^{\mathrm{(OBJ)}} 
&\mydefn 
\frac{c}{\strongconvex} \cdot \frac{\holderconstprime}{  \strongconvex^2} \cdot \bigg( \frac{\sigstartild^2}{\numobs} \bigg)^{3/2} + \frac{c}{\lammin (\hessianstar)} \cdot \frac{\holderconstprime^2}{ \strongconvex^4} \cdot \bigg( \frac{\sigstartild^2}{\numobs} \bigg)^2
.&
\end{align}
\end{subequations}
For simplicity we only consider the multi-epoch \ROOTSGD as specified in Theorem~\ref{theo_rootHolder_multi}, where we conclude the following Corollary:

\begin{corollaryprime}[Nonasymptotic bounds in alternative metrics, multi-epoch \ROOTSGD]\label{theo_other_bounds}
Under the setup of Theorem~\ref{theo_rootHolder_multi}, the multi-epoch \ROOTSGD algorithm with the optimal step-size choice of $\stpsz\asymp \frac{1}{\smoothness} \wedge \frac{1}{\sglip \sqrt{\numobs}}$ produces an estimator that satisfies the following bound for $\numobs \geq \shortepoch \Epoch + 1$:
\begin{subequations}
\begin{align}
\Exs\vecnorm{\theta_\numobs^{\text{\rm final}} - \thetastar}{2}^2
-
\frac{1}{\numobs} \mathrm{Tr} \left( (\hessianstar)^{-1} \SigStar (\hessianstar)^{-1} \right)
&\le
c\left\{\frac{\sglip}{\strongconvex\sqrt{\nsamples}}
+
\frac{\smoothness}{\strongconvex \nsamples} 
\right\} \frac{\sigstar^2 \log \numobs}{\lammin (\hessianstar)^2 \numobs} + \highordertilde_\numobs^{(\mathrm{MSE})}
,\label{eq:main-mse-bound}\\
\Exs \left[ F \big( \theta_\numobs^{\text{\rm final}} \big) - F (\thetastar)  \right]
-
\frac{1}{2\numobs} \mathrm{Tr} \left( (\hessianstar)^{-1} \SigStar \right)
&\le
c\left\{\frac{\sglip}{\strongconvex\sqrt{\nsamples}}
+
\frac{\smoothness}{\strongconvex \nsamples} 
\right\} \frac{\sigstar^2 \log \numobs}{\lammin (\hessianstar) \numobs} + \highordertilde_\numobs^{(\mathrm{OBJ})}
.\label{eq:main-val-bound}
\end{align}
\end{subequations}
\end{corollaryprime}
\noindent 
See \S\ref{subsec:proof-theo-other-bounds} for the proof of this result, where the key technical addition lies on the adoption of a generic matrix-induced bound on the stochastic processes.
Obviously, the optimal step-size choice in alternative metrics is same in magnitude as the one in squared gradient norm metric in Theorem~\ref{theo_rootHolder_multi}.
We can compare the bounds in Corollary~\ref{theo_other_bounds} with the gradient norm bound~\eqref{rootHolder_multi_accu} induced by Theorem~\ref{theo_rootHolder_multi}, discussed term-by-term as follows:

\begin{enumerate}[label=(\roman*)]
%\bcar
\item
The leading-order terms in the bound~\eqref{eq:main-mse-bound} and~\eqref{eq:main-val-bound}, specified in the subtracted second terms on the left hands, are both optimal in a local asymptotic minimax sense with unity pre-factor. 
In particular, they are exactly the asymptotic risk of the limiting Gaussian random variable $\mathcal{N} \left( \thetastar, \frac{1}{\numobs} (\hessianstar)^{-1} \SigStar (\hessianstar)^{-1} \right)$ in corresponding metrics.
We also note that in the special case of well-specified maximal-likelihood estimation, Fisher's identity $\hessianstar = \SigStar$ holds true, and the leading-order terms in Eq.~\eqref{eq:main-mse-bound} and~\eqref{eq:main-val-bound} become $\frac{1}{\numobs} \mathrm{Tr} \left( (\hessianstar)^{-1} \right) $ and $\frac{d}{2n}$, respectively. 
This is in accordance with the classical asymptotic theory for M-estimators (c.f.~\cite{VAN[derVaart]}. \S5.3);

\item
The dominating term $\left\{\frac{\sglip}{\strongconvex\sqrt{\nsamples}}
+
\frac{\smoothness}{\strongconvex \nsamples} 
\right\}
\frac{\sigstar^2 \log \numobs}{\nsamples}$ in the gradient norm bound is correspondingly multiplied by a factor of $\lammin (\hessianstar)^{-2}$ (resp.~$\lammin (\hessianstar)^{-1}$) in the MSE (resp.~objective gap) bound on the right hand \eqref{eq:main-mse-bound} (resp.~\eqref{eq:main-val-bound}), which is intuitively consistent with conversions in metric;

\item
While the dominating high-order term on the right hand matches the corresponding optimal statistical risk and the higher-order terms altogether scale as $O(\numobs^{- 3/2})$, the linearization-induced error terms $\highordertilde_\numobs^{\mathrm{MSE}}$ and $\highordertilde_\numobs^{\mathrm{OBJ}}$ both decay at a rate of $O \big( \numobs^{-3/2} \big)$ as long as $\holderconstprime$ is bounded away from zero i.e.~the objective is essentially nonquadratic.
This is vastly different from the bound in gradient norm \eqref{rootHolder_multi_raw} of Theorem~\ref{theo_rootHolder_multi} where the linearization-induced terms are all incorporated in $O(\numobs^{-7/4})$, primarily due to that the pre-factor $
\left\{ \frac{\sgliptild}{\strongconvex \sqrt{\numobs}} + \frac{\smoothness}{\strongconvex \numobs} \right\}^{1/2}
$ in the linearization-induced term $\highordertilde_\numobs$ in \eqref{rootHolder_multi_accu} is replaced by unity.%
\footnote{This is because the Hessian-Lipschitz assumption plays a key role in relating MSE and objective gap to the underlying noise structure in the stochastic optimization problem, paying for larger linearization error; whereas in the gradient norm bound, the Hessian-Lipschitz assumption is employed only to mitigate the effect correlation that appears at even higher-order terms in the bound.}
In consistency with the metric conversion, the linearization-induced terms $\highordertilde_\numobs^{\mathrm{MSE}}$ and $\highordertilde_\numobs^{\mathrm{OBJ}}$ also incur additional factors related to the smallest eigenvalue of $\hessianstar$ on top of the linearization-induced error term $\highordertilde_\numobs$ in \eqref{rootHolder_multi_accu}.

\end{enumerate}
%\end{subtheorem}

\pb\subsection{Asymptotic results}\label{sec_algoresult_asy_hessiansmooth}
In this subsection, we study the asymptotic behavior of our \ROOTSGD algorithm.
We aim to prove the asymptotic efficiency of the multi-epoch estimator of Algorithm~\ref{algo_multiepoch} under minimal assumptions. 
In this case, Assumptions~\ref{assu_StrcvxSmooth},~\ref{assu_noisethetastar} and~\ref{assu_smoothnoise} are the standard ones needed for proving asymptotic normality of M-estimators and Z-estimators (see e.g.~\cite[Theorem 5.21]{VAN[derVaart]}). 
We first introduce our one-point Hessian continuity condition as follows as the qualitative counterpart of the continuity Assumption~\ref{assu_holderhessian}:

\begin{assumption}[One-point Hessian continuity]\label{assu_hessian_cts}
The Hessian mapping $\nabla^2 F(\theta)$ is continuous at the minimizer $\thetastar$, i.e.,
\begin{align*}
\lim_{\theta \rightarrow \thetastar} \opnorm{\nabla^2 F(\theta) - \hessianstar} = 0
.
\end{align*}
\end{assumption}
Note in Assumption~\ref{assu_hessian_cts} we assume only the continuity of Hessian matrix at $\thetastar$ without posing any bounds on its modulus of continuity. 
This is a much weaker condition than a Lipschitz or H\"{o}lder condition posted on the Hessian matrix as required in the analysis of the Polyak-Ruppert averaging procedure~\cite{polyak1992acceleration}.

With this setup, we are ready to state our weak convergence asymptotic efficiency result for $\theta^{\text{\rm final}}$ in the following theorem,%
\footnote{We emphasize our estimator's dependency on the step-size $\stpsz$ by explicitly bracketing it in the superscript.}
whose proof is provided in \S\ref{sec_proof,theo_asymptotic}:

\begin{theorem}[Asymptotic efficiency, multi-epoch \ROOTSGD]\label{theo_asymptotic_iteration}
Under
Assumptions~\ref{assu_StrcvxSmooth},~\ref{assu_noisethetastar},~\ref{assu_smoothnoise}
and~\ref{assu_hessian_cts}, suppose
that we run the multi-epoch Algorithm~\ref{algo_multiepoch} with
burn-in time $ \burnin = \frac{24}{\stpsz\strongconvex}$,
short-epoch length $ \inneriters = \frac{7340}{\stpsz\strongconvex}
$ and number of short epochs $ \Epoch = \left\lceil
\frac12\log\left(\frac{e\normb{\naF{\theta_0}}}{\stpsz\strongconvex\sigstar^2}\right)
\right\rceil $.  
Then as $\nsamples \to \infty$, $\stpsz\to 0$ such
that $\stpsz (\nsamples - \inneriters\Epoch)\to \infty$ and
$\inneriters\Epoch/\nsamples \to 0$, the estimate $\theta^{\text{\rm
    final},(\stpsz)}_\nsamples$ satisfies the weak convergence
\begin{align}
\label{asymptotic_iteration}
\sqrt{\nsamples} \left( \theta^{\text{\rm
    final},(\stpsz)}_\nsamples - \thetastar \right) \xrightarrow{d}
\mathcal{N} \left( 0, [\nabla^2 F(\thetastar)]^{-1}\SigStar[\nabla^2
  F(\thetastar)]^{-1} \right) ,
\end{align}
where $\SigStar\mydefn \Exs\left[\nabla f(\thetastar; \xi) \nabla f(\thetastar;\xi)^\top\right]$ is the covariance of the stochastic
gradient at the minimizer.
\end{theorem}

We remark that Theorem \ref{theo_asymptotic_iteration} holds under the mere additional assumption of one-point continuity on the Hessian matrix, which is usually the minimal assumption needed for an asymptotic efficiency result to hold.
Here we are adopting the multi-epoch \ROOTSGD with the same algorithmic specifications as in Theorem~\ref{theo_rootHolder_multi}, and we achieve the asymptotic convergence to the Gaussian limit that matches the Cram\'{e}r-Rao lower bound.  
The asymptotic covariance matrix in Eq.~\eqref{asymptotic_iteration}, however, carries significantly more information than the (scalar) optimal asymptotic risk.  
Our asymptotic result is in a triangular-array format: we let the fixed constant step-size scale down with $\nsamples$ where the scaling condition is essentially $\stpsz\to 0$, $\nsamples\to\infty$ with $\frac{\stpsz \nsamples}{\log(\stpsz^{-1})} \to \infty$, which is satisfied when $\stpsz \asymp \frac{1}{\nsamples^{c_1}}$ for any fixed $c_1\in (0,1)$.
Although not directly comparable, the range of step-size asymptotics is broader than \cite{polyak1992acceleration} and accordingly hints at potential advantages over PRJ, primarily due to our de-biasing corrections in our algorithm design and is consistent with our improved higher-order term in nonasymptotic result (Theorem~\ref{theo_rootHolder_multi} and Corollary~\ref{theo_other_bounds}).
In Appendix \S\ref{sec_algoresult_asy_hessiansmooth_additional}, we establish an additional asymptotic normality result for \ROOTSGD with fixed \emph{constant} step-size, which exhibits exactly the same limiting behavior as constant-step-size \emph{linear} stochastic approximation with PRJ averaging procedure under comparable asssumptions~\cite{mou2020linear}.

We end this subsection by remarking that Theorem~\ref{theo_asymptotic_iteration} only requires strong convexity, smoothness, and a set of noise moment assumptions standard in asymptotic statistics, but not any higher-order smoothness other than the continuity of Hessian matrices at $\thetastar$.
This matches the assumptions for asymptotic efficiency results in classical statistics literature \cite{VAN[derVaart],VAN-WELLNER}.

%%%%%%%%%%%%%%%%%%%%%%%%%%%%%%%%%%%%%%%%%%%%%%%%%%%%%%%%%%%%%%%%%%%%%%%%%%%%%%%%%%%%%%%%%%%%

\pb\section{Future directions}\label{sec_summaryfuture}
We have shown that \ROOTSGD enjoys favorable asymptotic and nonasymptotic behavior for solving the stochastic optimization problem \eqref{stoch_opt} in the smooth, strongly convex case.
With this result in hand, several promising future directions arise.
First, it is natural to extend the results for \ROOTSGD to non-strongly convex and nonconvex settings, for both nonasymptotic and asymptotic analyses.
Second, it would also be of significant interest to investigate both the nonasymptotic bounds and asymptotic efficiency of the variance-reduced estimator of \ROOTSGD in Nesterov's acceleration setting, in the hope of achieving all regime optimality in terms of the sample complexity to the stochastic first-order oracle.
Finally, for statistical inference using online samples, the near-unity nonasymptotic and asymptotic results presented in this work can potentially yield confidence intervals and other inferential assertions for the use of \ROOTSGD estimators.

\pb\section*{Acknowledgements}
We thank Peter Bartlett, Nicolas Flammarion, Koulik Khamaru, and Nicolas Le Roux for helpful discussions. 
This work was done in part while the authors were participating in the Theory of Reinforcement Learning program at the Simons Institute for the Theory of Computing.
This research was supported in part by Office of Naval Research Grant DOD ONR-N00014-18-1-2640 and N00014-21-1-2842, NSF-DMS grant 2015454 and 1612948, as well as NSF-IIS grant 1909365 to MJW, and also by the Mathematical Data Science program of the Office of Naval Research under grant number N00014-18-1-2764 and by the Vannevar Bush Faculty Fellowship program under grant number N00014-21-1-2941 and NSF grant IIS-1901252 to MIJ.
This research was also generously supported by NSF grant NSF-FODSI 202350 to MJW and MIJ.
%%%%%%%%%%%%%%%%%%%%%%%%%%%%%%%%%%%%%%%%%%%%%%%%%%%%%%%%%%%%%%%%%%%%%%

\pb
\bibliographystyle{alpha}
\bibliography{SAILreferences}

%%%%%%%%%%%%%%%%%%%%%%%%%%%%%%%%%%%%%%%%%%%%%%%%%%%%%%%%%%%%%%%%%%%%%%%%%

\newpage\appendix

\pb\section*{Appendices}
\noindent
In this appendix, we provide deferred proofs for theorems and lemmas in the main text organized as follows.
\S\ref{sec_morerelatedworks} provides additional work related to us.
\S\ref{sec_discussion} provides additional discussion on comparison of our results with concurrent work.
\S\ref{sec_general} proves the main results for both the nonasymptotic and the asymptotic convergence properties of our \ROOTSGD algorithm.
\S\ref{sec_algoresult_asy_hessiansmooth_additional} complements our asymptotic efficiency result in \S\ref{sec_algoresult_asy_hessiansmooth} and establishes an additional asymptotic result for constant-step-size \ROOTSGD.
\S\ref{sec:lem_holder} presents auxiliary lemmas stated in \S\ref{sec_proof,theo_rootHolder_single}, \S\ref{sec_proof,theo_rootHolder_multi} and \S\ref{subsec:proof-theo-other-bounds}.
Finally, \S\ref{SecAsympAuxiliary} proves necessary lemmas for the proof of Proposition \ref{theo_asymptotic_iteration_const_step_size} (in \S\ref{SecProofthm:asymptotic-const-step-size}).

\iffalse
\paragraph{Additional notation:}
Given vector-valued martingales $(X_t)_{t \geq \burnin}, (Y_t)_{t \geq \burnin}$ adapted to the filtration $(\filtration_t)_{t \geq \burnin}$, we use the following notation for cross variation for $t \geq \burnin$:
\begin{align*}
\crossvar{X}{Y}{t} \mydefn \sum_{s = \burnin + 1}^{t} \binprod{ X_t - X_{t - 1}}{Y_t - Y_{t - 1} }.
\end{align*}
We also define $\quadvar{X}{t} \mydefn \crossvar{X}{X}{t}$ to be the quadratic variation of the process $(X_t)_{t \geq \burnin}$.

We adopt some tensor notation.
For two matrices $A, B$, we use $A \otimes B$ to denote their Kronecker product. 
When it is clear from the context, we slightly overload the notation to let $A \otimes B$ denote the 4-th-order tensor produced by taking the tensor product of $A$ and $B$. 
Note that Kronecker product is just a flattened version of the tensor. 
For a $k$-th order tensor $T$, matrix $M$ and vector $v$, we use $T[M]$ to denote the $(k - 2)$-th order tensor obtained by applying $T$ to matrix $M$, and similarly, we use $T[v]$ to denote the $(k - 1)$-th order tensor obtained by applying $T$ to vector $v$.
\fi

\pb\section{Additional related work}\label{sec_morerelatedworks}
\paragraph{SGD and Polyak-Ruppert-Juditsky averaging procedure}
The theory of the stochastic approximation method has a long history since its birth in the 1950s \cite{robbins1951stochastic,bottou2004large,zhang2004solving,nemirovski2009robust,bottou2010large,bubeck2017convex} and recently regains its attention due to its superb performance in real-world application practices featured by deep learning \cite{GOODFELLOW-BENGIO-COURVILLE}, primarily due to its exceptional handling of the online samples.
Classics on this topic include
\cite{BERTSEKAS-TSITSIKLIS,BENVENISTE-METIVIER-PRIOURET,LJUNG-PFLUG-WALK,BORKAR} and many more.
Specially on the study of asymptotic normality which can trace back to \cite{fabian1968asymptotic}, the general idea of iteration averaging is based on the analysis of two-time-scale iteration techniques and it achieves asymptotic normality with an optimal covariance \cite{ruppert1988efficient,polyak1990new,polyak1992acceleration}.
Recent work along this line includes \cite{bach2011non,bach2013non,bach2014adaptivity,defossez2015averaged,flammarion2015averaging,dieuleveut2016nonparametric,duchi2021asymptotic,dieuleveut2017harder,allen2018make,dieuleveut2020bridging,asi2019stochastic}, presenting attractive asymptotic and nonasymptotic properties under a variety of settings and assumptions.
\cite{agarwal2012information,woodworth2016tight} provide minimax lower bounds for stochastic first-order algorithms.
\cite{jain2017parallelizing,jain2018markov,jain2018accelerating} analyze SGD and its acceleration with \emph{tail-averaging} that simultaneously achieves exponential forgetting and optimal statistical risk up to a constant, nonunity pre-factor.
It is also worth mentioning that iteration averaging provides robustness and adaptivity \cite{lei2020adaptivity}.
Instead of averaging the iterates, our \ROOTSGD algorithm averages the past stochastic gradients with proper de-biasing corrections and achieves competitive asymptotic performance.%
\footnote{A related but fundamentally different idea was proposed in \cite{nesterov2009primal,xiao2010dual,lee2012manifold} called \emph{dual averaging} for optimizing the regularized objectives. In contrast to their method, we focus in this work on the smooth objective setting and augment our estimator with de-biasing corrections. See also \cite{duchi2021asymptotic,tripuraneni2018averaging} for more on first-order optimization methods on Riemannian manifolds.}
For statistical inferential purposes, recent work \cite{chen2016statistical,su2018uncertainty} presents confidence interval assertions via online stochastic gradient with Polyak-Ruppert-Juditsky averaging procedure; analogous results for the \ROOTSGD algorithm are hence worth exploring, building upon the asymptotic normality that our work has established.

\paragraph{Variance-reduced gradient methods}
In the field of smooth and convex stochastic optimization, variance-reduced gradient methods represented by, but not limited to,
SAG~\cite{roux2012stochastic},
SDCA~\cite{shalev2013stochastic},
SVRG~\cite{johnson2013accelerating,konevcny2017semi,babanezhad2015stop,lei2017less},
SAGA~\cite{defazio2014saga},
SARAH~\cite{nguyen2017sarah,nguyen2021inexact}
have been proposed to improve the theoretical convergence rate of (stochastic) gradient descent.
Accelerated variants of SGD provide further improvements in convergence rate~\cite{lin2015universal,shalev2016sdca,allen2017katyusha,lan2018optimal,kulunchakov2019estimate,lan2019unified}. 
More recently, a line of work on recursive variance-reduced stochastic first-order algorithms have been studied in the nonconvex stochastic optimization literature \cite{fang2018spider,zhou2018stochastic,wang2019spiderboost,nguyen2021inexact,pham2020proxsarah,li2021page}. 
These algorithms, as well as their hybrid siblings~\cite{cutkosky2019momentum,tran2021hybrid}, achieve optimal iteration complexities for an appropriate class of nonconvex functions and in particular are faster than SGD under mild additional smoothness assumption on the stochastic gradients and Hessians \cite{arjevani2020second}.
Limited by space, we refer interested readers to a recent survey article by \cite{gower2020variance}, and while our \ROOTSGD algorithm can be viewed as a variant of variance-reduced algorithms, our goal is substantially different: we aim to establish for strongly convex objectives both a sharp, unity pre-factor nonasymptotic bound and asymptotic normality with Cram\'{e}r-Rao optimal asymptotic covariance that matches the local asymptotic minimax optimality \cite{duchi2021asymptotic}.

\paragraph{Sharp nonasymptotics and asymptotic efficiency}
When the objective admits additional smoothness, nonasymptotic rate analyses for either SGD with iteration averaging or variance-reduced stochastic first-order algorithms have been studied in various settings.
\cite{bach2011non} presents a nonasymptotic analysis of SGD with PRJ averaging procedure showing that, after processing $\numobs$ samples, the algorithm achieves a nonasymptotic rate that matches the Cram\'{e}r-Rao lower bound with a pre-factor equal to one with the additional term being $O(\numobs^{-7/6})$ (see the discussions in \S\ref{sec_discussion}).
\cite{xu2011towards,gadat2017optimal} improves the additional term to $O(\numobs^{-5/4})$ under comparable assumptions.
\cite{defossez2015averaged,dieuleveut2016nonparametric,duchi2021asymptotic,asi2019stochastic} achieves either sharp nonasymptotic bounds (in the quadratic case) or asymptotic efficiency that matches the local asymptotic minimax lower bound.
The asymptotic efficiency of variance-reduced stochastic approximation methods, however, has been less studied.
More related to this work, \cite{frostig2015competing} establishes the nonasymptotic upper bounds on the objective gap for an online variant of the SVRG algorithm \cite{johnson2013accelerating}, where the leading-order nonasymptotic bound on the excess risk matches the optimal asymptotic behavior of the empirical risk minimizer under certain \emph{self-concordant condition} posed on the objective function; the additional higher-order term reported is at least $\Omega(\numobs^{-8/7})$.

\paragraph{Other related work}
\cite{lakshminarayanan2018linear,mou2020linear} studies fixed-constant-step-size linear stochastic approximation with PRJ averaging procedure beyond an optimization algorithm \cite{polyak1992acceleration}, which includes many interesting applications in minimax game and reinforcement learning.
To be specific, \cite{lakshminarayanan2018linear} provides general nonasymptotic bounds which suffer from a nonunity pre-factor on the optimal statistical risk, and \cite{mou2020linear} studies the PRJ averaging procedure for general linear stochastic approximation and precisely characterizes the asymptotic limiting Gaussian distribution, delineating the additional term that adds onto the Cram\'{e}r-Rao asymptotic covariance and which vanishes as $\stpsz\to 0$ \cite{dieuleveut2020bridging}, and further establishes sharp concentration inequalities under stronger moment conditions on the noise.
\cite{arnold2019reducing} proposes an extrapolation-smoothing scheme of \emph{Implicit Gradient Transportation} to reduce the variance of the algorithm and provides convergence rates for quadratic objectives, which is further generalized to nonconvex optimization to improve the convergence rate of normalized SGD \cite{cutkosky2020momentum}.
For the policy evaluation problem in reinforcement learning, \cite{khamaru2020temporal} establishes an instance-dependent non-asymptotic upper bound on the $\ell_\infty$ estimation error, for a variance-reduced stochastic approximation algorithm. Their bound matches the risk of optimal Gaussian limit up to constant or logarithmic factors.
Recently, \cite{mou2022optimal} extends the algorithmic idea in this work and proposes the recursive variance-reduced stochastic approximation in span seminorm, which is applicable for generative models in reinforcement learning.

\pb\section{Comparison to related works}\label{sec_discussion}

In this section, we provide a careful comparison of our convergence results to those for stochastic first-order gradient algorithms.
For all nonasymptotic results, we compare our algorithm results with that of vanilla stochastic gradient descent, possibly equipped with iteration averaging and variance-reduced stochastic first-order optimization algorithms.
In the Lipschitz continuous Hessian case, we can achieve asymptotic unity.
We compare our \ROOTSGD convergence result with comparative work along with the following discussions in three aspects:

\paragraph{Comparison with classical results on SGD and its acceleration}
SGD is known to be worst-case optimal for optimizing smooth and strongly-convex objectives up to a constant pre-factor.
By way of contrast, our convergence metric in use, oracle query model and assumption on stochasticity are fundamentally different.
Despite this, a recent work due to~\cite{nguyen2019new} building upon earlier analysis surveyed by \cite{bottou2018optimization} makes a comparable noise assumption that allows the noise variance to grow at most quadratically with the distance to optimality and applies to SGD.
\cite{nguyen2019new} shows that for appropriate diminishing step-sizes $\stpsz_t$ we can conclude a guarantee of $
\Exs\|\theta^{\text{SGD}}_{\niters}-\thetastar\|_2^2
    \lesssim
\frac{\sigstar^2}{\strongconvex^2 \niters}
$.
With additional smoothness and noise assumptions we aim to achieve fine-grained non-asymptotic and asymptotic local asymptotic minimax optimality \emph{with unity pre-factor}.
It is straightforward to observe that the convergence rate bound of SGD under shared assumptions is in no regime better than that of \ROOTSGD presented in \eqref{complexity_finitebdd_single}.

We turn to compare our \ROOTSGD convergence result with the existing arts on stochastic accelerated gradient descent for strongly convex objectives~\cite{ghadimi2012optimal,ghadimi2013optimal}.
With an appropriate multi-epoch design, their guarantees on the objective gap are worst-case optimal for optimizing smooth and strongly-convex objectives in terms of the dependency on \emph{all} terms of condition number $\smoothness / \strongconvex$ and a uniform upper bound on the noise variance. 
Our preliminary nonasymptotic guarantee for single-epoch \ROOTSGD in Theorem \ref{theo_finitebdd_single}, in contrast, does not require uniform boundedness on the variance and depends solely on the variance at the minimizer $\thetastar$.%
\footnote{When measuring the risk via gradient norm, the optimal risk is characterized by the gradient noise variance $\sigstar$ at $\thetastar$.}
%%%%This is often inapplicable in the statistical learning setting (Nguyen et al., 2019), even though can be mitigated via a careful bootstrapping bound.
That being said, our result cannot not admit an accelerated rate in terms of condition number $\sqrt{\smoothness / \strongconvex}$ even with the help of multi-epoch designs. 
It is an important direction of future research to incorporate acceleration mechanism into our framework so as to achieve all-regime optimality.

\paragraph{Comparison with near-optimal guarantees in gradient norm}
\cite{allen2018make} develops a multi-epoch variant of SGD with averaging (under the name SGD3) via recursive regularization techniques and achieved a near-optimal rate for attaining an estimator of $O(\varepsilon)$-gradient norm. Our assumptions are not comparable in general: on the one hand, we assume a second-moment version of stochastic Lipschitz assumption (assumption~\ref{assu_smoothnoise}), which makes it possible to establish guarantees that depends on the noise variance $\sigstar$ at the optimum $\thetastar$; on the other hand,~\cite{allen2018make} makes no assumption on the modulus of continuity for the stochastic gradient, while their bound depends on a uniform upper bound $\sigma$ on the noise variance. Besides, it is worth noticing that as $\varepsilon \rightarrow 0^+$, the leading-order term in their bound scales as $\tfrac{\sigma^2}{\varepsilon^2} \log^3 \big( \tfrac{\smoothness}{\strongconvex} \big)$, which is sub-optimal by a polylogarithmic factor even if the uniform boundedness assumption on the variance is satisfied.
In a subsequent work, \cite{foster2019complexity} applies the idea of recursive regularization to AC-SA \cite{ghadimi2012optimal} and achieves an accelerated rate to find an approximate minimizer with $\varepsilon$-gradient norm, while a lower bound analysis provided further justifies the necessity of a multiplicative logarithmic factor in the nonstrongly convex, local oracle setting.

It is also worth-mentioning the (Inexact) SARAH algorithm and its analysis developed by~\cite{nguyen2021inexact}, which also achieves a near-optimal complexity upper bound of $O\left(\frac{\sigmastarsq}{\varepsilon^2} \log\left(\frac{1}{\varepsilon}\right) + \frac{\Lmax}{\strongconvex} \log\left(\frac{1}{\varepsilon}\right)\right)$ to obtain an approximate minimizer with $\varepsilon$-gradient norm.
Note that the setting for that result is slightly different (their setting is dubbed as the \ISCspace case in our extended analysis of Theorem \ref{theo_finitebdd_single_complete} in \S\ref{sec_proof,theo_finitebdd_single}).
The algorithm of~\cite{nguyen2021inexact} requires random output and burn-in batches that is inversely dependent on the desired accuracy $\varepsilon$, yielding a logarithmic pre-factor on top of the statistical error corresponding to the Cram\'{e}r-Rao lower bound; in comparison, our preliminary single-epoch \ROOTSGD result has a leading-order term in complexity bound that removes a logarithmic factor, attaining an optimal non-asymptotic guarantee up to a nonunity pre-factor.

\paragraph{Nonasymptotic guarantees matching local asymptotic minimax with near-unity pre-factor}%
\footnote{For convenience we include all comparable results in Table~\ref{tab_comparison_unity}.}
For SGD with PRJ averaging procedure, \cite{bach2011non} present a convergence rate that provides a useful point of comparison, although the assumptions are different (no Lipschitz gradient, bounded variance). In particular, when choosing the step-size $\stpsz_t = Ct^{-\bar{\alpha}}$ for $\bar{\alpha}\in (1/2,1)$, \cite{bach2011non} show that the following bound holds true for the averaged iterates $\bar{\theta}_{\niters}$ for the PRJ:
\begin{align*}
\sqrt{\Exs\vecnorm{\bar{\theta}_{\niters} - \thetastar}{2}^2} 
-
\sqrt{\frac{\mathrm{Tr}\left((\hessianstar)^{-1} \SigStar (\hessianstar)^{-1}\right)}{\niters}}
\le
\frac{c_0}{\niters^{2/3}}
,
\end{align*}
which corresponds to an $O(\nsamples^{-7/6})$ additional term in the squared estimation error metric (\eqref{eq:main-mse-bound} in Corollary~\ref{theo_other_bounds}).
Here, the constant $c_0$ depends on the initial distance to optimum, smoothness and strong convexity parameters of second- and third-order derivatives, as well as higher-order moments of the noise.
\cite{xu2011towards,gadat2017optimal} further improves the higher-order term from $O(\nsamples^{-7/6})$ to $O(\nsamples^{-5/4})$.
The convergence rate of (single-loop) \ROOTSGD is similar to SGD with PRJ averaging procedure in the nature of the leading term and the high-order terms, but our convergence rate bound of \ROOTSGD is comparatively cleaner and easier to interpret.

The work by \cite{frostig2015competing} proposes the Streaming SVRG algorithm that provides nonasymptotic guarantees in terms of the objective gap. 
Under a slightly different setting where smoothness and convexity assumptions are imposed on the individual function, their objective gap bound asymptotically matches the optimal risk achieved by the empirical risk minimizer under an additional self-concordance condition, with a multiplicative constant that can be made arbitrarily small. 
In particular, via our notations their results take the following form:
\begin{align*}
\Exs \big[ F(\thetahat_\numobs) - F (\thetastar) \big] 
\leq 
\Big(1 + \frac{5}{b} \Big)\frac{1}{2\numobs} \mathrm{Tr} \left( (\hessianstar)^{-1} \SigStar \right) + \mbox{high-order terms}
,
\end{align*}
where they require $\numobs \geq b^{2p + 3}$ for some $p \geq 2$. 
In order to achieve the sharp pre-factor, the additional term in this bound is at least $\Omega (\numobs^{- 8/7})$, a worse rate than our Corollary~\ref{theo_other_bounds}. 
Additionally, to get the corresponding nonasymptotic guarantees under such a setting, their bound requires a scaling condition $\niters \gtrsim \frac{\Lmaxsq}{\strongconvex^2}$ where $\Lmax$ denotes the smoothness of the individual function, which is larger than our burn-in sample size.
Without the self-concordance condition, the convergence rate bound of Streaming SVRG suffers from an extra multiplicative factor belonging to the interval $\left[1, \frac{\Lmax}{\strongconvex}\right]$, and its leading-order term thereby admits a dependency on the condition number worse than SGD.

% \iffalse
\begin{table}[!tb]
\centering
%\resizebox{\textwidth}{!}{%
\begin{tabular}{|c|c|c|l|}
\hline
Algorithm							&	Assumption					&	Additional Term											&	Reference
\\ \hline
PRJ								&	Hessian Lipschitz				&	$O\left(\frac{1}{\numobs^{7/6}}\right)$					&	\cite{bach2011non}
\\
PRJ								&	Hessian Lipschitz				&	$O\left(\frac{1}{\numobs^{5/4}}\right)$					&	\cite{xu2011towards,gadat2017optimal}
\\
Streaming SVRG					&	Self-concordant					&	multiplicative%
\footnote{Note that the paper~\cite{frostig2015competing} achieves a risk bound whose leading-order term is a $1 + O (b^{-1})$-multiplicative approximation to the optimal risk, with some additional terms (See Corollary 5 in their paper). Since this result requires $b^7\leq \numobs$, the additional term is at least $\Omega(\numobs^{-8/7})$.}	&	\cite{frostig2015competing}
\\
\ROOTSGD						&	Hessian Lipschitz				&	$O\left(\frac{1}{\numobs^{3/2}}\right)$					&	(This work)
%%%%%%%%%%%%%%%%%%%%%
\\ \hline
\end{tabular}
%}
\caption{Comparison of our results with comparative work. For the unity pre-factor nonasymptotic result, we only characterize the additional term to the optimal risk.}
\label{tab_comparison_unity}
\end{table}
% \fi

\pb\section{Proofs of nonasymptotic and asymptotic results}\label{sec_general}
We provide the convergence rate analysis and the proofs of our theorems in this section.
In our analysis we utilize the central object the \emph{tracking error process} $z_t$ defined as in \eqref{ztdefn}, and we heavily use the fact that the process $(t z_t)_{t \geq \burnin}$ is a martingale adapted to the natural filtration.

\pb\subsection{Proof of Theorem~\ref{theo_finitebdd_single} and extended analysis}
\label{sec_proof,theo_finitebdd_single}
This subsection is devoted to an (extended) analysis and proof of Theorem~\ref{theo_finitebdd_single}.
In part of our analysis, as an alternative to our Lipschitz stochastic noise Assumption~\ref{assu_smoothnoise}, we can impose the following \emph{individual convexity and smoothness} condition~\cite{roux2012stochastic,johnson2013accelerating,defazio2014saga,nguyen2017sarah}:
\begin{assumption}[Individual convexity{/}smoothness]  \label{assu_Lmax_indcvx}
Almost surely, the (random) function $\theta \mapsto f(\theta; \xi)$ is convex, twice continuously differentiable and satisfies the Lipschitz condition
\begin{align}
\vecnorm{\nabla f(\theta; \xi) - \nabla f(\theta'; \xi)}{2}
\le
\Lmax \vecnorm{\theta - \theta'}{2}
\ a.s.,
\quad 
\text{for all pairs $\theta; \theta' \in \real^d$.}
\end{align}
\end{assumption}
\noindent
All Assumptions~\ref{assu_StrcvxSmooth} and~\ref{assu_noisethetastar} along with either Assumption~\ref{assu_smoothnoise} or~\ref{assu_Lmax_indcvx}, are standard in the stochastic optimization literature (cf.~\cite{nguyen2019new,asi2019stochastic,lei2020adaptivity}).
Note that Assumption~\ref{assu_Lmax_indcvx} implies Assumption~\ref{assu_smoothnoise} with constant $\Lmax$; in many statistical applications, the quantity $\Lmax$ can be significantly larger than $\sqrt{L^2 + \sglip^2}$ in magnitude.

With these assumptions in place, let us formalize the two cases in which we analyze the \ROOTSGD algorithm.  
We refer to these cases as the \emph{Lipschitz Stochastic Noise} case (or \LSNspace for short), and the \emph{Individually Smooth and Convex} case (or \ISCspace for short).
\begin{description}
\item[\LSN\ Case:]
Suppose that Assumptions~\ref{assu_StrcvxSmooth},~\ref{assu_noisethetastar} and~\ref{assu_smoothnoise} hold, and define
\begin{align}
\label{indnoncvx_case}
\stpszmax \defn \frac{1}{4L} \land \frac{\mu}{8\sglip^2}
.
\end{align}
\item[\ISC\ Case:]
Suppose that Assumptions~\ref{assu_StrcvxSmooth},~\ref{assu_noisethetastar} and~\ref{assu_Lmax_indcvx} hold, and define
\begin{align}
\label{indcvx_case}
\stpszmax \defn \frac{1}{4 \Lmax}
.
\end{align}
\end{description}
As the readers shall see immediately, $\kappamax$ is a key quantity that plays a pivotal role in our analysis for both cases.

\begin{theorem}[Unified nonasymptotic results, single-epoch \ROOTSGD]\label{theo_finitebdd_single_complete}
Suppose that the conditions in either the \LSNspace or \ISCspace Case are in force, and let the step sizes be chosen according to the protocol~\eqref{step-sizet} for some $\stpsz \in (0, \stpszmax]$, and assume that we use the following burn-in time:
\begin{align}\label{Tburnin_complete}
\Tburnin \mydefn \Big \lceil \frac{24}{\stpsz\mu} \Big \rceil.
\end{align}
Then, for any iteration $\niters \geq 1$, the iterate $\theta_T$ from Algorithm~\ref{algo_singleepoch} satisfies the bound
\begin{align}
\label{finitebdd_single_complete}
\Exs\| \nabla F(\theta_T) \|_2^2 
& \leq 
\frac{2700 \; \|\nabla F(\theta_0)\|_2^2}{\eta^2 \mu^2 (T+1)^2} 
+
\frac{28 \; \sigmastarsq}{T+1}
.
\end{align}
\end{theorem}
We provide the proof of Theorem \ref{theo_finitebdd_single_complete} in both the \LSNspace and \ISCspace cases; the \LSNspace case corresponds to Theorem \ref{theo_finitebdd_single}.
In accordance with the discussion in \S\ref{sec_intro}, our nonasymptotic convergence rate upper bound \eqref{finitebdd_single_complete} for the expected squared gradient norm consists of the addition of two terms.
The first term, $\frac{\sigstar^2}{\niters}$, corresponds to the \emph{nonimprovable statistical error} depending on the noise variance at the minimizer.
The second term, which is equivalent to $\frac{\|\nabla F(\theta_0)\|_2^2\burnin^2}{\niters^2}$, corresponds to the \emph{bias} or \emph{optimization error} that indicates the polynomial forgetting from the initialization.
Theorem \ref{theo_finitebdd_single_complete} copes with a wide range of step sizes $\eta$: fixing the number of online samples $\niters$, \eqref{finitebdd_single_complete} asserts that the optimal asymptotic risk $\frac{\sigstar^2}{\niters}$ for the squared gradient holds up to an absolute constant whenever $\niters	\gtrsim	\frac{1}{\stpsz\mu}	\lor	\frac{\|\nabla F(\theta_0)\|_2^2}{\stpsz^2 \strongconvex^2 \sigstar^2}$.

Converting the convergence rate bound in \eqref{finitebdd_single_complete}, we can achieve a tight upper bound on the sample complexity to achieve a statistical estimator of $\thetastar$ with gradient norm bounded by $O(\varepsilon)$:%
\footnote{Indeed, we choose $\niters$ in Eq.~\eqref{finitebdd_single} to be sufficiently large such that it satisfies the inequalities $
T\ge	\Tburnin	=	\lceil \frac{24}{\eta\mu }\rceil
$, as well as $
\frac{2700 \|\nabla F(\theta_0)\|_2^2}{\eta^2 \strongconvex^2 \niters^2}	\le	\frac{\varepsilon^2}{2} 
$ and $
\frac{28\sigstar^2}{\niters}			\le	\frac{\varepsilon^2}{2}
$.
Here and on, we assume without loss of generality that $\varepsilon^2 \le \|\nabla F(\theta_0)\|_2^2$.
It is then straightforward to see that \eqref{complexity_finitebdd_single} serves as a tight sample complexity upper bound.
}
\begin{equation}\label{complexity_finitebdd_single_appendix}
\begin{aligned}
C_{\ref{theo_finitebdd_single}}(\varepsilon)
&=%%%%
\max\left\{
\frac{74}{\stpszmax\mu}\cdot\frac{\|\nabla F(\theta_0)\|_2}{\varepsilon}
, 
\frac{56\sigstar^2}{\varepsilon^2} 
\right\}
\\&\asymp%%%%
\begin{cases}
\max\left\{
\left(\frac{\smoothness}{\strongconvex} + \frac{\sglip^2}{\strongconvex^2} \right) \cdot \frac{\|\nabla F(\theta_0)\|_2}{\varepsilon}
,
\frac{\sigstar^2}{\varepsilon^2}
\right\}
,		& \text{for the \LSN case}
,
\\
\max\left\{
\frac{\Lmax}{\mu} \cdot \frac{\|\nabla F(\theta_0)\|_2}{\varepsilon}
,
\frac{\sigstar^2}{\varepsilon^2}
\right\}
,		& \text{for the \ISC case}
.
\end{cases}
\end{aligned}\end{equation}
In above, the step size $\stpsz = \stpszmax$ is optimized as in \eqref{indnoncvx_case} for the \LSN and \eqref{indcvx_case} for the \ISC case, separately, and where the asymptotics holds as $\varepsilon$ tends to zero while $\sigmastar$ is bounded away from zero.
In both cases, the leading-order term of $C_{\ref{theo_finitebdd_single}}(\varepsilon)$ in either case is $\asymp \frac{\sigstar^2}{\varepsilon^2}$ which matches the optimal statistical error up to universal constants, first among comparable literature in both cases.

\paragraph{Detailed proof.}
The rest of this subsection devotes to prove Theorem \ref{theo_finitebdd_single_complete}.
It is straightforward to show first \eqref{finitebdd_single_complete} automatically holds for $\niters < \burnin$ since for these $\niters$, $\theta_{\niters} = \theta_0$ and hence $\Exs\| \nabla F(\theta_{\niters}) \|_2^2  = \Exs\| \nabla F(\theta_0) \|_2^2$, so we only need to prove the result for $\niters \ge \burnin$.

We first define $\kappamax$ which is a key quantity in our analysis in this section for both cases, as follows
\begin{equation}\label{kappamax}
\kappamax	\defn		\begin{cases}
\frac{2\sglip^2}{\strongconvex^2},	& \text{for \LSN case}
,\\ 
\frac{2\Lmax}{\mu},		& \text{for \ISC case}
.
\end{cases}
\end{equation}
A central object in our analysis is the iteration of \emph{tracking error}, defined as
\begin{equation}\label{ztdefn}
z_t	\mydefn	v_t - \nabla F(\theta_{t-1})
,\qquad
\text{for $t\ge \burnin$}
.
\end{equation}
At a high level, this proof involves analyzing the evolution of the quantities $v_t$ and $z_t$, and then bounding the norm of the gradient $\nabla F(\theta_{t-1})$ using their combination.
From the updates~\eqref{algoROOTSGD}, we can identify a martingale difference structure for the quantity $tz_t$: its difference decomposes as the sum of \emph{pointwise stochastic noise}, $\noise_t(\theta_{t-1})$, and the \emph{incurred displacement noise}, $
(t-1)\left[\noise_t(\theta_{t-1}) - \noise_t(\theta_{t-2}) \right]
$.
The expression of the martingale structure is expressed as
\begin{align}\label{eq:zt_defn}
\begin{aligned}
tz_t
= 
t\left(v_t - \nabla F (\theta_{t - 1}) \right)
&= 
\noise_t(\theta_{t - 1}) + (t-1) (v_{t - 1} - \nabla F(\theta_{t - 2})) + (t-1) (\noise_t (\theta_{t - 1}) - \noise_t (\theta_{t - 2})) 
\\&= 
\noise_t(\theta_{t - 1}) + (t-1) z_{t - 1} + (t-1) (\noise_t (\theta_{t - 1}) - \noise_t (\theta_{t - 2}))
.
\end{aligned}\end{align}
Unwinding this relation recursively yields the decomposition
\begin{align}\label{eq:main-z-decomposition}  
\begin{aligned}
tz_{t} - \burnin z_{\burnin} 
&= 
\sum_{s = \burnin+1}^{t} \noise_s(\theta_{s - 1})
+ 
\sum_{s = \burnin+1}^{t} (s - 1) (\noise_s(\theta_{s - 1}) - \noise_s(\theta_{s - 2}))
.
\end{aligned}\end{align}

We now turn to the proofs of the three auxiliary lemmas that allow us to control the relevant quantities and the main theorem, as follows:
\begin{lemma}[Recursion involving $z_t$]
\label{lemm_estimationerror}
Under the conditions of Theorem~\ref{theo_finitebdd_single_complete}, for all $t \geq
\Tburnin + 1$, we have
\begin{subequations}
\begin{align}\label{estimationerror}  
t^2 \Exs\| z_t \|_2^2 
\le 
(t-1)^2 \Exs \| z_{t-1} \|_2^2 
+ 
2 \Exs \|\noise_t(\theta_{t-1}) \|_2^2 
+ 
2(t-1)^2 \Exs\| \noise_t(\theta_{t-1}) - \noise_t(\theta_{t-2}) \|_2^2
.
\end{align}
On the other hand, for $t = \Tburnin$, we have
\begin{align}\label{tequalburnin}
\burnin^2 \Exs\| v_\Tburnin \|_2^2 - \burnin^2 \Exs\vecnorm{\nabla F(\theta_0)}{2}^2 
= 
\burnin^2 \Exs \| z_\Tburnin \|_2^2 
= 
\Tburnin \Exs \| \noise_\Tburnin(\theta_0) \|_2^2
.
\end{align}
\end{subequations}
\end{lemma}
\noindent See \S\ref{SecProoflemm_estimationerror} for the proof of this claim. 
Note we have $z_\Tburnin = v_\Tburnin - \nabla F(\theta_0)$ which is simply the arithmetic average of $\Tburnin$ i.i.d.~noise terms at $\theta_0$, $\noise_1(\theta_{0}), \dots, \noise_\Tburnin(\theta_{0})$.

Our next auxiliary lemma characterizes the evolution of the sequence $(v_t: t \geq \Tburnin)$ in terms of the quantity $\Exs \|v_t\|_2^2$.
\begin{lemma}[Evolution of $v_t$]\label{lemm_Deltathetabdd2}
Under the settings of Theorem \ref{theo_finitebdd_single_complete}, for any $\stpsz \in (0, \stpszmax]$, we have
\begin{subequations}
\begin{align}\label{Deltathetabdd2}
t^2 \Exs\|v_t\|_2^2 - 2 t \Exs \langle v_t, \nabla F(\theta_{t-1}) \rangle + \Exs\left\| \nabla F(\theta_{t-1}) \right\|_2^2 
& =
\Exs\left\| t v_t - \nabla F(\theta_{t-1}) \right\|_2^2
,
\end{align}
and
\begin{multline}\label{Deltathetabdd2_inq}
\Exs \left \| t v_t - \nabla F(\theta_{t-1}) \right\|_2^2 
\le 
(1 - \stpsz\strongconvex) \cdot (t-1)^2 \Exs\|v_{t-1}\|_2^2 + 2\Exs\|\noise_t(\theta_{t-1})\|_2^2 
\\
- 2 (t-1)^2 \Exs \left \|\noise_t(\theta_{t-1}) - \noise_t(\theta_{t-2}) \right\|_2^2
,
\end{multline}
\end{subequations}
for all $t \ge \Tburnin + 1$.
\end{lemma}
\noindent See \S\ref{SecProoflemm_Deltathetabdd2} for the proof of this claim. 
\\

Our third auxiliary lemma bounds the second moment of the stochastic noise.
\begin{lemma}[Second moment of pointwise stochastic noise]\label{lemm_noisebatch}
Under the conditions of Theorem~\ref{theo_finitebdd_single_complete}, we have
\begin{align}
\label{noisebatch}
\Exs \|\noise_t(\theta_{t-1})\|_2^2 
& \le 
\kappamax \Exs\| \nabla F(\theta_{t-1}) \|_2^2 + 2 \sigstar^2
,
\qquad 
\mbox{for all $t \geq \Tburnin + 1$.}
\end{align}
\end{lemma}
\noindent See \S\ref{SecProoflemm_noisebatch} for the proof of this claim. 
\\

Equipped with these three auxiliary results, we are now ready to prove Theorem~\ref{theo_finitebdd_single_complete}. 

\begin{proof}[Proof of Theorem~\ref{theo_finitebdd_single_complete}]
Our proof proceeds in two steps.

\noindent\textbf{Step 1.}
We begin by applying the Cauchy-Schwarz and Young inequalities to the inner product $\inprod{v_t}{\nabla F(\theta_{t-1})}$. 
Doing so yields the upper bound
\begin{align*}
2 t \inprod{v_t}{\nabla F(\theta_{t-1})} 
\le 
2 \left[ t\| v_t \|_2\cdot \|\nabla F(\theta_{t-1}) \|_2 \right] 
\le 
\stpsz\strongconvex t^2 \|v_t\|_2^2 + \frac{1}{\stpsz\strongconvex} \| \nabla F(\theta_{t-1})\|_2^2 
.
\end{align*}
Taking the expectation of both sides and applying the bound~\eqref{Deltathetabdd2} from Lemma~\ref{lemm_Deltathetabdd2} yields
\begin{align*}
(1 - \stpsz\strongconvex) t^2 \Exs\| v_t\|_2^2 - \frac{1 - \stpsz\strongconvex}{\stpsz\strongconvex} \Exs\| \nabla F(\theta_{t-1})\|_2^2 & \le t^2 \Exs\|v_t\|_2^2 - 2t
\Exs\langle v_t, \nabla F(\theta_{t-1}) \rangle + \Exs\left\| \nabla F(\theta_{t-1}) \right\|_2^2 
\\& \le 
(1 - \stpsz\strongconvex) \cdot (t-1)^2 \Exs\|v_{t-1}\|_2^2 + 2\Exs\|\noise_t(\theta_{t-1})\|_2 ^2 
\\ & \hspace{1in}
- 2(t-1)^2 \Exs \|\noise_t(\theta_{t-1}) - \noise_t(\theta_{t-2}) \|_2^2
.
\end{align*}
Moreover, since we have $\stpsz \le \stpszmax \le \frac{1}{4\strongconvex}$ under condition~\eqref{eq:step-size-and-burnin}, we can multiply both sides by
$(1-\stpsz\strongconvex)^{-1}$, which lies in $[1, \frac{3}{2}]$.
Doing so yields the bound
\begin{align*}
  t^2 \Exs\| v_t\|_2^2 - \frac{1}{\stpsz\strongconvex} \Exs\| \nabla
  F(\theta_{t-1})\|_2^2 \le (t-1)^2 \Exs\|v_{t-1}\|_2^2 + 3
  \Exs\|\noise_t(\theta_{t-1})\|_2^2 - 2 (t-1)^2
  \Exs\left\|\noise_t(\theta_{t-1}) -
  \noise_t(\theta_{t-2})\right\|_2^2.
\end{align*}
Combining with the bound~\eqref{estimationerror} from
Lemma~\ref{lemm_estimationerror} gives
\begin{multline*}
t^2 \Exs\| z_t \|_2^2 + t^2 \Exs\| v_t\|_2^2 - (t-1)^2 \Exs \| z_{t-1} \|_2^2 - (t-1)^2 \Exs\|v_{t-1}\|_2^2 
\\\le 
5 \Exs\|\noise_t(\theta_{t-1})\|_2^2 + \frac{1}{\stpsz\strongconvex} \Exs\|\nabla F(\theta_{t-1})\|_2^2 
.
\end{multline*}
By telescoping this inequality from $\Tburnin+1$ to $\niters$, we find that
\begin{multline}  \label{Deltathetabdd2v2}
\niters^2 \Exs\| z_{\niters} \|_2^2 + \niters^2 \Exs\| v_{\niters}\|_2^2 - \burnin^2\Exs\| z_\Tburnin \|_2^2 -  \burnin^2 \Exs\|v_\Tburnin\|_2^2 
\\ \le 
\sum_{t=\Tburnin+1}^\niters\left[
   5 \Exs\|\noise_t(\theta_{t-1})\|_2^2 + \frac{1}{\stpsz\strongconvex} \Exs\|   \nabla F(\theta_{t-1})\|_2^2 
\right]
.
\end{multline}
Next, applying the result~\eqref{tequalburnin} from
Lemma~\ref{lemm_estimationerror} yields
\begin{align*}
%&\quad\,%%%%
\lefteqn{
\frac{\niters^2}{2} \Exs\| \nabla F(\theta_{\niters-1}) \|_2^2 
\le 
\niters^2 \Exs\| z_{\niters} \|_2^2 + \niters^2 \Exs\| v_{\niters}\|_2^2 
}
\\&\le
\burnin^2\Exs\| z_\Tburnin \|_2^2 
+ 
\burnin^2\Exs\|v_\Tburnin\|_2^2 
+ 
\sum_{t=\Tburnin+1}^\niters\left[
5\Exs\|\noise_t(\theta_{t-1})\|_2^2 + \frac{1}{\stpsz\strongconvex} \Exs\| \nabla F(\theta_{t-1})\|_2^2 
\right]
%%%%
\\&= 
\burnin^2 \|\nabla F(\theta_0) \|_2^2 
+ 
2\Tburnin \Exs\| \noise_\Tburnin(\theta_0) \|_2^2 
+ 
5\sum_{t=\Tburnin+1}^\niters \Exs\|\noise_t(\theta_{t-1})\|_2^2 
+ 
\frac{1}{\stpsz\strongconvex} \sum_{t=\Tburnin+1}^\niters \Exs\| \nabla F(\theta_{t-1})\|_2^2 
.
\end{align*}

Following some algebra, we find that
\begin{multline}
\label{gradient2ndmom_imp}
\Exs\| \nabla F(\theta_{\niters-1}) \|_2^2 
\le 
\frac{2\burnin^2 \|\nabla F(\theta_0) \|_2^2 + 4\Tburnin \Exs\| \noise_\Tburnin(\theta_0) \|_2^2}{\niters^2}  
\\ 
+
\frac{10}{\niters^2}\sum_{t=\Tburnin+1}^\niters \Exs\|\noise_t(\theta_{t-1})\|_2^2
+ 
\frac{2}{\stpsz\strongconvex \niters^2} \sum_{t=\Tburnin+1}^\niters \Exs\| \nabla F(\theta_{t-1})\|_2^2 
.
\end{multline}
Combining inequality~\eqref{gradient2ndmom_imp} with the bound~\eqref{noisebatch} from Lemma~\ref{lemm_noisebatch} gives
\begin{multline*}
\Exs\| \nabla F(\theta_{\niters-1}) \|_2^2 
\le \frac{2\burnin^2 \|\nabla F(\theta_0) \|_2^2 + 4\Tburnin
\left[ \kappamax \Exs\| \nabla F(\theta_0) \|_2^2 + 2\sigstar^2
 \right]}{\niters^2} \\
\qquad \qquad \qquad + \frac{10}{\niters^2}\sum_{t=\Tburnin+1}^\niters \left[
\kappamax \Exs\| \nabla F(\theta_{t-1}) \|_2^2 + 2\sigstar^2
\right] + \frac{2}{\stpsz\strongconvex \niters^2} \sum_{t=\Tburnin+1}^\niters \Exs\| \nabla
F(\theta_{t-1})\|_2^2 \\
%%%%
\le \frac{(4\kappamax + 2\Tburnin) \Tburnin \Exs\| \nabla
F(\theta_0) \|_2^2}{\niters^2} + \frac{10\kappamax +
2\strongconvex^{-1}\stpsz^{-1}}{\niters^2} \sum_{t=\Tburnin+1}^\niters \Exs\| \nabla
F(\theta_{t-1})\|_2^2 + \frac{20\sigstar^2}{\niters},
\end{multline*}
concluding  the following key gradient bound that controls the evolution of the gradient norm $\|\nabla F(\theta_{\niters-1})\|_2$:
\begin{align}  \label{zestimatebdd}
\Exs \| \nabla F(\theta_{\niters-1}) \|_2^2 & \le \frac{1}{\niters^2} \left \{
\aone \Exs\| \nabla F(\theta_0) \|_2^2 + \atwo \sum_{t=\Tburnin+1}^\niters
\Exs\| \nabla F(\theta_{t-1})\|_2^2 \right \} +
\frac{20\sigstar^2}{\niters},
\end{align}
where $\aone \mydefn (4 \kappamax + 2\Tburnin) \; \Tburnin$ and $\atwo
\mydefn 10 \kappamax + \frac{2}{\stpsz\strongconvex}$.

\noindent\textbf{Step 2.}
Based on the estimation bound \eqref{zestimatebdd}, the proof of Theorem~\ref{theo_finitebdd_single_complete} relies on a bootstrapping argument in order to remove the dependence of the right-hand side of Eq.~\eqref{zestimatebdd} on the quantity $\Exs\| \nabla F(\theta_{t-1})\|_2^2$.  
Let $\Tarb \ge \Tburnin+1$ be arbitrary.  
Telescoping the bound~\eqref{zestimatebdd} over the iterates $\niters=\Tburnin+1, \dots, \Tarb$ yields
\begin{align*}
  \sum_{\niters=\Tburnin+1}^{\Tarb} \Exs\| \nabla F(\theta_{\niters-1}) \|_2^2 &
  \le \underbrace{\aone \sum_{\niters=\Tburnin+1}^{\Tarb}
    \frac{\|\nabla F(\theta_0) \|_2^2 }{\niters^2}}_{\Quant_1} +
  \underbrace{\sum_{\niters=\Tburnin+1}^{\Tarb} \frac{\atwo}{\niters^2}
    \sum_{t=\Tburnin+1}^\niters \Exs \|\nabla F(\theta_{t-1}) \|_2^2}_{\Quant_2} +
    \underbrace{\sum_{\niters=\Tburnin+1}^{\Tarb} \frac{20\sigstar^2}{\niters}}_{\Quant_3}.
\end{align*}
Let us deal with each of these quantities in turn, making use of the integral inequalities
\begin{align}  \label{EqnIntegral}
\sum_{\niters=\Tburnin+1}^{\Tarb} \frac{1}{\niters^2} 
\stackrel{(i)}{\le}
\int_{\Tburnin}^\Tarb \frac{d \tau}{\tau^2} \le \frac{1}{\Tburnin}
,
\qquad \mbox{and} \quad 
\sum_{\niters=\Tburnin+1}^{\Tarb} \frac{1}{\niters}
\stackrel{(ii)}{\le} 
\int_{\Tburnin}^\Tarb \frac{d \tau}{\tau} = \log \Big(\frac{\Tarb}{\Tburnin} \Big)
.
\end{align}
We clearly have 
\begin{align*}
\Quant_1 
& \le 
\frac{\aone}{\Tburnin} \vecnorm{\nabla F(\theta_0)}{2}^2
\; = \; 
(4 \kappamax + 2\Tburnin) \; \vecnorm{\nabla F(\theta_0)}{2}^2
.
\end{align*}
Moreover, by using the fact that $\Tarb \geq T$, interchanging the
order of summation, and then using inequality~\eqref{EqnIntegral}(i)
again, we have
\begin{align*}
\Quant_2 
\le 
\sum_{\niters=\Tburnin+1}^{\Tarb} \frac{\atwo}{\niters^2} \sum_{t=\Tburnin+1}^\Tarb \Exs \|\nabla F(\theta_{t-1}) \|_2^2 
& =
\sum_{t=\Tburnin+1}^\Tarb \Big( \sum_{\niters=\Tburnin+1}^{\Tarb} \frac{\atwo}{\niters^2} \Big) \Exs \|\nabla F(\theta_{t-1}) \|_2^2 
\\& \le 
\frac{\atwo}{\Tburnin} \sum_{t=\Tburnin+1}^\Tarb  \Exs \|\nabla F(\theta_{t-1}) \|_2^2
.
\end{align*}
Finally, turning to the third quantity, we have $\Quant_3 \le 20\sigstar^2 \log \big( \frac{\Tarb}{\Tburnin} \big)$, where we have used inequality~\eqref{EqnIntegral}(ii).  
Putting together the pieces yields the upper bound
\begin{align*}
\sum_{\niters=\Tburnin+1}^{\Tarb} \Exs\| \nabla F(\theta_{\niters-1}) \|_2^2 
&\le 
(4\kappamax + 2\Tburnin) \|\nabla F(\theta_0) \|_2^2 
+
\frac{\atwo}{\Tburnin}\sum_{t=\Tburnin+1}^{\Tarb} \Exs\left\|\nabla F(\theta_{t-1}) \right\|_2^2 
+ 
20\sigstar^2 \log\left(\frac{\Tarb}{\Tburnin}\right)
.
\end{align*}
Eqs.~\eqref{eq:step-size-and-burnin} imply that, for either case under consideration, we have the bound $\kappamax \le \frac{1}{\stpsz\strongconvex}$, and, since $0< \stpsz \strongconvex \le \frac{1}{4} < 1$, we have from \eqref{eq:step-size-and-burnin} that $
\Tburnin
=
\left\lceil \frac{24}{\stpsz\strongconvex}\right\rceil
\le
\frac{1}{\stpsz\strongconvex}$, resulting in
\begin{align*}
4\kappamax + 2\Tburnin
&\le
\frac{4}{\stpsz\strongconvex} + 2\left( \frac{1}{\stpsz\strongconvex} \right)
=
\frac{54}{\stpsz\strongconvex}
,
\end{align*}
where we have the choice of burn-in time $\Tburnin$ from Eq.~\eqref{eq:step-size-and-burnin}.
Similarly, we have $
\atwo 
	= 
10 \kappamax + \frac{2}{\stpsz\strongconvex} 
	\le 
\frac{12}{\stpsz\strongconvex}
	\le
\frac{\Tburnin}{2}
$.
Putting together the pieces yields
\begin{align}\label{EqnBootstrapBound}
\frac{1}{2} \sum_{t=\Tburnin+1}^{\Tarb} \Exs \|\nabla F(\theta_{t-1}) \|_2^2 
& \le 
\frac{54}{\stpsz\strongconvex} \Exs \|\nabla F(\theta_0) \|_2^2 
+ 
20 \sigstar^2 \log \left(\frac{\Tarb}{\Tburnin} \right)
.
\end{align}
Now substituting the inequality~\eqref{EqnBootstrapBound} back into the earlier bound~\eqref{zestimatebdd} with $\Tarb = T$ allows us to obtain a bound on $\Exs \| \nabla F(\theta_{\niters-1}) \|_2$.  
In particular, for any $\niters \ge \Tburnin + 1$, we have
\begin{align*}
\Exs\| \nabla F(\theta_{\niters-1}) \|_2^2 
&\le 
\frac{54 \Tburnin}{\stpsz\strongconvex}\cdot \frac{\Exs\| \nabla F(\theta_0) \|_2^2}{\niters^2} 
+
\frac{\Tburnin}{\niters^2} \cdot \frac{1}{2} \sum_{t=\Tburnin+1}^\niters \Exs\|\nabla F(\theta_{t-1})\|_2 ^2 + \frac{20\sigstar^2}{\niters} 
\\&\le%%%%
\frac{54 \Tburnin}{\stpsz\strongconvex}\cdot \frac{\Exs\| \nabla F(\theta_0) \|_2^2}{\niters^2} 
+ 
\frac{\Tburnin}{\niters^2} \left[ \frac{54}{\stpsz\strongconvex} \|\nabla F(\theta_0) \|_2^2 
+ 
20 \sigstar^2 \log\left(\frac{\niters}{\Tburnin}\right) \right] + \frac{20 \sigstar^2}{\niters} 
\\& \le%%%%
\frac{2(54) \, \Tburnin}{\stpsz\strongconvex}\cdot \frac{\Exs\| \nabla F(\theta_0) \|_2^2}{\niters^2} 
+ 
\frac{20 \sigstar^2}{\niters} \left[ 1 + \frac{\Tburnin}{\niters} \log\left(\frac{\niters}{\Tburnin} \right) \right]
.
\end{align*}
Using the inequality $\frac{\log(x)}{x}\le \frac{1}{e}$, valid for $x\ge 1$, we
conclude that
\begin{align*}
\Exs\| \nabla F(\theta_{\niters-1}) \|_2^2 
&\le
\frac{2(54)\Tburnin}{\stpsz\strongconvex}\cdot \frac{\Exs\| \nabla F(\theta_0) \|_2^2}{\niters^2} 
+
\frac{20 \sigstar^2}{\niters} \left[ 1+\frac{\Tburnin}{\niters}\log
\left(\frac{\niters}{\Tburnin}\right) \right] 
\\& \le%%%%
\frac{108}{\stpsz\strongconvex}\cdot \frac{1}{\stpsz\strongconvex}\cdot \frac{\Exs\| \nabla F(\theta_0) \|_2^2}{\niters^2} 
+ 
\frac{20 \sigstar^2}{\niters} \left[ 1 + \frac{1}{e} \right] 
\\& \le%%%%
\frac{2700\; \vecnorm{\nabla F(\theta_0)}{2}^2}{\stpsz^2 \strongconvex^2 \niters^2} + \frac{28 \sigstar^2}{\niters}
.
\end{align*}
Shifting the subscript forward by one yields Theorem~\ref{theo_finitebdd_single_complete}.
\end{proof}

\pb\subsubsection{Proof of Lemma~\ref{lemm_estimationerror}}
\label{SecProoflemm_estimationerror}

The claim~\eqref{tequalburnin} follows from the definition along with some basic probability. 
In order to prove the claim~\eqref{estimationerror}, recall from the \ROOTSGD update rule for $v_t$ in the first line of \eqref{algoROOTSGD} that for $t\ge \Tburnin+1$ we have:
\begin{align}
tv_t &= (t-1)v_{t-1} + t \nabla f(\theta_{t-1}; \xi_t) - (t-1)\nabla
f(\theta_{t-2}; \xi_t) .
\end{align}
Subtracting the quantity $t\nabla F(\theta_{t-1})$ from both sides
yields
\begin{align*}
t z_t
& = 
(t-1)v_{t-1} + t\nabla f(\theta_{t-1}; \xi_t) - (t-1) \nabla f(\theta_{t-2}; \xi_t) - t\nabla F(\theta_{t-1})
.
\end{align*}
Thus, we arrive at the following recursion for the estimation error
$z_t$:
\begin{align*}
t z_t 
&= 
(t-1)\left[ v_{t-1}  - \nabla F(\theta_{t-2}) \right] 
\\&\quad\,
+ 
t \left[ \nabla f(\theta_{t-1};\xi_{t}) - \nabla F(\theta_{t-1}) \right] - (t-1) \left[ \nabla f(\theta_{t-2}; \xi_{t}) - \nabla   F(\theta_{t-2})\right] 
\\&= 
(t-1) z_{t-1} +  \noise_t(\theta_{t-1}) + (t-1)\left[ \noise_t(\theta_{t-1})    -\noise_t(\theta_{t-2})\right] 
.
\end{align*}
Observing that the variable $\noise_t(\theta_{t-1}) + (t-1)\left[  \noise_t(\theta_{t-1}) -\noise_t(\theta_{t-2})\right]$, defines
an $L^2$-martingale-difference
sequence, we see that
\begin{align*}
t^2\Exs\| z_t \|_2^2 
&= 
\Exs\left\|(t-1) z_{t-1}\right\|_2^2 
+  
\Exs\left\|\noise_t(\theta_{t-1}) + (t-1)\left[ \noise_t(\theta_{t-1}) - \noise_t(\theta_{t-2})\right]\right\|_2^2 
\\ &\le 
(t-1)^2\Exs\| z_{t-1} \|_2^2 
+
2\Exs\| \noise_t(\theta_{t-1}) \|_2^2 
+ 
2(t-1)^2 \Exs\| \noise_t(\theta_{t-1}) - \noise_t(\theta_{t-2}) \|_2^2 
,
\end{align*}
where in the last step follows from Young's inequality.
Computing the constants out completes the proof of the claim~\eqref{estimationerror}.

%%%%%%%%%%%%%%%%%%%%%%%%%%%%%%%%%%%%%%%%%%%%%%%%%%%%%%%%%%%%%%%%%%%%%%%%%%%%%%%%%%%%%%%%%%%

\pb\subsubsection{Proof of Lemma~\ref{lemm_Deltathetabdd2}}\label{SecProoflemm_Deltathetabdd2}
Eq.~\eqref{Deltathetabdd2} follows in a straightforward manner by expanding the square and taking an expectation.
As for the inequality~\eqref{Deltathetabdd2_inq}, from the update rule~\eqref{algoROOTSGD} for $v_t$, we have
\begin{align*}
  t v_t - \nabla F(\theta_{t-1}) & = t\nabla f(\theta_{t-1}; \xi_t) +
  (t-1)\left[ v_{t-1} - \nabla f(\theta_{t-2}; \xi_t) \right] - \nabla
  F(\theta_{t-1}) \\&= (t-1)v_{t-1} + (t-1)\left[ \nabla
    f(\theta_{t-1};\xi_t) - \nabla f(\theta_{t-2};\xi_t) \right] +
  \noise_t(\theta_{t-1}) .
\end{align*}
Using this relation, we can compute the expected squared Euclidean norm as
\begin{multline*}
\Exs\| t v_t - \nabla F(\theta_{t-1}) \|_2^2  = \Exs\left\| (t-1)
v_{t-1} + (t-1) \left[ \nabla f(\theta_{t-1};\xi_t) - \nabla
  f(\theta_{t-2}; \xi_t) \right] + \noise_t(\theta_{t-1}) \right\|_2^2
\nonumber \\
 = \Exs\left\|(t-1)v_{t-1} \right\|_2^2 + \Exs\left\|(t-1)\left[
   \nabla f(\theta_{t-1};\xi_t) - \nabla f(\theta_{t-2};\xi_t) \right]
 + \noise_t(\theta_{t-1})\right\|_2^2 \nonumber \\ + 2\Exs\left\langle
 (t-1)v_{t-1}, (t-1)\left[ \nabla f(\theta_{t-1};\xi_t) - \nabla
   f(\theta_{t-2};\xi_t) \right] + \noise_t(\theta_{t-1})\right\rangle.
 \nonumber
\end{multline*}
Further rearranging yields
\begin{multline}
  \label{vtsubt}
\Exs\| t v_t - \nabla F(\theta_{t-1}) \|_2^2 = (t-1)^2
\Exs\|v_{t-1}\|_2^2 + 2(t-1)^2 \Exs\left\|\nabla f(\theta_{t-1};\xi_t)
- \nabla f(\theta_{t-2};\xi_t)\right\|_2^2 \\ + 2\Exs\left\|
\noise_t(\theta_{t-1}) \right\|_2^2 + 2(t-1)^2 \Exs \left\langle
v_{t-1}, \nabla f(\theta_{t-1};\xi_t) - \nabla
f(\theta_{t-2};\xi_t)\right \rangle.
\end{multline}
We split the remainder of our analysis into two cases, corresponding to the \LSN case or the \ISC case.
The difference in the analysis lies in how we handle the term $\inprod{v_{t-1}}{\nabla F (\theta_{t-1}) - \nabla F (\theta_{t - 2})}$.

%%%%%%%%%%%%%%%%%%%%%%%%%%%%%%%%%%%%%%%%%%%%%%%%%%%%%%%%%%%%%%%%%%%%%%%%%%%%%%%%%%%%%%%
\paragraph{Analysis in the \LSN case:}
From $L$-Lipschitz smoothness of $F$ in Assumption~\ref{assu_StrcvxSmooth}, we have
\begin{equation}\label{nessmooth}
\begin{aligned}
\binprod{v_{t-1}}{\nabla F(\theta_{t-1}) - \nabla F(\theta_{t-2})} 
&=
-\frac{1}{\eta} \binprod{\theta_{t-1} - \theta_{t-2}}{\nabla F(\theta_{t-1}) - \nabla F(\theta_{t-2})}
\\
&\le
-\frac{1}{\eta L} \|\nabla F(\theta_{t-1}) - \nabla F(\theta_{t-2})\|_2^2
.
\end{aligned}\end{equation}
Now consider the inner product term $\inprod{v_{t - 1}}{\nabla  F(\theta_{t - 1}) - \nabla F (\theta_{t - 2})}$ in Eq.~\eqref{vtsubt}. 
We split it into two terms, and upper bound them using equations~\eqref{nesstrong} and~\eqref{nessmooth} respectively. 
Doing so yields:
\begin{align*}
&\quad\,%%%%
\Exs\| t v_t - \nabla F(\theta_{t-1}) \|_2^2 
%%%%%%%%
\\&\le
(t-1)^2\Exs\|v_{t-1}\|_2^2 + 2 (t-1)^2\Exs\left\|\nabla f(\theta_{t-1};\xi_t) - \nabla f(\theta_{t-2};\xi_t)\right\|_2^2 
\\&\quad\,
+ 2 \Exs\left\| \noise_t(\theta_{t-1}) \right\|_2^2 + 2(t-1)^2 \Exs\left\langle v_{t-1}, \nabla F(\theta_{t-1}) - \nabla F(\theta_{t-2}) \right\rangle 
%%%%%%%%
\\& \le 
(t-1)^2 \Exs\|v_{t-1}\|_2^2 + 2(t-1)^2 \Exs\left\|\nabla F(\theta_{t-1}) - \nabla F(\theta_{t-2})\right\|_2^2 + 2(t-1)^2 \Exs\left\|\noise_t(\theta_{t-1}) - \noise_t(\theta_{t-2})\right\|_2^2 
\\&\quad\,
+ 2\Exs\|\noise_t(\theta_{t-1})\|_2^2 - \frac{3\eta\mu}{2} (t-1)^2 \Exs\|v_{t-1}\|_2^2 - \frac{1}{2\eta L}(t-1)^2 \Exs\|\nabla F(\theta_{t-1}) - \nabla F(\theta_{t-2})\|_2^2 
%%%%%%%%
\\&\le
\left(1 - \frac{3\eta\mu}{2}\right) (t-1)^2 \Exs\|v_{t-1}\|_2^2 
+ 
2\Exs\|\noise_t(\theta_{t-1})\|_2^2 
+ 
4(t-1)^2 \Exs\left\|\noise_t(\theta_{t-1}) - \noise_t(\theta_{t-2})\right\|_2^2
\\&\quad\,%
- 2(t-1)^2 \Exs\left\|\noise_t(\theta_{t-1}) - \noise_t(\theta_{t-2})\right\|_2^2	
%%%%%%%%
\\&\le
\left(1 - \frac{3\eta\mu}{2} + 4\eta^2\sglip^2 \right) (t-1)^2 \Exs\|v_{t-1}\|_2^2 
+ 
2\Exs\|\noise_t(\theta_{t-1})\|_2^2 
- 
2(t-1)^2 \Exs\left\|\noise_t(\theta_{t-1}) - \noise_t(\theta_{t-2})\right\|_2^2
.
\end{align*}
From the condition \eqref{indcvx_case}, we have \mbox{$
1 - \frac{3}{2} \eta \mu + 4 \eta^2 \sglip^2 \le 1 - \eta \mu
$,} which
completes the proof.

%%%%%%%%%%%%%%%%%%%%%%%%%%%%%%%%%%%%%%%%%%%%%%%%%%%%%%%%%%%%%%%%%%%%%%%%%%%%%%%%%%%%%%

\paragraph{Analysis in the \ISC case:}
We deal with the last summand in the last line of Eq.~\eqref{vtsubt}, where we use the iterated law of expectation to achieve
\begin{align*}
\Exs \left\langle 
v_{t-1}, \nabla f(\theta_{t-1};\xi_t) - \nabla f(\theta_{t-2};\xi_t)
\right\rangle 
&= 
\Exs \left\langle 
v_{t-1}, \Exs\left[\nabla f(\theta_{t-1};\xi_t) - \nabla f(\theta_{t-2};\xi_t) \mid \cF_{t-1}\right] 
\right\rangle
\\&= 
\Exs \left\langle 
v_{t-1}, \nabla F(\theta_{t-1}) - \nabla F(\theta_{t-2})
\right\rangle
.
\end{align*}
The update rule for $v_t$ implies that $
v_{t-1}	=	-\frac{\theta_{t-1} - \theta_{t-2}}{\eta}
$ for all $t \ge \Tburnin + 1$.
The following analysis uses various standard inequalities (c.f.~\S 2.1 in~\cite{NESTEROV[Lectures]}) that hold for individually convex and $\Lmax$-Lipschitz smooth functions.
First, we have
\begin{equation}\label{deallast}  
\begin{aligned}
\left\langle v_{t-1}, \nabla f(\theta_{t-1};\xi_t) - \nabla f(\theta_{t-2};\xi_t) \right\rangle 
&= 
-\frac{1}{\eta} \left\langle \theta_{t-1} - \theta_{t-2}, \nabla f(\theta_{t-1};\xi_t) - \nabla f(\theta_{t-2};\xi_t) \right\rangle 
\\& \le 
-\frac{1}{\eta\Lmax} \left\| \nabla f(\theta_{t-1};\xi_t) - \nabla f(\theta_{t-2};\xi_t) \right\|_2^2
,
\end{aligned}\end{equation}
where the inequality follows from the Lipschitz condition.  
On the other hand, the $\mu$-strong convexity of $F$ implies that
\begin{equation}\label{nesstrong}
\begin{aligned}
\left\langle v_{t-1}, \nabla F(\theta_{t-1}) - \nabla F(\theta_{t-2})\right\rangle 
&= 
-\frac{1}{\eta} \left\langle \theta_{t-1} - \theta_{t-2}, \nabla F(\theta_{t-1}) - \nabla F(\theta_{t-2}) \right\rangle 
\\
& \le 
-\frac{\mu}{\eta} \|\theta_{t-1} - \theta_{t-2}\|_2^2 
= 
-\eta \mu \|v_{t-1}\|_2^2
.
\end{aligned}\end{equation}
Plugging the bounds~\eqref{deallast} and~\eqref{nesstrong} into Eq.~\eqref{vtsubt} yields
\begin{align*}
&\quad\,%%%%
\Exs\| tv_t - \nabla F(\theta_{t-1}) \|_2^2 
%%%%%%%%
\\&\le
(t-1)^2\Exs\|v_{t-1}\|_2^2 
+ 
2 (t-1)^2 \Exs\left\|\nabla f(\theta_{t-1};\xi_t) - \nabla f(\theta_{t-2};\xi_t)\right\|_2^2 
+ 
2 \Exs\left\| \noise_t(\theta_{t-1}) \right\|_2^2 
\\&\quad\,%%%%
+ 
(t-1)^2 \Exs \left\langle v_{t-1}, \nabla F(\theta_{t-1}) - \nabla F(\theta_{t-2}) \right\rangle 
+ 
(t-1)^2 \Exs \left\langle v_{t-1}, \nabla f(\theta_{t-1};\xi_t) - \nabla f(\theta_{t-2};\xi_t) \right\rangle 
%%%%%%%%
\\& \le 
(t-1)^2 \Exs\|v_{t-1}\|_2^2 
+ 
2(t-1)^2 \Exs\left\|\nabla f(\theta_{t-1};\xi_t) - \nabla f(\theta_{t-2};\xi_t)\right\|_2^2 
+
2\Exs\|\noise_t(\theta_{t-1})\|_2^2 
\\&\quad\,%%%%
- 
\eta \mu (t-1)^2 \Exs \|v_{t-1}\|_2^2 
- 
\frac{1}{\eta \Lmax} (t-1)^2 \Exs \|\nabla f(\theta_{t-1};\xi_t) 
- 
\nabla f(\theta_{t-2};\xi_t) \|_2^2 
%%%%%%%%
\\&\le 
(1 - \eta \mu) (t-1)^2 \Exs\|v_{t-1}\|_2^2 
+ 
2 \Exs\|\noise_t(\theta_{t-1})\|_2^2 
- 
2(t-1)^2 \Exs \|\nabla f(\theta_{t-1};\xi_t) - \nabla f(\theta_{t-2};\xi_t) \|_2^2
,
\end{align*}
where in the last inequality relies on the fact that $\eta \in (0, \frac{1}{4 \Lmax}]$ (see Eq.~\eqref{indcvx_case}), leading to the bound~\eqref{Deltathetabdd2_inq}.

%%%%%%%%%%%%%%%%%%%%%%%%%%%%%%%%%%%%%%%%%%%%%%%%%%%%%%%%%%%%%%%%%%%%%%%%%%%%%%%%%%%%%%

\pb\subsection{Proof of Lemma~\ref{lemm_noisebatch}}
\label{SecProoflemm_noisebatch}

We again split our analysis into two cases, corresponding to the \LSNspace and \ISC cases.
Recall that the main difference is whether the Lipschitz stochastic noise condition holds (cf.\ Assumption~\ref{assu_smoothnoise}), or the functions are individually convex and smooth (cf.\  Assumption~\ref{assu_Lmax_indcvx}).

\paragraph{\textbf{Analysis in the \LSN case:}}

From the $\sglip$-Lipschitz smoothness of the stochastic gradients (Assumption~\ref{assu_smoothnoise}) and the $\mu$-strong-convexity of $F$ (Assumption~\ref{assu_StrcvxSmooth}), we have
\begin{equation}\label{TriangleStandard}\begin{aligned}
\Exs \|\noise_t(\theta_{t-1})\|_2^2 & \le
2\Exs\|\noise_t(\theta_{t-1}) - \noise_t(\thetastar)\|_2^2 +
2\Exs\|\noise_t(\thetastar)\|_2^2 \\
& \le 2\sglip^2\Exs\|\theta_{t-1} - \thetastar\|_2^2 +
2\Exs\|\noise_t(\thetastar)\|_2^2 \\
& \le \frac{2\sglip^2}{\strongconvex^2} \Exs\| \nabla F( \theta_{t-1}) \|_2^2 +
2 \sigstar^2,
\end{aligned}\end{equation}
which establishes the claim.

%%%%%%%%%%%%%%%%%%%%%%%%%%%%%%%%%%%%%%%%%%%%%%%%%%%%%%%%%%%%%%%%%%%%%%%%%%%%%%%%%%%%%%%%%%%%%

\paragraph{\textbf{Analysis in the \ISC case:}}

Using Assumption~\ref{assu_Lmax_indcvx} and standard inequalities for
$\Lmax$-smooth and convex functions yields
\begin{align*}
f(\thetastar;\xi) + \langle \nabla f(\thetastar;\xi), \theta \rangle +
\frac{1}{2\Lmax} \|\nabla f(\theta;\xi) - \nabla
f(\thetastar;\xi)\|_2^2 \le f(\theta;\xi).
\end{align*}
Taking expectations in this inequality and performing some
algebra%
\footnote{In performing this algebra, we assume exchangeability
  of gradient and expectation operators, which is guaranteed because
  the function $x \mapsto \nabla f(x; \xi)$ is $\Lmax$-Lipschitz for
  a.s.~$\xi$.} yields
\begin{align*}
\Exs\|\nabla f(\theta;\xi) - \nabla f(\thetastar;\xi)\|_2^2 & =
2\Lmax\langle \Exs[\nabla f(\thetastar;\xi)], \theta \rangle +
\Exs\|\nabla f(\theta;\xi) - \nabla f(\thetastar;\xi)\|_2^2 \\
& \le 2\Lmax\Exs\left[ f(\theta;\xi) - f(\thetastar;\xi) \right] \\
& = 2\Lmax \left[ F(\theta) - F(\thetastar) \right].
\end{align*}
Recall that $\nabla F(\thetastar) = 0$ since $\thetastar$ is a
minimizer of $F$.  
Using this fact and the $\mu$-strong convexity condition, we have $F(\theta) - F(\thetastar) \le \frac{1}{2 \mu} \|\nabla F(\theta)\|_2^2$.  
Substituting back into our earlier inequality yields
\begin{align*}
\Exs\|\nabla f(\theta;\xi) - \nabla f(\thetastar;\xi)\|_2^2     & \le
\frac{\Lmax}{\mu} \| \nabla F(\theta) \|_2^2.
\end{align*}
We also note that%
\footnote{This proof strategy is forklore and appears elsewhere in the variance-reduction literature; see, e.g., the proof of Theorem 1 in \cite{johnson2013accelerating}, and also adopted by \cite{nguyen2019new,nguyen2021inexact}.}
\begin{align*}
  \Exs \left \| \noise_t(\theta_{t-1}) - \noise_t(\thetastar)
  \right\|_2^2 & = \Exs \left \| \nabla f(\theta_{t-1};\xi_t) - \nabla
  f(\thetastar;\xi_t) - [\nabla F(\theta_{t-1}) - \nabla
    F(\thetastar)] \right \|_2^2 \\
& \le \Exs \left\| \nabla f(\theta_{t-1};\xi_t) - \nabla
  f(\thetastar;\xi_t) \right\|_2^2 \\
& \le \frac{\Lmax}{\mu} \Exs\| \nabla F(\theta_{t-1}) \|_2^2 .
\end{align*}
Finally, applying the argument of \eqref{TriangleStandard} yields the claim~\eqref{noisebatch}.

\pb\subsection{Intermediate result Proposition \ref{theo_rootHolder_single} and its proof}\label{sec_proof,theo_rootHolder_single}
En route our proof of Theorem \ref{theo_rootHolder_multi} we state and prove an intermediate upper-bound result for single-epoch version of \ROOTSGD.
For our convenience we forgo tracking the universal constants (which can be change at each appearance) due to complications of our derivations.

%\begin{subtheorem}{theorem}\label{theo_rootHolder}
\begin{theorem2}[Improved nonasymptotic upper bound, single-epoch \ROOTSGD]\label{theo_rootHolder_single}
Under Assumptions \ref{assu_StrcvxSmooth}, \ref{assu_holderhessian}, \ref{assu_noise_bdd_holder}, \ref{assu_noise_smooth_holder}, suppose that we run Algorithm~\ref{algo_singleepoch} with step-size $\stpsz \in \big(0, \frac{1}{56\smoothness} \wedge \frac{\strongconvex}{64\sgliptild^2}\big]$.  
Then for any $\niters \geq 1$, the iterate $\theta_{\niters}$ satisfies the bound
\begin{align}
\Exs \normb{\naF{\theta_{\niters}}} - \frac{\sigstar^2}{\niters} 
&\leq 
C\Big\{ \frac{\sglip^2\stpsz}{\strongconvex} 
+
\frac{\log \niters}{\stpsz\strongconvex \niters} 
+ 
\frac{\sglip^2 \log\niters}{\strongconvex^2\niters} \Big\} \frac{\sigstar^2}{\niters}
%\\&\quad\,%%%%
+
\frac{C \holderconst \sigstartild^{3}}{\stpsz^{1/2}\strongconvex^{5/2} \niters^{2}} \notag
\\&\qquad%%%%
+
\frac{C\normb{\naF{\theta_0}}}{\stpsz^2 \strongconvex^2 \niters^2}
+
\frac{C \holderconst \norm{\naF{\theta_0}}^{3} }{\stpsz^{7/2} \strongconvex^{11/2}\niters^{7/2}} 
\label{rootHolder_single}
\end{align}
\end{theorem2}

A few remarks are in order.  
When setting $\niters\to \infty$ the leading-order term $(1 + \frac{C\sglip^2\stpsz}{\strongconvex}) \frac{\sigstar^2}{\niters}$ of the nonasymptotic bound~\eqref{rootHolder_single} nearly matches the optimal statistical risk for the gradient norm with unit pre-factor when $\stpsz$ is prescribed as positively small, and as will be seen later it matches the asymptotic Proposition \ref{theo_asymptotic_iteration_const_step_size} under a shared
umbrella of assumptions.
It can be observed that the dependence on the initial gradient norm $\vecnorm{\nabla F(\theta_0)}{2}$ decays polynomially, which is generally unavoidable for single-epoch \ROOTSGD, as the gradient noise at the initial point $\theta_0$ is also averaged along the iterates. 
However, as we will see anon, an improved guarantee can be obtained by appropriately re-starting the algorithm, leading to near-optimal guarantees in terms of the gradient norm. 
In addition, we note that the high-order terms of Eq.~\eqref{rootHolder_single} contains terms that depend on the step-size $\stpsz$ at opposite directions which demands a trade-off.
%To avoid unnecessary complications, we leave the question of optimizing the step-size to the upcoming multi-epoch algorithm.
We forgo optimizing the step-size as is the conduct in our multi-epoch result.

For the rest of \S\ref{sec_proof,theo_rootHolder_single} we prepare to prove Proposition \ref{theo_rootHolder_single}.
From the discussions in \S\ref{sec_hidden} we decomposes $\Exs \normb{\naF{\theta_{t - 1}}}$ as the summation of three terms:
\begin{equation}\label{decompvz}
\Exs\vecnorm{\naF{\theta_{t-1}}}{2}^2 
=
\Exs\vecnorm{v_t - z_t}{2}^2 
=
\Exs\vecnorm{v_t}{2}^2 
+
\Exs\vecnorm{z_t}{2}^2 
-
2 \Exs \binprod{v_t}{z_t}
.
\end{equation}
En route our proof, we provide estimations for $\Exs \normb{tv_t}$, $\Exs \normb{tz_t}$ and $\Exs \binprod{tz_t}{tv_t}$ separately, where our main focus will be on bounding the cross term.
On a very intuitive and high-level viewpoint, when comparing with the Polyak-Ruppert-Juditsky analysis, we can roughly think of the $(\stpsz t v_t: t\ge 0)$ process acts like a last-iterate SGD (as it is in the quadratic minimization case) and is \emph{fast} and \emph{small}. 
The $tz_t$ process more resembles random walk at a slower rate driven by the same noise sequence.
The two timescale intuitions beneath is that two fast-slow discounted random walks processes driven by the same noise has an inner product that is approximately the second moment of the fast process. 
In our case this results in the "asymptotically independence" of the two processes in the sense that $\Exs \binprod{tz_t}{tv_t}$ scales as $\Exs \normb{tv_t}$, so $\naF{\theta_{t - 1}} = v_t - z_t$ is approximately of the same scale as $z_t$ in its first and second orders.
%Also, we can inject randomness on $\theta_0$ as long as it is $\cF_0$-measurable.

We first introduce the following lemma which is an essential part of the proof:
\begin{lemma}[Sharp bound on $v_t$]\label{lemm_vtsharpbdd}
Under the setting of Theorem \ref{theo_finitebdd_single}, there exists a universal constant $c > 0$, such that for $\niters\ge \burnin + 1$, we have:
\begin{align}\label{vtsharpbdd}
\Exs \normb{v_{\niters}} 
\le
\frac{c\sigstar^2}{\stpsz\strongconvex \niters^2} 
+
\frac{c}{\stpsz^4 \strongconvex^4 \niters^4} \normb{\nabla F(\theta_0)}
.
\end{align}
\end{lemma}
\noindent 
We defer the proof of Lemma~\ref{lemm_vtsharpbdd} to \S\ref{sec_vtsharpbdd}.
This lemma, along with Theorem~\ref{theo_finitebdd_single}, helps conclude the following bound on $z_t$ that has a leading-order term of near-unity pre-factor, that is, $(1 + o(1)) \frac{\sigma_*}{\sqrt{t}}$:

\begin{lemma}[Sharp bound on $z_t$]\label{lemm_ztsharpbdd}
Under settings of Theorem \ref{theo_finitebdd_single}, the following bounds hold true for $\niters \ge \burnin + 1$:
\begin{align}
\Exs\vecnorm{z_\niters}{2}^2 - \frac{\sigstar^2}{\niters}
\leq
c \Big\{
\frac{\sglip^2\stpsz}{\strongconvex}
+
\frac{\sglip}{\strongconvex \sqrt{\niters}}
+
\frac{\log\left( \frac{\niters}{\burnin} \right)\sglip^2}{\strongconvex^2\niters}
\Big\} \frac{\sigstar^2}{\niters} 
+ 
c \frac{\sglip \sigstar}{\mu} \cdot \frac{\burnin}{\niters^2} \vecnorm{\nabla F(\theta_0)}{2}
+
c \frac{\burnin^2}{ \niters^2}\normb{\naF{\theta_0}}
, \label{ztsharpbdd}
\end{align}
for some universal constant $c > 0$.
\end{lemma}
\noindent See \S\ref{sec_ztsharpbdd} for the proof of this lemma.

Finally, we need the following lemma, which bounds the cross term $\Exs\inprod{v_t}{z_t}$. 
Under the Lipschitz condition on the Hessian matrix and additional moment conditions, this lemma provides significant sharper bound than the na\"{i}ve bound obtained by applying the Cauchy-Schwartz inequality and invoking the previous two lemmas.
\begin{lemma}[Sharp bound on the cross term]\label{lemm_cross-term-bound}
Under settings of Theorem \ref{theo_finitebdd_single}, we have the following bound for any $\niters \geq \burnin + 1$:
\begin{multline}
\left| \Exs \binprod{v_{\niters}}{z_{\niters}} \right|
\leq 
c\left(\frac{\sigstar^2}{\stpsz \strongconvex \niters^2}
+
\frac{ \normb{\naF{\theta_0}}}{\stpsz^4 \strongconvex^4 \niters^4} \right)\log \niters
%\\
+
c  \holderconst
\left( 
\frac{\sigstartild^{3} }{\stpsz^{1/2}\strongconvex^{5/2}\niters^{2}} 
+
\frac{ \norm{\naF{\theta_0}}^{3} }{\stpsz^{7/2} \strongconvex^{11/2}\niters^{7/2}} 
\right)
,
\end{multline}
for some universal constant $c > 0$.
\end{lemma}
\noindent See \S\ref{subsubsec:proof-cross-term-bound} for the proof of this lemma.

Taking the aforementioned lemmas as given, we are ready to prove the sharp bound. In particular, by substituting these three lemmas into the decomposition~\eqref{decompvz}, we have the following bound
\begin{align*}
&\Exs \normb{\nabla F(\theta_{\niters-1})}
-
\frac{\sigstar^2}{\niters}
=
\Exs \normb{v_{\niters} - z_{\niters}}
-
\frac{\sigstar^2}{\niters}
=
\left(\Exs \normb{z_{\niters}} - \frac{\sigstar^2}{\niters}\right)
+ 
\Exs \normb{v_{\niters}} - 2\Exs \binprod{v_{\niters}}{z_{\niters}}
	\notag \\&\leq 
C \left(\frac{\sglip^2\stpsz}{\strongconvex}
+\frac{\sglip}{\strongconvex \sqrt{\niters}}
+\frac{\log\left( \frac{\niters}{\burnin} \right)\sglip^2}{\strongconvex^2\niters}
\right)
\frac{\sigstar^2}{\niters}
%\\&\quad\,%%%%
+
C \left(
\frac{\sglip\sigstar}{\strongconvex}\cdot \frac{\burnin}{\niters^2} \norm{\naF{\theta_0}}
%\\&\quad\,%%%%
+
\frac{\sglip^2}{\strongconvex^2}\cdot\frac{\burnin}{\niters^2}\normb{\naF{\theta_0}}
\right)
	\notag \\&\quad\,%%%%
+C \left(\frac{\sigstar^2}{\stpsz \strongconvex \niters^{2}} + \frac{ \normb{\naF{\theta_0}}}{\stpsz^{4} \strongconvex^{4} \niters^{4}} \right) \log \niters
	 +
6 C_0  \holderconst
\left( 
\frac{\sigstartild^{3} }{\stpsz^{1/2}\strongconvex^{5/2}\niters^{2}} 
+
\frac{ \norm{\naF{\theta_0}}^{3} }{\stpsz^{7/2} \strongconvex^{11/2}\niters^{7/2}} 
\right)
	\notag \\&\leq 
C \left(\frac{\sglip^2\stpsz}{\strongconvex}
+\underbrace{
\frac{\sglip}{\strongconvex \sqrt{\niters}}
}
+\frac{\log\niters}{\stpsz\strongconvex \niters}
+\frac{\log\left( \frac{\niters}{\burnin} \right)\sglip^2}{\strongconvex^2\niters}
\right)
\frac{\sigstar^2}{\niters}
%\\&\quad\,%%%%
	+
\frac{C\holderconst \sigstartild^{3}}{\stpsz^{1/2}\strongconvex^{5/2} \niters^{2}}
	\notag \\&\quad\,%%%%
+C \left(
\frac{\normb{\naF{\theta_0}}}{\stpsz^2 \strongconvex^2 \niters^2}
+\underbrace{
\frac{\sglip\sigstar}{\strongconvex}\cdot \frac{\burnin}{\niters^2} \norm{\naF{\theta_0}}
}
%\\&\quad\,%%%%
\right)
	 +
C  
\frac{ \holderconst \norm{\naF{\theta_0}}^{3} }{\stpsz^{7/2} \strongconvex^{11/2}\niters^{7/2}} 
.
\end{align*}
Absorbing the bracketed cross terms into corresponding sum of the squares, this gives Eq.~\eqref{rootHolder_single} and concludes Proposition \ref{theo_rootHolder_single}.

\pb\subsubsection{Proof of Lemma \ref{lemm_vtsharpbdd}}\label{sec_vtsharpbdd}
Our main technical tools is the following Lemma~\ref{lemm_vtrecursive}, which recursively bound the second moments of $v_t$:
\begin{lemma}\label{lemm_vtrecursive}
Under the setting of Theorem~\ref{theo_finitebdd_single}, we have the following bound for $t\ge \burnin+1$
\begin{align}
t^2\Exs\vecnorm{v_t}{2}^2
\le
\left( 1 - \frac{\stpsz\strongconvex}{2} \right) (t-1)^2\Exs\vecnorm{v_{t - 1}}{2}^2
+
\frac{10}{\stpsz\strongconvex} \Exs\vecnorm{\nabla F (\theta_{t - 1})}{2}^2
+ 4\sigstar^2
.
\label{eq_vtrecursive}
\end{align}
\end{lemma}
\noindent 
See \S\ref{sec_proof,lemm_vtrecursive} for the proof of this lemma.

On the other hand, invoking Theorem~\ref{theo_finitebdd_single}, we have that
\begin{align*}
\Exs\| \nabla F(\theta_{t-1}) \|_2^2 
\le 
\frac{2700 \; \vecnorm{\nabla F(\theta_0)}{2}^2}{\stpsz^2 \strongconvex^2 t^2} + \frac{28 \; \sigstar^2}{t}
,\qquad\text{for  $t\ge \burnin + 1$}
.
\end{align*}
Now, to combine everything together, we conclude from \eqref{finitebdd_single} and \eqref{eq_vtrecursive} that
\begin{align}
t^2\Exs\vecnorm{v_t}{2}^2
&\le%%%%
\left( 1 - \frac{\stpsz \strongconvex}{2} \right) (t-1)^2\Exs\vecnorm{v_{t - 1}}{2}^2
+
\frac{10}{\stpsz\strongconvex} 
\left[
\frac{2700 \; \vecnorm{\nabla F(\theta_0)}{2}^2}{\stpsz^2 \strongconvex^2 t^2} + \frac{28 \; \sigstar^2}{t}
\right]
+ 4 \sigstar^2 \nonumber
\\&\le
\left( 1 - \frac{\stpsz \strongconvex}{2} \right) (t-1)^2\Exs\vecnorm{v_{t - 1}}{2}^2
+ c\frac{  \vecnorm{\nabla F(\theta_0)}{2}^2}{\stpsz^3 \strongconvex^3 t^2} 
+ c\sigstar^2
. \label{vtrec}
\end{align}
Multiplying both sides by $t^2$, we obtain that
\begin{align*}
t^4\Exs\vecnorm{v_t}{2}^2
&\le%%%%
\left( 1 - \frac{\stpsz \strongconvex}{2} \right) t^2(t-1)^2\Exs\vecnorm{v_{t - 1}}{2}^2
+\frac{ c \vecnorm{\nabla F(\theta_0)}{2}^2}{\stpsz^3 \strongconvex^3}
+  c\sigstar^2 t^2
%%%%%%%%
\\&\leq
\left( 1 - \frac{\stpsz \strongconvex}{6} \right) (t-1)^4\Exs\vecnorm{v_{t - 1}}{2}^2
+
\frac{c \vecnorm{\nabla F(\theta_0)}{2}^2}{\stpsz^3 \strongconvex^3}
+c \sigstar^2 t^2,
\end{align*}
for time index $t$ satisfying $t \geq \burnin \geq \frac{6}{\stpsz \strongconvex}$.
This gives, by solving the recursion,
\begin{align*}
\niters^4\Exs\vecnorm{v_{\niters}}{2}^2
&\leq
\left( 1 - \frac{\stpsz \strongconvex}{6} \right)^{\niters-\burnin} \burnin^4\Exs\vecnorm{v_{\burnin}}{2}^2
+
c \sum_{t=\burnin+1}^\niters
\left( 1 - \frac{\stpsz \strongconvex}{6} \right)^{\niters-t}\left(
\frac{ \; \vecnorm{\nabla F(\theta_0)}{2}^2}{\stpsz^3 \strongconvex^3} 
+ \sigstar^2 T^2
\right)\\
&\leq \left( 1 - \frac{\stpsz \strongconvex}{6} \right)^{\niters-\burnin} \burnin^4\Exs\vecnorm{v_{\burnin}}{2}^2
+ 6c \frac{ \vecnorm{\nabla F(\theta_0)}{2}^2}{\stpsz^4 \strongconvex^4}
+ 6c \frac{\sigstar^2}{\stpsz \strongconvex} \niters^2.
\end{align*}
It suffices to bound the initial condition $\Exs\vecnorm{v_{\burnin}}{2}^2$. Recall that
$ v_{\burnin}	=	\frac{1}{\burnin} \sum_{s = 1}^{\burnin} \nabla f(\theta_0; \xi_s)$,
which is average of $\mathrm{i.i.d.}$ random vectors. It immediately follows from Assumptions~\ref{assu_noisethetastar} and~\ref{assu_smoothnoise} that:
\begin{align*}
    \Exs\vecnorm{v_{\burnin}}{2}^2 \leq \vecnorm{\nabla F(\theta_0)}{2}^2 + \frac{2 \sigstar^2}{\burnin} + \frac{2 \sglip^2 \vecnorm{\nabla F(\theta_0)}{2}^2 }{\strongconvex^2 \burnin}.
\end{align*}
Putting them together, we complete the proof of this lemma.

\pb\subsubsection{Proof of Lemma~\ref{lemm_ztsharpbdd}}\label{sec_ztsharpbdd}
Recalling that the recursive update rule of $z_t$ reveals an underlying martingale structure
\begin{align*}
t z_t
=
(t - 1) z_{t - 1} + (t - 1) (\noise_t(\theta_{t - 1}) - \noise_t(\theta_{t - 2})) + \noise_t(\theta_{t - 1})
.
\end{align*}
Adding and subtracting the $\noise_t(\thetastar)$ term in the above display we express the noise increment as
\begin{align*}
t z_t - (t - 1) z_{t - 1} 
= \noise_t(\thetastar) +  \underbrace{
(t - 1) (\noise_t(\theta_{t - 1}) - \noise_t(\theta_{t - 2})) + \noise_t(\theta_{t - 1}) - \noise_t(\thetastar)
}_{=:~ \remnoise_t }.
\end{align*}
In words, the increment of $t z_t$ splits into two parts: the additive part $\noise_t(\thetastar)$ and the multiplicative part $\remnoise_t$.
Taking expectation on the squared norm in above and using the property of square-integrable martingales, we have via further expanding the square on the right hand
\begin{align*}
t^2 \Exs\vecnorm{z_t}{2}^2
-
(t - 1)^2 \Exs\vecnorm{z_{t - 1}}{2}^2
&=
\Exs\vecnorm{\noise_t(\thetastar) + \remnoise_t}{2}^2 =
\Exs\vecnorm{\noise_t(\thetastar)}{2}^2
+
\Exs\vecnorm{\remnoise_t}{2}^2
+
2 \Exs\inprod{ \noise_t(\thetastar)}{\remnoise_t}
.
\end{align*}
Telescoping the above equality for $t=\burnin+1,\dots,T$ gives
\begin{align*}
\niters^2 \Exs\vecnorm{z_\niters}{2}^2
-
\burnin^2 \Exs\vecnorm{ z_{\burnin} }{2}^2
=
\sum_{t = \burnin + 1}^{\niters} \Exs\vecnorm{ \noise_t(\thetastar)}{2}^2
+
\sum_{t = \burnin + 1}^{\niters} \Exs\vecnorm{\remnoise_t}{2}^2
+
2 \sum_{t = \burnin + 1}^{\niters} \Exs \inprod{ \noise_t(\thetastar)}{ \remnoise_t}.
\end{align*}
By Assumption \ref{assu_noise_bdd_holder}, we have $\Exs\vecnorm{\noise_t(\thetastar)}{2}^2 = \sigma_*^2$. 
For the additional noise $\remnoise_t$, Young's inequality leads to the bound
\begin{align}
\Exs\vecnorm{\remnoise_t}{2}^2&\leq
\left(
\sglip (t - 1)\sqrt{\Exs\vecnorm{ \theta_{t - 1} - \theta_{t - 2}}{2}^2}
+
\sglip \sqrt{\Exs\vecnorm{\theta_{t - 1} - \thetastar}{2}^2}
\right)^2\nonumber
%%%%%%%%
\\&\leq
\left(
\stpsz\sglip (t - 1) \sqrt{\Exs\vecnorm{v_{t - 1}}{2}^2}
+
\frac{\sglip}{\strongconvex} \sqrt{\Exs \left\|\nabla F(\theta_{t - 1}) \right\|_2^2}
\right)^2\nonumber
%%%%%%%%
\\&\le
2\stpsz^2\sglip^2 (t - 1)^2 \Exs\vecnorm{v_{t - 1}}{2}^2
+ 
\frac{2\sglip^2}{\strongconvex^2} \Exs\vecnorm{\nabla F(\theta_{t - 1}) }{2}^2
.\label{eq:remnoise-bound}
\end{align}
It remains to bound the summation of the cross term. Observing that:
\begin{align*}
\sum_{t = \burnin + 1}^{\niters} \Exs
\inprod{\noise_t(\thetastar)}{ \remnoise_t} &=
\sum_{t = \burnin + 1}^{\niters} \Exs \inprod{
\noise_t(\thetastar)}{ \noise_t(\theta_{t - 1}) - \noise_t(\thetastar)
+ (t - 1) \left(\noise_t(\theta_{t - 1}) - \noise_t(\theta_{t - 2})\right)}
\\
&= \sum_{t = \burnin + 1}^{\niters} \Big\{
t \Exs \inprod{
\noise_t(\thetastar)}{ \noise_t(\theta_{t - 1}) - \noise_t(\thetastar)}
- 
(t - 1) \Exs \inprod{
\noise_t(\thetastar) }{ \noise_t(\theta_{t - 2}) - \noise_t(\thetastar)} \Big\}.
\end{align*}
Since the random samples $(\xi_t)_{t \geq 1}$ are $\mathrm{i.i.d.}$ and the iterate $\theta_{t - 2}$ is independent of the sample $\xi_{t - 1}$, we have that
\begin{align*}
\Exs \inprod{\noise_t(\thetastar) }{ \noise_t(\theta_{t - 2}) - \noise_t(\thetastar)} 
= 
\Exs \inprod{\noise_{t - 1}(\thetastar) }{ \noise_{t - 1} (\theta_{t - 2}) - \noise_{t - 1} (\thetastar)}
.
\end{align*}
Consequently, we can re-write the quantity of interests as a telescope sum, leading to the following identity:
\begin{align*}
\sum_{t = \burnin + 1}^{\niters} \Exs
\inprod{\noise_t(\thetastar)}{ \remnoise_t} &= \sum_{t = \burnin + 1}^{\niters} \Big\{
t \Exs \inprod{
\noise_t(\thetastar)}{ \noise_t(\theta_{t - 1}) - \noise_t(\thetastar)}
- 
(t - 1) \Exs \inprod{
\noise_t(\thetastar) }{ \noise_t(\theta_{t - 2}) - \noise_t(\thetastar)} \Big\}\\
%%%%%%%%
&=
\sum_{t = \burnin + 1}^{\niters} \Big\{
t \Exs \inprod{
\noise_t(\thetastar)}{ \noise_t(\theta_{t - 1}) - \noise_t(\thetastar)}
- 
(t - 1) \Exs \inprod{
\noise_{t - 1} (\thetastar) }{ \noise_{t - 1} (\theta_{t - 2}) - \noise_{t - 1} (\thetastar)} \Big\}\\
%%%%%%%%
&=
\niters \cdot  \Exs \inprod{
\noise_\niters (\thetastar)}{ \noise_\niters(\theta_{\niters - 1}) - \noise_\niters(\thetastar)} - \burnin \cdot \Exs \inprod{
\noise_{\burnin} (\thetastar)}{ \noise_{\burnin}(\theta_{\burnin - 1}) - \noise_{\burnin}(\thetastar)
}
.
\end{align*}
In order to bound the inner product terms, we invoke the Cauchy-Schwartz inequality and Assumption~\ref{assu_smoothnoise}. For each $t \geq \burnin$, we have that
\begin{multline}
\abss{t \cdot \Exs \inprod{
\noise_t (\thetastar)}{ \noise_t(\theta_{t - 1}) - \noise_t(\thetastar)} } 
\leq 
t \cdot \sqrt{\Exs\vecnorm{\noise_t(\thetastar)}{2}^2} \cdot \sqrt{\Exs\vecnorm{\noise_t(\theta_{t - 1}) - \noise_t(\thetastar)}{2}^2 }
\\ \leq 
t \sigstar \sglip \sqrt{\Exs\vecnorm{\theta_{t - 1} - \thetastar}{2}^2 } 
\leq 
\frac{t \sigstar \sglip}{\strongconvex} \sqrt{\Exs\vecnorm{\nabla F (\theta_{t - 1})}{2}^2}
.\label{eq:zt-cross-term-bound-after-telescope}
\end{multline}
Plugging $t = \burnin$ and $t = \niters$ into Eq.~\eqref{eq:zt-cross-term-bound-after-telescope} separately and combining with Eq.~\eqref{eq:remnoise-bound}, we have that
\begin{align}
\niters^2 \Exs\vecnorm{z_\niters}{2}^2 
&\leq 
\burnin^2 \Exs\vecnorm{ z_{\burnin} }{2}^2 
+ 
(\niters - \burnin) \sigstar^2 
+
2\stpsz^2\sglip^2\sum_{t = \burnin + 1}^{\niters} (t - 1)^2 \Exs\vecnorm{v_{t - 1}}{2}^2
+
\frac{2\sglip^2}{\strongconvex^2}\sum_{t = \burnin + 1}^{\niters} \Exs\vecnorm{\nabla F(\theta_{t - 1}) }{2}^2 
\nonumber\\&\qquad
+
\frac{2\sglip\sigstar}{\strongconvex}\left(
\niters\sqrt{\Exs\|\naF{\theta_{\niters-1}}\|_2^2}
+
\burnin \norm{\naF{\theta_0}}
\right)
.\label{eq:zt-sharp-lemma-ready-form}
\end{align}

The above bound involves the second moments of the vectors $v_t$ and $\nabla F (\theta_t)$. We recall the following bounds from Theorem~\ref{theo_finitebdd_single} and Lemma~\ref{lemm_vtsharpbdd}, for each $t \geq \burnin$:
\begin{align*}
&
\Exs \normb{v_t} 
\leq 
\frac{c\sigstar^2}{\stpsz\strongconvex t^2} 
+
\frac{c}{\stpsz^4 \strongconvex^4 t^4} \normb{\nabla F(\theta_0)}
, \quad \mbox{and}
\\& 
\Exs\vecnorm{\nabla F(\theta_t) }{2}^2 
\leq 
\frac{c \vecnorm{\nabla F(\theta_0)}{2}^2}{\stpsz^2\strongconvex^2t^2} +
\frac{c \sigstar^2}{t}
.
\end{align*}
Substituting these bounds to Eq.~\eqref{eq:zt-sharp-lemma-ready-form}, we note that
\begin{align*}
&
\sum_{t = \burnin + 1}^{\niters} (t - 1)^2 \Exs\vecnorm{v_{t - 1}}{2}^2 
\leq 
\frac{ c \sigstar^2 \niters}{\stpsz\strongconvex} 
+ 
\frac{ c \vecnorm{\nabla F(\theta_0)}{2}^2}{\stpsz^4\strongconvex^4\burnin}
, \quad \mbox{and}\\&
\sum_{t = \burnin + 1}^{\niters} \Exs\vecnorm{\nabla F(\theta_{t - 1}) }{2}^2 
\leq 
c \sigstar^2 \log \Big( \frac{\niters}{\burnin} \Big) 
+ 
\frac{c \vecnorm{\nabla F(\theta_0)}{2}^2}{\stpsz^2\strongconvex^2\burnin}
.
\end{align*}
Finally, for the burn-in period we note that $
z_{\burnin}	=	\sum_{s=1}^{\burnin} \noise_s(\theta_0)
$ is a sum of $\burnin$ i.i.d.~random vectors, and hence the following estimation holds
\begin{align*}
\lefteqn{
\burnin^2 \Exs\vecnorm{z_{\burnin}}{2}^2 
=
\burnin \Exs\vecnorm{\noise_1(\theta_0)}{2}^2 
}\\&\le%%%%
\burnin \Exs\vecnorm{\noise_1(\theta_0) - \noise_1(\thetastar)}{2}^2 
+ 
2\burnin \sqrt{\Exs\vecnorm{\noise_1(\thetastar)}{2}^2 \Exs\vecnorm{\noise_1(\theta_0) - \noise_1(\thetastar)}{2}^2}
+ 
\burnin \Exs\vecnorm{\noise_1(\thetastar)}{2}^2 
%%%%%%%%
\\&\le
\burnin \sglip^2 \vecnorm{\theta_0 - \thetastar}{2}^2 
+ 
2\burnin \sglip\sigstar \vecnorm{\theta_0 - \thetastar}{2}
+ 
\burnin \sigstar^2
\le%%%%
\frac{\burnin \sglip^2}{\strongconvex^2} \vecnorm{\nabla F(\theta_0)}{2}^2 
+
\frac{2\burnin \sglip \sigstar}{\strongconvex} \vecnorm{\nabla F(\theta_0)}{2} 
+ 
\burnin \sigstar^2
.
\end{align*}
Some algebra yields \eqref{ztsharpbdd} and hence the whole Lemma \ref{lemm_ztsharpbdd}.

\pb\subsubsection{Proof of Lemma~\ref{lemm_cross-term-bound}}\label{subsubsec:proof-cross-term-bound}
First, by Cauchy-Schwartz inequality, we can easily observe that:
\begin{align*}
\abss{\Exs \inprod{v_\niters}{z_\niters}}
\le
\sqrt{\Exs\vecnorm{v_\niters}{2}^2} \cdot \sqrt{\Exs\vecnorm{z_\niters}{2}^2} 
\le
c \Big( \frac{\sigstar}{\niters \sqrt{ \stpsz \strongconvex}} 
+ 
\frac{\vecnorm{\nabla F(\theta_0)}{2}}{\stpsz^2 \strongconvex^2 \niters^2} \Big) \cdot \Big( \frac{\sigstar}{\sqrt{\niters}} 
+ 
\frac{\vecnorm{\nabla F(\theta_0)}{2}}{\stpsz \strongconvex \niters} \Big)
.
\end{align*}
So for $\niters \leq c \burnin \log \burnin$, the conclusion of this lemma is automatically satisfied. For the rest of this section, we assume that $\frac{\niters}{\log \niters} > c \burnin$ for some universal constant $c > 0$.

The proof requires some bounds on the fourth moment of the stochastic process defined by the algorithm. 
In particular, we need the following two lemmas. 
The first lemma is analogous to the bound in Theorem \ref{theo_finitebdd_single}:
\begin{lemma}[Higher-order-moment bound on $\naF{\theta_{t - 1}}$]\label{lemm_finitebddHolder}
Suppose Assumptions~\ref{assu_StrcvxSmooth},~\ref{assu_noise_bdd_holder} and~\ref{assu_noise_smooth_holder} hold. 
Let the step-size $\stpsz \le \frac{1}{56\smoothness} \wedge \frac{\strongconvex}{64\sgliptild^2}$ and the burn-in time $\burnin = \left\lceil\frac{24}{\stpsz\strongconvex}\right\rceil$. 
Then for any $\niters \geq \burnin$, the estimator $\theta_{\niters}$ produced by the \ROOTSGD algorithm satisfies the bound 
\begin{align}\label{finitebdd_holder}
\para{\Exs \norm{\naF{\theta_{\niters}}}^4}^{1/2}
\le
\frac{140\sigstartild^2}{T + 1}
+
\frac{60 \normb{\naF{\theta_0}}}{\stpsz^2 \strongconvex^2 (T + 1)^2}
.
\end{align}
\end{lemma}
\noindent Proof can be found in \S\ref{sec_proof,lemm_finitebddHolder}.

We also need a lemma on the fourth-moment bound of $v_t$, analogous to Lemma~\ref{lemm_vtsharpbdd}:
\begin{lemma}[sharp higher-order-moment bound on $v_t$]\label{lemm_vtsharpbdd_holder}
Under the setting of Proposition \ref{theo_rootHolder_single} we have the following bound for $\niters\ge \burnin+1$
\begin{align}\label{vtsharpbdd_holder}
\sqrt{\Exs \normd{v_\niters}} 
	\le
\frac{4484 \sigstartild^2}{\stpsz\strongconvex \niters^2} 
+
\frac{1359375}{\stpsz^4 \strongconvex^4 \niters^4} \normb{\nabla F(\theta_0)}
.
\end{align}
\end{lemma}
\noindent Proof can be found in \S\ref{sec_proof,lemm_vtsharpbdd_holder}.

Taking these two lemmas as given, we proceed with the proof.
Following the two-time-scale intuition discussed in Section~\ref{sec_hidden}, the process $v_t$ moves faster than the averaging process $z_t$. 
Therefore, it is reasonable to expect the correlation between $v_t$ and $z_{t - \Twindow}$ to be small, for sufficiently large time window $\Twindow > 0$. 
For the rest of this section, we choose the window size:
\begin{align}
\Twindow 
= 
\frac{c}{\strongconvex \stpsz} \log \niters
, \quad\mbox{for some universal constant $c > 0$.}\label{eq:window-size-in-cross-term-proof}
\end{align}
Since we have assumed without loss of generality that $\tfrac{\niters}{\log \niters} > c \burnin =\frac{24 c}{\stpsz\strongconvex}$, the window size guarantees the relation $\niters - \Twindow > \niters / 2$.

We subtract off a $(t - \Twindow) z_{t - \Twindow}$ term the $t z_t$ expression above, and decompose the absolute value of the cross term $|\Exs \binprod{ v_t}{tz_t}|$ as:
\begin{align}\label{eq:decomp_holder}
|\Exs \binprod{tz_t}{v_t} |
&\leq 
(t - \Twindow) \underbrace{\left|\Exs \binprod{ z_{t - \Twindow}}{v_t} \right|}_{=: I_1} 
+
\underbrace{\left| \Exs \binprod{tz_t - (t - \Twindow) z_{t - \Twindow}}{v_t} \right|}_{=: I_2}
.
\end{align} 
For bounding the term $I_2$, we make use of the recursive rule of $tz_t$ to obtain the bound
\begin{align*}
&
\Exs \normb{tz_t - (t - \Twindow) z_{t - \Twindow}} 
= 
\Exs \normb {\sum_{s = t - \Twindow+1}^t\Big\{ 
(s - 1) (\noise_s(\theta_{s - 1}) - \noise_s(\theta_{s - 2})) + \noise_s(\theta_{s - 1}) - \noise_s(\theta^*) + \noise_s(\theta^*)
\Big\}}
\\&\leq
\sum_{s = t - \Twindow + 1}^t \Big\{ (s - 1)^2 \stpsz^2 \sglip^2 \Exs \normb{v_{s - 1}} + \frac{\sglip^2}{\strongconvex^2} \Exs \normb{\naF{\theta_{s - 1}}} + \sigstar^2  \Big\}
\leq 
\Twindow \cdot \left( 
\sigstar^2 + \frac{\normb{\naF{\theta_0}}}{\stpsz^3 \strongconvex^3 t^2}
\right)
.
\end{align*}
Consequently, we have the bound
\begin{multline*}
\left| \Exs \binprod{ tz_t - (t - \Twindow) z_{t - \Twindow}}{v_t} \right|
\leq 
\sqrt{\Exs \normb{tz_t - (t - \Twindow) z_{t - \Twindow}}} \cdot 
\sqrt{\Exs \normb{v_t}}\\
\leq 
c \sqrt{\Twindow} \left(\frac{\sigstar}{\sqrt{\stpsz\strongconvex} t} + \frac{\norm{\naF{\theta_0}}}{\stpsz^2 \strongconvex^2 t^2} \right) \left(\sigstar + \frac{\norm{\naF{\theta_0}}}{ \stpsz^{3/2} \strongconvex^{3/2} t} \right) 
\leq 
c \left(\frac{\sigstar^2}{\stpsz \strongconvex t } + \frac{\normb{\naF{\theta_0}}}{\stpsz^{4}\strongconvex^{4} t^{3} } \right) \sqrt{\log t}
.
\end{multline*}

The bound for the term $I_1$ in the decomposition~\eqref{eq:decomp_holder} is given by the following analysis: law of iterated expectations gives
\begin{align}\label{eq:cauchy_decomp}
\left|\Exs \binprod{ z_{t - \Twindow}}{v_t } \right| 
=
\left| \Exs \binprod{z_{t - \Twindow}}{\Exs (v_t \mid \cF_{t - \Twindow}) } \right| 
\leq 
\sqrt{\Exs \norm{z_{t - \Twindow}}^2} \cdot \sqrt{\Exs \norm{\Exs(v_t \mid \cF_{t - \Twindow})}^2}
,
\end{align}
where the last inequality comes from applying the Cauchy-Schwarz inequality.

The second moment for $z_{t - \Twindow}$ is relatively easy to estimate using Lemma~\ref{lemm_ztsharpbdd}. 
It suffices to study the conditional expectation $\Exs(v_t \mid \cF_{t - \Twindow})$. We claim the following bound:
\begin{align}
\sqrt{\Exs\vecnorm{ \Exs(v_t \mid \cF_{t - \Twindow}) }{2}^2}
\leq 
\frac{c \holderconst}{\strongconvex} \left(\frac{\sigstartild}{\sqrt{\stpsz\strongconvex} t} 
+ 
\frac{\norm{\naF{\theta_0}}}{\stpsz^2 \strongconvex^2 t^2} \right) \left( \frac{\sigstartild}{\strongconvex \sqrt{t}} 
+ 
\frac{\norm{\naF{\theta_0}}}{\stpsz\strongconvex^2 t}\right) 
.\label{eq:conditional-expectation-bound}
\end{align}
We prove this inequality at the end of this section. 
Taking this bound as given, we now proceed with the proof for Lemma~\ref{lemm_cross-term-bound}.

Bringing this back to the inequality~\eqref{eq:cauchy_decomp} and by utilizing the $z_t$ bound by Lemma~\ref{lemm_ztsharpbdd}, we have
\begin{align*}
\Exs \normb{z_t} 
\le
C \left(\frac{\sigstar^2}{t} + \frac{\normb{\naF{\theta_0}}}{\stpsz^2 \strongconvex^2 t^2} \right)
,
\end{align*}
and thus
\begin{align*}
\lefteqn{
\left| \Exs \binprod{ z_{t - \Twindow}}{v_t } \right| 
\le
\sqrt{\Exs \normb{z_{t - \Twindow}}} \cdot 
\sqrt{\Exs \normb{\Exs \left(v_t \mid \cF_{t - \Twindow} \right)}}
}\notag \\&\leq 
\frac{c \holderconst}{\strongconvex} \left( 
\frac{\sigstar}{\sqrt{t}} + \frac{\norm{\naF{\theta_0}}}{\stpsz\strongconvex t}
\right)\left(\frac{\sigstartild}{\sqrt{\stpsz\strongconvex} t} + \frac{\normb{\naF{\theta_0}}}{\stpsz^2 \strongconvex^2 t^2} \right) \left( \frac{\sigstartild}{\strongconvex\sqrt{t}} + \frac{\normb{\naF{\theta_0}}}{\stpsz\strongconvex^2 t}\right) 
\notag \\&\leq 
\frac{c \holderconst}{\strongconvex} 
\left( 
\frac{\sigstartild^{3}}{\stpsz^{1/2}\strongconvex^{3/2} t^{2}} 
	+
\frac{ \norm{\naF{\theta_0}}^{3} }{\stpsz^{7/2} \strongconvex^{9/2} t^{7/2}}
\right)
.
\end{align*}
Combining the bounds for $I_1$ and $I_2$ together, we estimate the cross term  as:
\begin{equation}\label{eq:last_two}
\begin{aligned}
\lefteqn{
\left| \Exs \binprod{ tz_t}{v_t } \right|
}
\\&\leq 
c (t - \Twindow)\frac{\holderconst}{\strongconvex} 
\left( 
\frac{\sigstartild^{3}}{\stpsz^{1/2}\strongconvex^{3/2} t^{2}} 
	+
\frac{ \norm{\naF{\theta_0}}^{3} }{\stpsz^{7/2} \strongconvex^{9/2} t^{7/2}}
\right)
+
c \left( 
\frac{\sigstar^2}{\stpsz \strongconvex t} + \frac{\normb{\naF{\theta_0}}}{\stpsz^{4}\strongconvex^{4} t^{3}}
\right)\sqrt{\log t}
. 
\end{aligned}\end{equation}
We conclude by dividing both sides of Eq.~\eqref{eq:last_two} by $\niters$ and arrive at the following bound:
\begin{align*}
\lefteqn{
\left| \Exs \binprod{v_{\niters}}{z_{\niters}} \right|
}\notag \\&\leq 
c\left(\frac{\sigstar^2}{\stpsz \strongconvex \niters^2}
+
\frac{ \normb{\naF{\theta_0}}}{\stpsz^4 \strongconvex^4 \niters^4} \right)\sqrt{\log \niters}
+
c  \holderconst
\left( 
\frac{\sigstartild^{3} }{\stpsz^{1/2}\strongconvex^{5/2}\niters^{2}} 
+
\frac{ \norm{\naF{\theta_0}}^{3} }{\stpsz^{7/2} \strongconvex^{11/2}\niters^{7/2}} 
\right)
.
\end{align*}
This finishes our bound on the cross term and conclude Lemma~\ref{lemm_cross-term-bound}.

\paragraph{Proof of Eq~\eqref{eq:conditional-expectation-bound}:}
We note the following expansion:
\begin{align*}
\naF{\theta_{t - 1}} - \naF{\theta_{t - 2}} 
=
\int_0^1 \nabla^2 F\left(\lambda \theta_{t - 2} + (1 - \lambda) \theta_{t - 1}\right) (\theta_{t - 1} - \theta_{t - 2}) d \lambda
,
\end{align*}
which leads to the following bound under the Lipschitz continuity condition for the Hessians (Assumption \ref{assu_holderhessian}):
\begin{align}\label{eq:hessian_approx}
&\lefteqn{
\norm{\naF{\theta_{t - 1}} - \naF{\theta_{t - 2}} - \nabla^2 F(\thetastar) (\theta_{t - 1} - \theta_{t - 2})}
}\notag \\&=
\int_0^1\norm{ \left( \nabla^2 F\left(\lambda \theta_{t - 2} + (1 - \lambda) \theta_{t - 1}\right) - \nabla^2 F(\thetastar)\right)   (\theta_{t - 1} - \theta_{t - 2})} d \lambda
\notag \\&\leq 
\stpsz \holderconst \norm{v_{t - 1}} \int_0^1 
\norm{\lambda (\theta_{t - 2} - \thetastar) + (1 - \lambda) (\theta_{t - 1} - \thetastar)} d \lambda
\notag \\&\leq 
\stpsz \holderconst \norm{v_{t - 1}} \cdot \max\left( 
\norm{\theta_{t - 1} - \thetastar},	\norm{\theta_{t - 2} - \thetastar}
\right)
.
\end{align} 
Since $\hessianstar = \nabla^2 F(\thetastar)$ we have
\begin{align}
\lefteqn{
t \norm{\Exs(v_t \mid \cF_{t - \Twindow})}
}\notag \\&= 
\norm{\Exs \left((t - 1)\left(v_{t - 1} + \naF{\theta_{t - 1}} - \naF{\theta_{t - 2}} \right) 
+ \nabla F(\theta_{t - 1}) \mid \cF_{t - \Twindow} \right) } 
\notag \\&= 
\left\|
\Exs \left((t - 1)\left(v_{t - 1} + H^* (\theta_{t - 1} - \theta_{t - 2}) \right.\right.\right.
\notag \\&~\quad 
+
\left.\left.\left. \naF{\theta_{t - 1}} - \naF{\theta_{t - 2}} - H^* (\theta_{t - 1} - \theta_{t - 2}) \right)
+ \nabla F(\theta_{t - 1}) \mid \cF_{t - \Twindow} \right)
\right\|_2
    \notag \\&\leq 
\vecnorm{ \Exs \left((t - 1)\left(v_{t - 1} + H^* (\theta_{t - 1} - \theta_{t - 2}) \right) \mid \cF_{t - \Twindow} \right)  }{2}
+
\vecnorm{ \Exs (v_t \mid \cF_{t - \Twindow}) }{2}
\notag \\&~\quad 
+
\vecnorm{ \Exs\left( (t - 1) \left( \naF{\theta_{t - 1}} - \naF{\theta_{t - 2}} - H^* (\theta_{t - 1} - \theta_{t - 2}) \right) \mid \cF_{t - \Twindow} \right)  }{2}
. \label{eq:cross_I1_decomp}
\end{align}
Further by rearranging the terms, and dividing both sides by $(t - 1)$, we obtain 
\begin{align*}
\lefteqn{
\vecnorm{ \Exs(v_t \mid \cF_{t - \Twindow}) }{2} 
}
\\&\le
\vecnorm{ \Exs  \left( v_{t - 1} + H^*(\theta_{t - 1} - \theta_{t - 2}) \mid \cF_{t - \Twindow} \right) }{2}
%\notag \\&\quad 
+ \vecnorm{ 
\Exs \left( \nabla F(\theta_{t - 1}) - \nabla F(\theta_{t - 2}) - H^* (\theta_{t - 1} - \theta_{t - 2}) \right) \mid \cF_{t - \Twindow}  
}{2}
\\&\le
(1 - \stpsz \strongconvex) \vecnorm{ \Exs \left( v_{t - 1} \mid \cF_{t - \Twindow} \right) }{2}
%\notag \\&\quad
+ 
\stpsz \holderconst \Exs \left( \vecnorm{v_{t - 1}}{2} \cdot \max \left( 
\vecnorm{\theta_{t - 1} - \theta^* }{2}
,  
\vecnorm{\theta_{t - 2} - \theta^* }{2} 
\right) \mid \cF_{t - \Twindow} \right)
,
\end{align*}
where in the last inequality we apply the result in Eq.~\eqref{eq:hessian_approx}.
Next by calculating the second moment of both the RHS and the LHS of the above quantity and the H\"{o}lder's inequality, we have
\begin{align*}
&
\sqrt{\Exs\vecnorm{\Exs(v_t \mid \cF_{t - \Twindow}) }{2}^2}
\\&\leq
(1 - \stpsz \strongconvex) \sqrt{\Exs\vecnorm{ \Exs(v_{t - 1} \mid \cF_{t - \Twindow}) }{2}^2} 
+
\stpsz \holderconst \sqrt{\Exs\vecnorm{\Exs\left(
\norm{v_{t - 1}}^2  \cdot \left( 
\norm{\theta_{t - 1} - \thetastar}^{2} + \norm{\theta_{t - 2} - \thetastar}^{2} 
\right) \mid \cF_{t - \Twindow}\right)}{2}^2
}
\notag \\&\leq 
(1 - \stpsz \strongconvex) \sqrt{\Exs\vecnorm{ \Exs(v_{t - 1} \mid \cF_{t - \Twindow}) }{2}^2}
%\notag \\&~\quad 
+
\stpsz L_{1}  \left( \Exs \norm{v_{t - 1}}^{4} \right)^{1/4}
\Big\{ 
\big( \Exs \norm{\theta_{t - 1} - \thetastar}^{4}\big)^{1/4} + \big( \Exs \norm{\theta_{t - 2} - \thetastar}^{4}\big)^{1/4}
\Big\}
.
\end{align*}
Recursively applying the above inequality from $t - \Twindow$ to $t$ and we have that 
\begin{equation}\label{eq:condition_vt_bound}
\begin{aligned}
\lefteqn{
\sqrt{\Exs\vecnorm{ \Exs(v_t \mid \cF_{t - \Twindow}) }{2}^2}
}\\&\leq 
(1 - \stpsz \strongconvex)^{\Twindow}  \Exs\vecnorm{ v_{t - \Twindow} }{2}^2
+ 
\frac{\holderconst}{\strongconvex} \max_{t - \Twindow \leq s \leq t} \left( \Exs \norm{v_{s}}^{4} \right)^{1/4} \cdot \max_{t - \Twindow \leq s \leq t} \big( \Exs \norm{\theta_{t - 2} - \thetastar}^{4}\big)^{1/4}
.
\end{aligned}\end{equation}

We recall from Lemmas~\ref{lemm_finitebddHolder} and~\ref{lemm_vtsharpbdd_holder} the following
$$\begin{aligned}
&
\left( \Exs \norm{v_{\niters}}^{4}\right)^{1/2} 
\le
C \left( 
\frac{\tilde{\sigstar}^2}{\stpsz \strongconvex \niters^2} 
+
\frac{ \normb{\naF{\theta_0}}}{\stpsz^4\strongconvex^4 \niters^4}
\right)
,&\text{and}\\&
\left(\Exs \norm{\naF{\theta_{\niters-1}}}^{4}\right)^{1/2}
\le
C \left(\frac{\tilde{\sigstar}^2}{\niters} + \frac{\normb{\naF{\theta_0}}}{ \stpsz^2 \strongconvex^2 \niters^2} \right)
.
&
\end{aligned}$$
Bringing this into Eq.~\eqref{eq:condition_vt_bound} and we have that 
\begin{multline*}
\sqrt{\Exs\vecnorm{ \Exs(v_t \mid \cF_{t - \Twindow}) }{2}^2 }
\leq 
(1 - \stpsz \strongconvex)^{\Twindow} \Exs\vecnorm{ v_{t - \Twindow} }{2}^2 
\notag \\ 
+
\frac{c \holderconst}{\strongconvex} \left(\frac{\sigstartild}{\sqrt{\stpsz\strongconvex} (t - \Twindow)} + \frac{\norm{\naF{\theta_0}}}{\stpsz^2 \strongconvex^2 (t - \Twindow)^2} \right) \left( \frac{\sigstartild}{\strongconvex \sqrt{t - \Twindow}} + \frac{\norm{\naF{\theta_0}}}{\stpsz\strongconvex^2 (t - \Twindow)}\right) 
.
\end{multline*}
Substituting with the window size $\Twindow$ defined in Eq~\eqref{eq:window-size-in-cross-term-proof}, the above inequality reduces as follows:
\begin{align*}
\sqrt{\Exs\vecnorm{ \Exs(v_t \mid \cF_{t - \Twindow}) }{2}^2}
\leq 
\frac{c \holderconst}{\strongconvex} \left(\frac{\sigstartild}{\sqrt{\stpsz\strongconvex} t} + \frac{\norm{\naF{\theta_0}}}{\stpsz^2 \strongconvex^2 t^2} \right) \left( \frac{\sigstartild}{\strongconvex \sqrt{t}} + \frac{\norm{\naF{\theta_0}}}{\stpsz\strongconvex^2 t}\right)
.
\end{align*}

\pb\subsection{Proof of Theorem \ref{theo_rootHolder_multi}}\label{sec_proof,theo_rootHolder_multi}
Utilizing the intermediate Proposition \ref{theo_rootHolder_single} in \S\ref{sec_proof,theo_rootHolder_single}, we now aim to improve the dependency on initialization and turn to the proof of our multi-epoch nonasymptotic result. 
Invoking Eq.~\eqref{finitebdd_holder} in Lemma~\ref{lemm_finitebddHolder}, we obtain for $b = 1, 2, \cdots, \Epoch$ the bound for $\inneriters \geq c \burnin$:
\begin{align*}
&\Exs \normb{\naF{\theta_0^{(b + 1)}}} 
\leq 
\frac1{e^2} \Exs \normb{\naF{\theta_0^{(b)}}} 
+
\frac{c \sigstar^2}{\inneriters}
,&\text{and}\\&
\sqrt{\Exs \normd{\naF{\theta_0^{(b + 1)}}}}
\leq 
\frac1{e^2} \sqrt{\Exs \normd{\naF{\theta_0^{(b)}}}} + \frac{c \sigstartild^2}{\inneriters}
,&
\end{align*}
where our setting of $\inneriters$ gives a discount factor of $1/e^2$.
Solving the recursion, we arrive at the bound:
\begin{align*}
&\Exs \normb{\naF{\theta_0^{(\Epoch + 1)}}} 
\leq 
\frac{c \sigstar^2}{\inneriters} + e^{-2\Epoch}  \normb{\naF{\theta_0}}
,&\text{and}\\&
\sqrt{\Exs \normd{\naF{\theta_0^{(\Epoch + 1)}}}}
\leq 
\frac{c \sigstartild^2}{\inneriters} +   e^{-2\Epoch} \sqrt{\Exs \normd{\naF{\theta_0}}}
.
\end{align*}
Our take is $\Epoch \geq \log \frac{\inneriters \sqrt{\Exs \normd{\naF{\theta_0}}}}{c\sigstar^2}$ such that $e^{-2\Epoch} \Exs \normb{\naF{\theta_0}} \leq \frac{\sigstar^2}{\inneriters}$ and $e^{-2\Epoch} \sqrt{\Exs \normd{\naF{\theta_0}}} \leq \frac{\sigstartild^2}{\inneriters}$ both hold.
Finally, we have 
\begin{align*}
&\Exs \normb{\naF{\theta_0^{(\Epoch + 1)}}} 
\le
e^{-2\Epoch} \Exs \normb{\naF{\theta_0}} 
+
\frac{c \sigstar^2}{\inneriters}
\le
\frac{c' \sigstar^2}{\inneriters}
,&\text{and}
\\
&\sqrt{\Exs \normd{\naF{\theta_0^{(\Epoch + 1)}}}}
\le
e^{-2\Epoch} \sqrt{\Exs \normd{\naF{\theta_0}}} + \frac{c \sigstartild^2}{\inneriters}
\leq 
\frac{c'\sigstartild^2}{\inneriters}
,&\text{and}
\\ 
&\Exs \norm{\naF{\theta_0^{(\Epoch + 1)}}}^{3}
\le
\left(\Exs \normd{\naF{\theta_0^{(\Epoch + 1)}}} \right)^{\frac{3}{4}}
\leq
\frac{c' \sigstartild^{3}}{(\inneriters)^{3/2}}
,&\text{and}
\\ 
&\Exs \norm{\naF{\theta_0^{(\Epoch + 1)}}} 
\le
\frac{c \sigstar}{(\inneriters)^{1/2}}
,
\end{align*}
where constants $c,c'$ change from line to line.
Substituting this initial condition into the bound~\eqref{rootHolder_single}, we obtain the final bound:
\begin{align*}
&\lefteqn{
\Exs \normb{\naF{\theta_{\niters}^{(\Epoch + 1)}}} 
-
\frac{\sigstar^2}{\niters}
}
%%%%%%%%
\notag \\&\leq 
C \left(\frac{\sglip^2\stpsz}{\strongconvex}
+\frac{\sglip}{\strongconvex \niters^{1/2}}
+\frac{\log \niters}{\stpsz\strongconvex \niters}
+\frac{\log\left( \frac{\niters}{\burnin} \right)\sglip^2}{\strongconvex^2\niters}
\right)
\frac{\sigstar^2}{\niters}
%\\&\quad\,%%%%
	+
\frac{C\holderconst \sigstartild^{3}}{\stpsz^{1/2}\strongconvex^{5/2} \niters^{2}}
	\notag \\&\quad\,%%%%
+C \left(
\frac{\normb{\naF{\theta_0}}}{\stpsz^2 \strongconvex^2 \niters^2}
	+
\frac{\sglip\sigstar\norm{\naF{\theta_0}}}{\stpsz \strongconvex^2 \niters^2 }
%\\&\quad\,%%%%
\right)
	 +
\frac{C \holderconst \norm{\naF{\theta_0}}^{3} }{\stpsz^{7/2} \strongconvex^{11/2}\niters^{7/2}} 
\notag \\&\leq 
C \left(\frac{\sglip^2\stpsz}{\strongconvex}
+
\frac{\sglip}{\strongconvex \niters^{1/2}}
+
\frac{\log \niters}{\stpsz\strongconvex \niters}
+
\frac{\log\left( \frac{\niters}{\burnin} \right)\sglip^2}{\strongconvex^2\niters}
\right)
\frac{\sigstar^2}{\niters}
%\\&\quad\,%%%%
	+
\frac{C\holderconst \sigstartild^{3}}{\stpsz^{1/2}\strongconvex^{5/2} \niters^{2}}
%	\notag \\&\quad\,%%%%
+C \left(
\frac{ \sigstar^2}{\stpsz \strongconvex \niters^2}
	+
\frac{\holderconst \sigstartild^{3}}{\stpsz^{2} \strongconvex^{4}\niters^{7/2}}
\right)
\notag \\&\leq 
C \left(\frac{\sglip^2\stpsz}{\strongconvex}
+\frac{\sglip}{\strongconvex \niters^{1/2}}
	+
\frac{\log \niters}{\stpsz\strongconvex \niters}
+\frac{\log\left( \frac{\niters}{\burnin} \right)\sglip^2}{\strongconvex^2\niters}
\right)
\frac{\sigstar^2}{\niters} 
	+
\frac{C\holderconst \sigstartild^{3}}{\stpsz^{1/2}\strongconvex^{5/2} \niters^{2}}
,
\end{align*}
which proves the bound~\eqref{rootHolder_multi_raw}.

Finally, substituting $\niters$ by the final epoch length $\nsamples - \Epoch \inneriters$ and adopt similar reasoning as the previous one, we arrive at the conclusion:
\begin{align*}
\Exs \normb{\naF{\theta^{\text{\rm final}}_\nsamples}} 
-
\frac{\sigstar^2}{\nsamples}
	\le
C \left(\frac{\sglip^2\stpsz}{\strongconvex} + \frac{\sglip}{\strongconvex \nsamples^{1/2}}
+
\frac{\log \numobs}{\stpsz\strongconvex \nsamples}
+\frac{\log\left( \frac{\nsamples}{\burnin} \right)\sglip^2}{\strongconvex^2\nsamples}
\right)
\frac{\sigstar^2}{\nsamples} 
+
\frac{C\holderconst \sigstartild^{3}}{\stpsz^{1/2}\strongconvex^{5/2} \nsamples^{2}}
,
\end{align*}
which proves the bound~\eqref{rootHolder_multi_raw}.
Plugging $\stpsz$ as given by $\eta = \frac{c}{\sglip \nsamples^{1/2}}
			\wedge \frac{1}{4\smoothness}$ with $c = 0.49$, we have
\begin{align*}
\Exs \left\| \nabla F(\theta_\nsamples^{\text{\rm final}}) \right\|_2^2 
    &\leq 
C \left(\frac{\sglip^2\stpsz}{\strongconvex} + \frac{\sglip}{\strongconvex \nsamples^{1/2}}
+
\frac{\log \numobs}{\stpsz\strongconvex \nsamples}
+\frac{\log\left( \frac{\nsamples}{\burnin} \right)\sglip^2}{\strongconvex^2\nsamples}
\right)
\frac{\sigstar^2}{\nsamples} 
+
\frac{C\holderconst \sigstartild^{3}}{\stpsz^{1/2}\strongconvex^{5/2} \nsamples^{2}}
    \\&\leq 
C \left(
\log\left( \frac{e\nsamples}{\burnin} \right)
\left(\frac{\sglip}{\strongconvex\sqrt{\nsamples}}\right)
+
\frac{\smoothness}{\strongconvex \nsamples} 
\right)
\frac{\sigstar^2}{\nsamples} %%%%
+
\frac{
C (\sglip^{1/2} \nsamples^{1/4} + \smoothness^{1/2})
\holderconst \sigstartild^{3}
}{
\strongconvex^{5/2} \nsamples^{2}}
.
\end{align*}
This concludes~\eqref{rootHolder_multi_accu} and hence Theorem~\ref{theo_rootHolder_multi}.

\pb\subsection{Proof of Corollary~\ref{theo_other_bounds}}\label{subsec:proof-theo-other-bounds}

The proof consists of two parts: bounds on the mean-squared error $\Exs\vecnorm{\theta_T - \thetastar}{2}^2$ and bounds on the expected objective gap $\Exs \big[ F (\theta_T) - F (\thetastar) \big]$. Two technical lemmas are needed in the proofs for both cases.

The first lemma is analogous to Lemma~\ref{lemm_ztsharpbdd}, which provides a sharp bound on $\testMat z_t$ for any matrix $\testMat \in \real^{d \times d}$.

\begin{lemma}\label{lemm_ztsharpbdd_testmat}
Under settings of Theorem \ref{theo_finitebdd_single}, for any matrix $\testMat \in \real^{d \times d}$, the following bounds hold true for $\niters \ge \burnin + 1$:
\begin{multline}
\Exs  \vecnorm{\testMat z_\niters}{2}^2  
\leq 
\frac{1}{\niters} \mathrm{Tr} \left( \testMat \SigStar \testMat^\top \right) 
+
c \opnorm{\testMat}^2 \Big\{
\frac{\sglip^2\stpsz}{\strongconvex}
+
\frac{\sglip}{\strongconvex \sqrt{\niters}}
+
\frac{\log\left( \frac{\niters}{\burnin} \right)\sglip^2}{\strongconvex^2\niters}
\Big\} \frac{\sigstar^2}{\niters} 
\\
+  c\opnorm{\testMat}^2 \Big\{\frac{\sglip \sigstar}{\mu} \cdot \frac{\burnin}{\niters^2} \vecnorm{\nabla F(\theta_0)}{2}
+ \frac{\burnin^2}{ \niters^2}\normb{\naF{\theta_0}} \Big\}
.
\end{multline}
for some universal constant $c > 0$.
\end{lemma}

The second lemma is analogous to Lemma~\ref{lemm_cross-term-bound}, and provides sharp bound on the cross term $\Exs \inprod{\testMat z_t}{\testMat v_t}$.

\begin{lemma}\label{lemm_cross-term-bound-testmat}
    Under settings of Theorem \ref{theo_finitebdd_single}, we have the following bound for any $\niters \geq \burnin + 1$:
\begin{multline}
\left| \Exs \binprod{\testMat v_{\niters}}{\testMat z_{\niters}} \right|
\leq 
c \opnorm{\testMat}^2 \left(\frac{\sigstar^2}{\stpsz \strongconvex \niters^2}
+
\frac{ \normb{\naF{\theta_0}}}{\stpsz^4 \strongconvex^4 \niters^4} \right) \log \niters\\
+
c \opnorm{\testMat}^2 \holderconstprime
\left( 
\frac{\sigstartild^{3} }{\stpsz^{1/2}\strongconvex^{5/2}\niters^{2}} 
+
\frac{ \norm{\naF{\theta_0}}^{3} }{\stpsz^{7/2} \strongconvex^{11/2}\niters^{7/2}} 
\right),
\end{multline}
for some universal constant $c > 0$.
\end{lemma}
\noindent See \S\ref{subsec:proof-test-mat-lemmas} for the proof of both lemmas.

Taking these two lemmas as given, we now proceed with the proof of Corollary~\ref{theo_other_bounds}.

\pb\subsubsection{Proof of the MSE bound~\eqref{eq:main-mse-bound}}
We start with the following decomposition of the gradient:
\begin{align*}
    \nabla F (\theta_\niters) = \int_0^1 \nabla^2 F \big(\rho \thetastar + (1 - \rho) \theta_\niters \big) (\theta_\niters - \thetastar) d \rho,
\end{align*}
which leads to the following bound under Assumption~\ref{assu_holderhessian}:
\begin{align}
    \vecnorm{(\hessianstar)^{-1} \nabla F (\theta_\niters) - (\theta_\niters - \thetastar)}{2} &\leq \int_0^1 \vecnorm{(\hessianstar)^{-1} \left( \nabla^2 F \big(\rho \thetastar + (1 - \rho) \theta_\niters \big)  - \hessianstar\right) (\theta_\niters - \thetastar)}{2} d \rho\nonumber\\
    &\leq \frac{\holderconstprime}{\lammin (\hessianstar)} \vecnorm{\theta_\niters - \thetastar}{2}^{2} \leq \frac{\holderconstprime}{\lammin (\hessianstar) \strongconvex^{2}} \vecnorm{\nabla F (\theta_\niters)}{2}^{2}.\label{eq:linearization-error}
\end{align}

We can then upper bound the mean-squared error using the processes $(z_t)_{t \geq \burnin}$ and $(v_t)_{t \geq \burnin}$:

\begin{align}
\Exs\vecnorm{\theta_\niters - \thetastar}{2}^2 &\leq  \Exs \left( \vecnorm{(\hessianstar)^{-1} \nabla F (\theta_\niters)}{2} + \frac{\holderconstprime}{\strongconvex^{2} \lammin (\hessianstar)} \vecnorm{\nabla F (\theta_\niters)}{2}^{2} \right)^2 
\nonumber\\&\leq 
\Exs\vecnorm{(\hessianstar)^{-1} \big(v_{\niters + 1} - z_{\niters + 1} \big)}{2}^2 + \frac{2 \holderconstprime}{\lammin (\hessianstar)^2 \strongconvex^{2}} \Exs\vecnorm{\nabla F (\theta_\niters )}{2}^{3}
\nonumber\\&\qquad 
+ \frac{\holderconstprime^2}{\lammin (\hessianstar)^2 \strongconvex^{4}}  \Exs\vecnorm{\nabla F (\theta_\niters )}{2}^{4}
.\label{eq:theta-mse-decomp}
\end{align}
The first term in the bound~\eqref{eq:theta-mse-decomp} admits the following decomposition:
\begin{align*}
&
\Exs\vecnorm{(\hessianstar)^{-1} \big(z_{\niters + 1} - v_{\niters + 1} \big)}{2}^2 
\\&= 
\Exs\vecnorm{(\hessianstar)^{-1} z_{\niters + 1}}{2}^2 + \Exs\vecnorm{(\hessianstar)^{-1} v_{\niters + 1}}{2}^2 - 2  \Exs \inprod{(\hessianstar)^{-1} z_{\niters + 1}}{(\hessianstar)^{-1}  v_{\niters + 1}}
 .
\end{align*}
Note that the re-starting scheme in Algorithm~\ref{algo_multiepoch} gives the initial conditions:
\begin{align}
\Exs\vecnorm{\nabla F(\theta_0)}{2}^2 
\leq 
\frac{c \sigstar^2}{\burnin} 
,\qquad \mbox{and} \quad 
\left(\Exs\vecnorm{\nabla F(\theta_0)}{2}^4 \right)^{1/2} 
\leq 
\frac{c \sigstartild^2}{\burnin}
.\label{eq:init-condition-in-mse-proof}
\end{align}
Using these initial conditions, and invoking the Lemma~\ref{lemm_ztsharpbdd_testmat} with test matrix $\testMat = (\hessianstar)^{-1}$, we obtain the bound:
\begin{align*}
\Exs \vecnorm{(\hessianstar)^{-1} z_\niters}{2}^2
\leq 
\frac{1}{\niters} \mathrm{Tr} \left( (\hessianstar)^{-1} \SigStar (\hessianstar)^{- \top}\right) 
+
\frac{c}{\lammin (\hessianstar)^2} \Big\{
\frac{\sglip^2\stpsz}{\strongconvex}
+
\frac{\sglip}{\strongconvex \sqrt{\niters}} + \frac{\burnin}{ \niters}
\Big\} \frac{\sigstar^2 \log \niters}{\niters}
.
\end{align*}
Similarly, invoking Lemma~\ref{lemm_cross-term-bound-testmat} with test matrix $\testMat = (\hessianstar)^{-1}$, we have that:
\begin{align*}
\left| \Exs \binprod{(\hessianstar)^{-1} v_{\niters}}{(\hessianstar)^{-1} z_{\niters}} \right|
\leq  
\frac{c \sigstar^2\sqrt{\log \niters}}{\lammin (\hessianstar)^2 \stpsz \strongconvex \niters^2}
+
\frac{c \holderconstprime}{\lammin (\hessianstar)^2 \strongconvex^{2}} \cdot
\frac{\sigstartild^{3} }{(\stpsz \strongconvex)^{1/2} \niters^{2}}
.
\end{align*}
For the term $\Exs\vecnorm{(\hessianstar)^{-1} v_{\niters}}{2}^2$, Lemma~\ref{lemm_vtsharpbdd} along with the initial condition yields:
\begin{align*}
\Exs\vecnorm{(\hessianstar)^{-1} v_{\niters}}{2}^2 
\leq 
\frac{1}{\lammin (\hessianstar)^2} \Exs\vecnorm{ v_{\niters}}{2}^2 
\leq 
\frac{c \sigstar^2}{\lammin (\hessianstar)^2 \strongconvex \stpsz \niters^2}
.
\end{align*}

Collecting above bounds, we conclude that
\begin{multline}
\Exs\vecnorm{(\hessianstar)^{-1} \nabla F (\theta_\niters)}{2}^2 
\leq 
\frac{1}{\niters} \mathrm{Tr} \left( (\hessianstar)^{-1} \SigStar (\hessianstar)^{- \top}\right) +
\frac{c}{\lammin (\hessianstar)^2} \Big\{
\frac{\sglip^2\stpsz}{\strongconvex}
+
\frac{\sglip}{\strongconvex \sqrt{\niters}} + \frac{1}{\strongconvex \stpsz \niters}
\Big\} \frac{\sigstar^2 \log \niters}{\niters}
\\
+ \frac{c \holderconstprime}{\lammin (\hessianstar)^2 \strongconvex^{2}} \cdot
\frac{\sigstartild^{3} }{(\stpsz \strongconvex)^{1/2} \niters^{2}}
.\label{eq:zt-vt-decomp-bound-in-mse-proof}
\end{multline}
In order to bound the last two terms of the decomposition~\eqref{eq:theta-mse-decomp}, we recall from Lemma~\ref{lemm_finitebddHolder} and the initial condition~\eqref{eq:init-condition-in-mse-proof} that:
\begin{align*}
\left(\Exs\vecnorm{\nabla F (\theta_\niters)}{2}^4 \right)^{1/2} \leq \frac{c \sigstartild^2}{\niters}
.
\end{align*}
Combining with Eq.~\eqref{eq:zt-vt-decomp-bound-in-mse-proof} and substituting into the decomposition~\eqref{eq:theta-mse-decomp}, we conclude that:
\begin{multline*}
\Exs\vecnorm{\theta_\niters - \thetastar}{2}^2
\leq 
\frac{1}{\niters} \mathrm{Tr} \left( (\hessianstar)^{-1} \SigStar (\hessianstar)^{- \top}\right) +
\frac{c}{\lammin (\hessianstar)^2} \Big\{
\frac{\sglip^2\stpsz}{\strongconvex}
+
\frac{\sglip}{\strongconvex \sqrt{\niters}} + \frac{1}{\stpsz \strongconvex \niters}
\Big\} \frac{\sigstar^2 \log \niters}{\niters}
\\
+ \frac{c \holderconstprime}{\lammin (\hessianstar)^2 \strongconvex^{2}} \cdot
\left\{ \frac{\sigstartild^{3} }{(\stpsz \strongconvex)^{1/2} \niters^{2}} + \frac{\sigstartild^{3}}{\niters^{3/2}} + \frac{\holderconstprime}{\strongconvex^{2}}\frac{\sigstartild^{4}}{\niters^{2}} \right\}
.
\end{multline*}
Note in the last line, the second $O(\niters^{-3/2})$ term is always no smaller than the previous first term.
Taking $\niters = \numobs - \Epoch \shortepoch$ with $\numobs \geq 2  \Epoch \shortepoch$, some algebra then completes the proof of the desired bound.

\pb\subsubsection{Proof of the objective gap bound~\eqref{eq:main-val-bound}}
Applying second-order Taylor expansion with integral remainder, for any $\theta \in \real^d$, we note the following identity.
\begin{align*}
F (\theta) 
&= 
F (\thetastar) + \inprod{\theta - \thetastar}{\nabla F (\thetastar)} + (\theta - \thetastar)^\top \int_0^1 \nabla^2 F \big( \rho \theta + (1 - \rho) \thetastar \big) d \rho \cdot (\theta - \thetastar)
.
\end{align*}
Noting that $\nabla F (\thetastar) = 0$ and invoking assumption~\ref{assu_holderhessian}, we have that:
\begin{equation}\label{eq:function-value-decomp}
\begin{aligned}
\lefteqn{
F (\theta) 
}\\&\leq 
F (\thetastar) + \frac{1}{2} (\theta - \thetastar)^\top \hessianstar (\theta - \thetastar) +  \vecnorm{\theta - \thetastar}{2} \cdot  \int_0^1 \opnorm{\nabla^2 F \big( \rho \theta + (1 - \rho) \thetastar \big) - \hessianstar } d \rho \cdot \vecnorm{\theta - \thetastar}{2} 
\\&\leq 
F (\thetastar) + \frac{1}{2} (\theta - \thetastar)^\top \hessianstar (\theta - \thetastar) + \holderconstprime \vecnorm{\theta - \thetastar}{2}^{3}
.
\end{aligned}\end{equation}
Analogous to Eq.~\eqref{eq:linearization-error}, we have the bound:
\begin{align*}
\lefteqn{
\vecnorm{(\hessianstar)^{1/2} (\theta_\niters - \thetastar) - (\hessianstar)^{-1/2} \nabla F (\theta_\niters)}{2} 
}\\&\leq 
\int_0^1 \vecnorm{(\hessianstar)^{-1/2} \left( \nabla^2 F \big(\rho \thetastar + (1 - \rho) \theta_\niters \big)  - \hessianstar\right) (\theta_\niters - \thetastar)}{2} d \rho
\nonumber\\&\leq 
\frac{\holderconstprime}{\sqrt{\lammin (\hessianstar)}} \vecnorm{\theta_\niters - \thetastar}{2}^{2} \leq \frac{\holderconstprime}{\strongconvex^{2} \sqrt{\lammin (\hessianstar)} } \vecnorm{\nabla F (\theta_\niters)}{2}^{2}
.
\end{align*}
Denote the residual $q_t \mydefn (\hessianstar)^{1/2} (\theta_t - \thetastar) - (\hessianstar)^{-1/2} \nabla F (\theta_t)$. Substituting into the bound~\eqref{eq:function-value-decomp}, we have that:
\begin{align}
&
\Exs \left[ F (\theta_\niters) - F (\thetastar) \right]
\nonumber\\&\leq
\frac{1}{2} \Exs\vecnorm{(\hessianstar)^{-1/2} \nabla F (\theta) + q_\niters}{2}^2 + \holderconstprime \Exs\vecnorm{\theta_\niters - \thetastar}{2}^{3}
\nonumber\\&\leq
\frac{1}{2} \Exs\vecnorm{(\hessianstar)^{-1/2} \nabla F (\theta_\niters)}{2}^2
+
\frac{1}{\sqrt{\lammin (\hessianstar)}} \Exs \Big[ \vecnorm{q_t}{2} \cdot \vecnorm{\nabla F (\theta_\niters)}{2} \Big]
+
\frac{1}{2} \Exs\vecnorm{q_t}{2}^2
+
\frac{\holderconstprime}{\strongconvex^{3}} \Exs\vecnorm{\nabla F (\theta_\niters)}{2}^{3}
.\label{eq:func-val-gap-decomp}
\end{align}
For the first term, by applying Lemma~\ref{lemm_cross-term-bound-testmat} and~\ref{lemm_ztsharpbdd_testmat} with $\testMat = (\hessianstar)^{-1/2}$, we can obtain the following bound analogous to Eq.~\eqref{eq:zt-vt-decomp-bound-in-mse-proof}:
\begin{multline*}
    \Exs\vecnorm{(\hessianstar)^{-1/2} \nabla F (\theta_\niters)}{2}^2 \leq \frac{1}{2 \niters} \mathrm{Tr} \left( (\hessianstar)^{-1} \SigStar\right) +
\frac{c}{\lammin (\hessianstar)} \Big\{
\frac{\sglip^2\stpsz}{\strongconvex}
+
\frac{\sglip}{\strongconvex \sqrt{\niters}} + \frac{1}{\strongconvex \stpsz \niters}
\Big\} \frac{\sigstar^2 \log \niters}{\niters}\\
+ \frac{c \holderconstprime}{\lammin (\hessianstar)\strongconvex^{2}} \cdot
\frac{\sigstartild^{3} }{(\stpsz \strongconvex)^{1/2} \niters^{2}}.
\end{multline*}
For the rest of the terms, we recall that Lemma~\ref{lemm_finitebddHolder} with the initial condition~\eqref{eq:init-condition-in-mse-proof} gives the bound $\big(\Exs\vecnorm{\nabla F (\theta_\niters)}{2}^4 \big)^{1/2} \leq \frac{c \sigstartild^2}{\niters}$.
Substituting into the decomposition~\eqref{eq:function-value-decomp}, we obtain that:
\begin{align*}
\Exs \left[ F (\theta_\niters) - F (\thetastar) \right]
&\leq
\frac{1}{2 \niters} \mathrm{Tr} \left( (\hessianstar)^{-1} \SigStar\right) + \frac{c}{\lammin (\hessianstar)} \Big\{
\frac{\sglip^2\stpsz}{\strongconvex}
+
\frac{\sglip}{\strongconvex \sqrt{\niters}} + \frac{1}{\strongconvex \stpsz \niters}
\Big\} \frac{\sigstar^2 \log \niters}{\niters} \\
&\qquad + \frac{c \holderconstprime}{\strongconvex^{2}} \cdot
\frac{\sigstartild^{3} }{\strongconvex \niters^{3/2}} + \frac{\holderconstprime^2}{\strongconvex^{4}} \cdot
\frac{\sigstartild^{4} }{\lammin (\hessianstar) \niters^{2}}.
\end{align*}
Noting that $\niters = \numobs - \Epoch \shortepoch$ with $\numobs \geq 2  \Epoch \shortepoch$, we completes the proof of the desired bound.

\pb\subsection{Proof of Theorem~\ref{theo_asymptotic_iteration}}\label{sec_proof,theo_asymptotic}
Here we provide a two-step proof of Theorem~\ref{theo_asymptotic_iteration}.
We continue to adopt the $v_t$---$z_t$ decomposition as earlier used, and we proceed with the proof in two steps:

\paragraph{Step 1:}
We first claim the following single-epoch result, Eq.~\eqref{asymptoticlemm}, that under the setting of Theorem~\ref{theo_asymptotic_iteration} along with $\|\nabla F(\theta_0)\| = O(\sqrt{\stpsz\strongconvex\sigstar^2})$, the single-epoch estimator produced by Algorithm~\ref{algo_singleepoch} with burn-in time $
\burnin	=	\left\lceil \frac{24}{\stpsz\strongconvex} \right\rceil
$, as $\niters\to\infty$, $\stpsz\to 0$ such that $\stpsz \niters\to \infty$ satisfies the following convergence in probability:
\begin{align}\label{asymptoticlemm}
\sqrt{\niters} z_{\niters} - \frac{1}{\sqrt{\niters}}\sum_{s = 1}^T \noise_s(\thetastar)	\xrightarrow{p}	0
.
\end{align}
Taking this as given, we now combine Eq.~\eqref{asymptoticlemm} with our multi-epoch design Algorithm~\ref{algo_multiepoch} we can essentially assume without loss of generality that $\|\nabla F(\theta_0)\| = O(\sqrt{\stpsz\strongconvex\sigstar^2})$.
Under the current scaling condition, the final long epoch in Algorithm~\ref{algo_multiepoch} will be triggered with length $\niters = \nsamples - \inneriters\Epoch$, and hence we apply Eq.~\eqref{vtsharpbdd} so for some $C\le 56$ we have the initial condition holds: $
\Exs\|\nabla F(\theta_0^{(\stpsz)})\|_2^2
\le
\frac{C \sigstar^2}{\inneriters} 
=
O(\stpsz\strongconvex\sigstar^2)
$, so that as $\stpsz\niters\to \infty$,
$$
\niters\Exs \normb{v_{\niters}}
\le
O\left(
\frac{\sigstar^2}{\stpsz\strongconvex \niters} 
+
\frac{\stpsz\strongconvex\sigstar^2}{\stpsz^4 \strongconvex^4 \niters^3}
\right)
\to 0
.
$$
Therefore, $\sqrt{\niters} v_{\niters}\xrightarrow{p}	0$ holds.

Now to put together the pieces, note that $
\frac{1}{\niters} \sum_{s = 1}^\niters \noise_s(\thetastar)
$ is the average of $\mathrm{i.i.d.}$ random vectors of finite second moment. 
By standard CLT, we have
$$
\frac{1}{\sqrt{\niters}}\sum_{s = 1}^\niters \noise_s(\thetastar)
	\xrightarrow{d}
\mathcal{N} (0, \SigStar)
.
$$
Consequently, replacing $\niters$ by $\nsamples - \inneriters\Epoch$ we can apply Slutsky's rule of weak convergence and obtain the desired weak convergence: as $\stpsz\to 0$, $\nsamples\to \infty$ such that $\stpsz(\nsamples - \inneriters\Epoch)\to \infty$
\begin{align*}
\sqrt{\niters} \nabla F(\theta_{\niters - 1}^{(\stpsz)})
&=
\sqrt{\niters} v_{\niters} - \sqrt{\niters} z_{\niters}
\\&=
\sqrt{\niters} v_{\niters} - \left(\sqrt{\niters}z_{\niters} - \frac{1}{\sqrt{\niters}}\sum_{s = 1}^\niters \noise_s(\thetastar)\right) - \frac{1}{\sqrt{\niters}}\sum_{s = 1}^\niters \noise_s(\thetastar)
\xrightarrow{d}
\mathcal{N} (0, \SigStar)
.
\end{align*}
Due to our additional scaling condition, we can further replace $\niters = \nsamples - \inneriters\Epoch$ by $\nsamples$, concluding Theorem~\ref{theo_asymptotic_iteration}.

\paragraph{Step 2:}
We proceed to prove Eq.~\eqref{asymptoticlemm} with the extra initialization condition $\|\nabla F(\theta_0)\| = O(\sqrt{\stpsz\strongconvex\sigstar^2})$.
By Eqs.~\eqref{vtsharpbdd} and \eqref{ztsharpbdd}, we have for $\niters\ge \burnin$ there exist constants $\asympconst_1, \asympconst_2, \asympconst_3 > 0$ independent of $\stpsz,T$ but depends on the problem parameters $(\strongconvex, \smoothness, \sglip, \sigstar,\theta_0, \alpha)$, such that
\begin{align*}
\Exs\vecnorm{z_{\niters}}{2}^2 
&\le
\frac{2\asympconst_2}{\niters}
,
\end{align*}
and consequently, we have from Eq.~\eqref{vtsharpbdd} that
$$
\Exs \normb{v_{\niters}} 
\le
\frac{752\sigstar^2}{\stpsz\strongconvex \niters^2} 
+
\frac{69175}{\stpsz^4 \strongconvex^4 \niters^4} \normb{\nabla F(\theta_0)}
\le
\frac{\asympconst_1}{\niters} \left( \frac{1}{\stpsz \niters} + \frac{\stpsz}{\stpsz^4 \niters^3} \right)
\le
\frac{2\asympconst_1}{\stpsz\niters^2}
,
$$
and hence
\begin{align*}
\Exs\vecnorm{\nabla F (\theta_{\niters-1})}{2}^2 
\le
2\left(\Exs\vecnorm{v_{\niters}}{2}^2 + \Exs\vecnorm{z_{\niters}}{2}^2 \right)
\le
\frac{4\asympconst_1}{\stpsz\niters^2}
+
\frac{4\asympconst_2}{\niters}
\le
\frac{4\asympconst_3}{\niters}
.
\end{align*}
Note from the definition in Eq.~\eqref{eq:zt_defn}
$$
tz_t
=
\noise_t(\theta_{t - 1}) + (t-1) z_{t - 1} + (t-1) (\noise_t (\theta_{t - 1}) - \noise_t (\theta_{t - 2}))
.
$$
By setting $
A_t = (t - 1) (\noiset(\theta_{t - 1}) - \noiset(\theta_{t - 2})) + \noiset(\theta_{t - 1}) - \noiset(\thetastar)
$, the process $
\niters z_{\niters} - \sum_{s = 1}^\niters \noise_s(\thetastar)
=
\sum_{s = 1}^\niters A_s
$ is a martingale. To conclude the bound~\eqref{asymptoticlemm}, we only need to show the following relation as $\niters\to\infty$ and $\stpsz\to 0$:
\begin{align}\label{EqnSecondMoment-prime}
\Exs\left\|
\sqrt{\niters} z_{\niters} - \frac{1}{\sqrt{\niters}}\sum_{s = 1}^\niters \noise_s(\thetastar)
\right\|_2^2
=
\frac{1}{\niters}\sum_{s = 1}^\niters \Exs \|A_s\|^2
\to
0
.
\end{align}
Since we have
\begin{align*}
&\quad\,%%%%
\Exs\left\|
\sum_{s = \burnin+1}^\niters (s - 1) (\noise_s(\theta_{s-1}) - \noise_s(\theta_{s - 2}))
\right\|_2^2 
%%%%%%%%
=
\sum_{s = \burnin+1}^\niters (s - 1)^2 \Exs\vecnorm{\noise_s (\theta_{s-1}) - \noise_s (\theta_{s - 2})}{2}^2 
\\&\le
\sglip^2 \sum_{s = \burnin+1}^\niters (s - 1)^2 \Exs\vecnorm{\theta_{s-1} - \theta_{s - 2}}{2}^2
=
\stpsz^2 \sglip^2 \sum_{s = \burnin+1}^\niters (s - 1)^2 \Exs\vecnorm{v_{s-1}}{2}^2 
\\&\le
\stpsz^2 \sglip^2 \sum_{s = \burnin+1}^\niters (s - 1)^2 \frac{2\asympconst_1}{\stpsz^4 (s-1)^4}
\le
\frac{2\asympconst_1 \sglip^2}{\stpsz^2 \burnin}
.
\end{align*}
We note that
\begin{align*}
\Exs \left\|
\sum_{s = \burnin+1}^\niters \left(\noise_s(\theta_{s-1}) - \noise_s(\thetastar)\right)
\right\|_2^2 
&= 
\sum_{s = \burnin+1}^\niters \Exs\vecnorm{\noise_s (\theta_{s-1}) - \noise_s (\thetastar)}{2}^2
\\&\le
\sglip^2 \sum_{s = \burnin+1}^\niters \Exs\vecnorm{\theta_{s-1} - \thetastar}{2}^2 
\le
\frac{\sglip^2}{\strongconvex^2}\cdot 4\asympconst_3 \log\left(\frac{\niters}{\burnin}\right)
.
\end{align*}
Therefore, combining this with $
\Exs\|A_t\|_2^2
=%%%%
\Exs\left\|
(t - 1) (\noiset(\theta_{t - 1}) - \noiset(\theta_{t - 2})) + \noiset(\theta_{t - 1}) - \noiset(\thetastar)
\right\|_2^2
%%%%%%%%
\le
2\sglip^2\stpsz^2 (t - 1)^2 \Exs\vecnorm{v_{t - 1}}{2}^2
+
\frac{2\sglip^2}{\strongconvex^2} \Exs \left\|\nabla F(\theta_{t - 1}) \right\|_2^2
$ we have as $\niters\to\infty$, $\stpsz\to 0$:
\begin{align*}
\frac{1}{\niters}\sum_{t = \burnin+1}^\niters \Exs\|A_t\|_2^2
&\le%%%%
\frac{1}{\niters}\sum_{t = \burnin+1}^\niters \Exs\left\|
(t - 1) (\noiset(\theta_{t - 1}) - \noiset(\theta_{t - 2})) + \noiset(\theta_{t - 1}) - \noiset(\thetastar)
\right\|_2^2
%%%%%%%%
\\&\le
2\sglip^2\stpsz^2 \cdot \frac{1}{\niters}\sum_{t = \burnin+1}^\niters (t - 1)^2 \Exs\vecnorm{v_{t - 1}}{2}^2
+
\frac{2\sglip^2}{\strongconvex^2} \cdot \frac{1}{\niters}\sum_{t = \burnin+1}^\niters \Exs \left\|\nabla F(\theta_{t - 1}) \right\|_2^2
%%%%%%%%
\\&\le
2\sglip^2\stpsz^2 \cdot \frac{1}{\niters}\sum_{t = \burnin+1}^\niters (t - 1)^2 \frac{2\asympconst_1}{\stpsz (t-1)^2}
+
\frac{2\sglip^2}{\strongconvex^2} \cdot \frac{1}{\niters}\sum_{t = \burnin+1}^\niters \frac{4\asympconst_3}{t}
%%%%%%%%
\\&=
4\asympconst_1\sglip^2\stpsz
+
\frac{2\sglip^2}{\strongconvex^2}\cdot\frac{4\asympconst_3\log\left(\frac{\niters}{\burnin}\right)}{\niters}
,
\end{align*}
i.e.~the limit \eqref{EqnSecondMoment-prime} holds, which implies $\sqrt{\niters} z_{\niters} - \frac{1}{\sqrt{\niters}}\sum_{s = 1}^\niters \noise_s(\thetastar)	\xrightarrow{p}	0$, completing our proof of Eq.~\eqref{asymptoticlemm}.

\pb\section{Asymptotic results for single-epoch fixed-step-size \ROOTSGD}\label{sec_algoresult_asy_hessiansmooth_additional}

In this section, we complement Theorem \ref{theo_asymptotic_iteration} in \S\ref{sec_algoresult_asy_hessiansmooth} and establish an additional asymptotic normality result for \ROOTSGD with large step-size.
Notably, the covariance of such asymptotic distribution is the sum of the optimal Gaussian limit and a correction term depending on the step-size, which exactly corresponds to existing results on fine-grained CLT for linear stochastic approximation with fixed step-size~\cite{mou2020linear}.

First, in order to obtain asymptotic results for single-epoch constant-step-size \ROOTSGD, we impose the following slightly stronger assumptions on the smoothness of stochastic gradients and Hessians:

\begin{description}
\item[(CLT.A)] 
For any $\theta \in \real^d$ we have
\begin{subequations}
\begin{align}
\sup_{v \in \sphere^{d - 1}}
\Exs\vecnorm{(\nabla^2 f (\theta; \xi) - \nabla^2 f (\thetastar; \xi) ) v}{2}^2 
\leq 
\beta^2 \vecnorm{\theta - \thetastar}{2}^2
.
\end{align}
\item[(CLT.B)]
The fourth moments of the stochastic gradient
vectors at $\thetastar$ exist, and in particular we have
\begin{align}
\Exs\vecnorm{\nabla f(\thetastar; \xi)}{2}^4
	<	\infty
,\qquad \mbox{and} \quad \sglip' \mydefn
\sup_{v \in \sphere^{d - 1}}\big(\Exs\vecnorm{\nabla^2 f (\thetastar; \xi) v}{2}^4 \big)^{1/4}
	<	\infty
.
\end{align}
\end{subequations}
\end{description}
Note that both conditions are imposed solely at the optimal point $\thetastar$; we do not impose globally uniform bounds in $\real^d$.

Defining the random matrix $\NoiseAplain (\theta) \mydefn \nabla^2 f(\theta; \xi) - \nabla^2 F (\theta)$ for any $\theta \in \real^d$, we consider the following matrix equation (a.k.a.~\emph{modified Lyapunov equation}):
\begin{align}
\label{eq:lambda-eq-constant-step-size}
\Lambda \hessianstar + \hessianstar \Lambda - \stpsz \Exs \big[ \NoiseAplain(\thetastar) \Lambda \NoiseAplain(\thetastar) \big] - \stpsz \hessianstar \Lambda \hessianstar 
= 
\stpsz \SigStar
.
\end{align}
in the symmetric matrix $\Lambda$. It can be shown that under the
given assumptions, this equation has a unique solution---denoted
$\Lambda_\stpsz$---which plays a key role in the following theorem.

\begin{theorem2}[Asymptotic efficiency, single-epoch \ROOTSGD]\label{theo_asymptotic_iteration_const_step_size}
Suppose that Assumptions~\ref{assu_StrcvxSmooth},~\ref{assu_noisethetastar},
and~\ref{assu_smoothnoise} are satisfied, as well as {\rm (CLT.A)} and {\rm (CLT.B)}. Then there exist constants $c_1, c_2$, given the step-size $\stpsz \in \Big(0, c_1 (\frac{\strongconvex}{\sglip^2} \wedge
\frac{1}{\smoothness} \wedge
\frac{\strongconvex^{1/3}}{\sglip'^{4/3}})\Big)$, and burn-in time
$\burnin = \frac{c_2}{\strongconvex \stpsz}$, we have
\begin{align}\label{asymptotic_iteration_const_step_size}
\sqrt{T} (\theta_T - \thetastar)
	\xrightarrow{d}
\mathcal{N}\left(0, 
(\hessianstar)^{-1} \left(\SigStar + \Exs\left[\NoiseAplain (\thetastar) \Lambda_\stpsz \NoiseAplain (\thetastar)\right]\right) (\hessianstar)^{-1} 
\right)
.
\end{align}
\end{theorem2}
\noindent See \S\ref{SecProofthm:asymptotic-const-step-size} for the proof of this theorem.

A few remarks are in order. 
First, we observe that the asymptotic covariance in Eq.~\eqref{asymptotic_iteration_const_step_size} is the sum of the matrix $(\hessianstar)^{-1} \SigStar (\hessianstar)^{-1}$ and an additional correction term defined in Eq.~\eqref{eq:lambda-eq-constant-step-size}. 
The asymptotic covariance of $(\hessianstar)^{-1} \SigStar (\hessianstar)^{-1}$ matches the standard Cram\'{e}r-Rao lower bound in the asymptotic statistics literature \cite{VAN-WELLNER,VAN[derVaart]} and matches the optimal rates achieved in the theory of stochastic approximation~\cite{KUSHNER-YIN,polyak1992acceleration,ruppert1988efficient}.
The correction term is of the same form as that of the constant-step-size linear stochastic approximation of the Polyak-Ruppert-Juditsky averaging prcedure as derived in~\cite{mou2020linear}, while our Proposition \ref{theo_asymptotic_iteration_const_step_size} is applicable to more general nonlinear stochastic problems. 
%\iffalse
For instance in our setting, the correction terms tends to zero as the (constant) step-size decreases to zero, which along with a trace bound leads to the following asymptotics as $\niters\to\infty$ (see~\cite{mou2020linear}):
\begin{align*}
T \Exs\vecnorm{\nabla F(\theta_T)}{2}^2
\sim
\mathrm{Tr} \left( 
\SigStar
+ 
\Exs\left[\NoiseAplain (\thetastar) \Lambda_\stpsz \NoiseAplain (\thetastar)\right]
\right) 
\le 
\left(1 + \frac{ \sglip^2 \stpsz}{\strongconvex} \right) \sigmastarsq
.
\end{align*}
The message conveyed by the last display is consistent with the leading two terms in our earlier nonasymptotic bound Eq.~\eqref{rootHolder_single} in Proposition \ref{theo_rootHolder_single}, and thanks to our additional smoothness assumptions {\rm (CLT.A)} and {\rm (CLT.B)} we are able to characterize this correction term in a more fine-grained fashion as in the asymptotic covariance of Eq.~\eqref{asymptotic_iteration_const_step_size}.
%\fi
Second, we note that Proposition \ref{theo_asymptotic_iteration_const_step_size} has an additional requirement on the step-size, needing it to be upper bounded by $\frac{\strongconvex^{1/3}}{\sglip’^{4/3}}$. 
This is a mild requirement on the step-size. 
In particular, for applications where the noises are light-tailed, $\sglip'$ and $\sglip$ are of the same order, and the additional requirement $\stpsz < \frac{c \strongconvex^{1/3}}{\sglip’^{4/3}}$ is usually weaker than the condition $\stpsz < \frac{c \strongconvex}{\sglip^2}$ needed in the previous section.

\pb\subsection{Proof of Proposition \ref{theo_asymptotic_iteration_const_step_size}}\label{SecProofthm:asymptotic-const-step-size}

Denote $H_t (\theta) \mydefn \nabla^2 f (\theta; \xi_t)$ and $\NoiseAt
(\theta) \mydefn H_t (\theta) - \nabla^2 F (\theta)$. Intuitively,
since the sequence $\theta_t$ is converging to $\thetastar$ at a $1 /
\sqrt{t}$ rate, replacing $\theta_{s - 1}$ with $\thetastar$ will only
lead to a small change in the sum. For the martingale $\Psi_t$, each
term can be written as:
\begin{align*}
    t (\noise_t (\theta_{t - 1}) - \noise_t (\theta_{t - 2})) = t
    \int_0^1 \NoiseAt \left(\rho \theta_{t - 2} + (1 - \rho)
    \theta_{t - 1} \right) (\theta_{t - 1} - \theta_{t - 2}) d \rho.
\end{align*}
By Assumption (CLT.A), this quantity should approach $\stpsz
\NoiseAt (\thetastar) \cdot (t v_{t - 1})$. If we can show the
convergence of the sequence $\{t v_t \}_{t \geq \Tburnin}$ to a
stationary distribution, then the asymptotic result follows from the
Birkhoff ergodic theorem and a martingale CLT. While the process $\{t v_t
\}_{t \geq \burnin}$ is not Markovian, we show that it can be
well-approximated by a time-homogeneous Markov process that we
construct in the proof.

In particular, consider the auxiliary process $\{y_t\}_{t \geq
  \burnin}$, initialized as $y_{\burnin} = \burnin v_{\burnin}$ and
updated as
\begin{align}
\label{eq:yt-defn}  
y_{t} = y_{t - 1} - \stpsz H_t (\thetastar) y_{t - 1} + \noise_t
(\thetastar), \quad \mbox{for all $t \geq \burnin + 1$.}
\end{align}
Note that $\{ y_t \}_{t \geq \burnin}$ is a time-homogeneous Markov
process that is coupled to $ \{ (\theta_t, v_t, z_t) \}_{t \geq
  \burnin}$. We have the following coupling estimate:
\begin{lemma}
\label{lemma-vt-coupling-estimate}
Supposing that Assumptions~\ref{assu_StrcvxSmooth},
~\ref{assu_noisethetastar} and~\ref{assu_smoothnoise}, as well as Conditions
(CLT.A) and (CLT.B) hold, then for any iteration $t \geq \burnin$ and
any step-size $\stpsz \in (0, \frac{1}{4 \smoothness} \wedge
\frac{\strongconvex}{8 \sglip^2})$, we have
\begin{align*}
  \Exs\vecnorm{t v_t - y_t}{2}^2 \leq \frac{\kcon}{\sqrt{t}},
\end{align*}
for a constant $\kcon$ depending on the smoothness and strong convexity parameters $\smoothness, \sglipfour, \strongconvex, \beta$ and the step-size $\stpsz$, but independent of $t$.
\end{lemma}
\noindent See \S\ref{subsubsec:proof-coupling-estimate} for the
proof of this lemma. \\

We also need the following lemma, which provides a convenient bound on
the difference $H_t(\theta) - H_t(\thetastar)$ for a vector
$\theta$ chosen in the data-dependent way.
\begin{lemma}
\label{lemma:hessian-noise-mixed-bound}
Suppose that Assumptions~\ref{assu_StrcvxSmooth},
~\ref{assu_noisethetastar} and~\ref{assu_smoothnoise}, as well as Conditions
(CLT.A) and (CLT.B) hold. Then for any iteration $t \geq \burnin$, any
step-size $\stpsz \in (0, \frac{1}{4 \smoothness} \wedge
\frac{\strongconvex}{8\sglip^2})$ and for any random vector
$\tilde{\theta}_{t - 1} \in \filtration_{t - 1}$, we have
\begin{align*}
\Exs\vecnorm{\left[ H_t (\tilde{\theta}_{t - 1}) - H_t (\thetastar) \right] y_{t - 1}}{2}^2 
\leq 
\Contwo \sqrt{\Exs\vecnorm{\tilde{\theta}_{t - 1} - \thetastar}{2}^2}
,
\end{align*}
where $\Contwo$ is a constant independent of $t$ and the choice of $\tilde{\theta}_{t - 1}$.
\end{lemma}
\noindent See \S\ref{subsubsec:proof-hessian-noise-mixed-bound}
for the proof of this lemma. \\

\noindent Finally, the following lemma characterizes the behavior of
the process $\{ y_t\}_{t \geq \burnin}$ defined in
Eq.~\eqref{eq:yt-defn}:
\begin{lemma}
 \label{lemma:yt-properties}
Suppose that Assumptions~\ref{assu_StrcvxSmooth},
~\ref{assu_noisethetastar} and~\ref{assu_smoothnoise}, as well as Conditions
(CLT.A) and (CLT.B) hold. Then for any iteration $t \geq \burnin$ and
any step-size $\stpsz \in (0, \frac{1}{4 \smoothness} \wedge
\frac{\strongconvex}{16\sglip^2} \wedge \frac{\strongconvex^{1/3}}{6 \sglipfourbracket^{4/3}})$, we have
\begin{align*} 
\Exs (y_t) = 0 \quad \mbox{for all $t \geq \burnin$},
\qquad\mbox{and}\quad
\sup_{t \geq \burnin} \Exs\vecnorm{y_t}{2}^4 < \aconprime
,
\end{align*}
for a constant $\aconprime > 0$, which is independent of $t$.
Furthermore, the process
$\{y_t\}_{t \geq 0}$ has a stationary distribution with finite second
moment, and a stationary covariance $Q_\stpsz$ that satisfies the
equation
\begin{align*}
\hessianstar Q_\stpsz + Q_\stpsz \hessianstar
-
\stpsz\left[
  \hessianstar Q_\stpsz \hessianstar + \Exs (\NoiseAplain(\thetastar) Q_\stpsz \NoiseAplain(\thetastar))
\right]
=
\frac{1}{\stpsz} \SigStar
.
\end{align*}
\end{lemma}
\noindent See \S\ref{subsubsec:yt-properties} for the proof of this lemma.
\\

Taking these three lemmas as given, we now proceed with the proof of Proposition \ref{theo_asymptotic_iteration_const_step_size}. 
We first define two auxiliary processes:
\begin{align*}
N_T \mydefn \sum_{t = \burnin+1}^T \noise_t (\thetastar),
\qquad
\Upsilon_T \mydefn \stpsz \sum_{t = \burnin+1}^T \NoiseAt
(\thetastar) y_{t - 1}
.
\end{align*}
Observe that both $N_T$ and $\Upsilon_T$ are martingales adapted to $(\filtration_t)_{t \geq \burnin}$. 
In the following, we first bound the differences $\vecnorm{M_T - N_T}{2}$ and $\vecnorm{\Psi_T -  \Upsilon_T}{2}$, respectively, and then show the limiting distribution results for $N_T + \Upsilon_T$.

\newcommand{\constA}{a_0}

By Theorem~\ref{theo_finitebdd_single}, define $\constA \mydefn \frac{28 \sigmastarsq}{\strongconvex^2} + \frac{2700}{\stpsz^2 \strongconvex^4 \burnin} \vecnorm{\nabla F(\theta_0)}{2}^2$, we have
\begin{align}\label{eq:theta-bound-in-proof-of-clt}  
  \Exs\vecnorm{\theta_t - \thetastar}{2}^2 \leq \frac{1}{\strongconvex^2} \Exs\vecnorm{ \nabla F(\theta_t)}{2}^2 & \leq \frac{2700 \; \vecnorm{\nabla F(\theta_0)}{2}^2}{\stpsz^2 \strongconvex^4 (t + 1)^2} +
\frac{28 \; \sigmastarsq}{\strongconvex^2(t + 1)} \leq \frac{\constA}{t + 1},\quad
  \mbox{for all $t \geq \burnin$.}
\end{align}
Applying the bound~\eqref{eq:theta-bound-in-proof-of-clt} with
Assumption~\ref{assu_smoothnoise}, we have
\begin{align}\label{eq:mn-bound-in-proof-of-clt}
\Exs\vecnorm{M_T - N_T}{2}^2 
= 
\sum_{t = \burnin+1}^T \Exs\vecnorm{\noise_t (\theta_{t - 1}) - \noise_t (\thetastar)}{2}^2
\leq 
\sglip^2 \sum_{t = \burnin+1}^T \Exs\vecnorm{\theta_{t - 1} - \thetastar}{2}^2 \leq \constA \sglip^2 \log\niters
. 
\end{align}
For the process $\Upsilon_T$, by the Cauchy-Schwartz inequality, we have
\begin{align*}
&
\Exs\vecnorm{\Psi_T - \Upsilon_T}{2}^2 
= 
\sum_{t = \burnin+1}^T \Exs\vecnorm{\stpsz \NoiseAt (\thetastar) y_{t - 1} - (t - 1)(\noise_t (\theta_{t - 1}) - \noise_t (\theta_{t - 2}) )}{2}^2
\notag \\& \leq 
\stpsz^2 \sum_{t = \burnin+1}^T \Exs \int_0^1 \vecnorm{\NoiseAt (\thetastar) y_{t - 1} - \NoiseAt (\rho \theta_{t - 1} + (1 - \rho) \theta_{t - 2}; \xi_t) (t - 1) v_{t - 1}}{2}^2 d \rho 
%\notag \\& 
\le
I_I + I_2
,
\end{align*}
where we define
\begin{align*}
I_1 & \mydefn 2 \stpsz^2 \sum_{t = \burnin+1}^T \Exs \int_0^1 \vecnorm{ \left(\NoiseAt (\thetastar) - \NoiseAt (\rho \theta_{t - 1} + (1 - \rho) \theta_{t - 2}) \right) y_{t - 1}}{2}^2 d \rho
, \quad \mbox{and} \\
I_2 & \mydefn 2 \stpsz^2 \sum_{t = \burnin+1}^T \Exs \int_0^1 \vecnorm{ \NoiseAt (\rho \theta_{t - 1} + (1 - \rho) \theta_{t - 2}) (y_{t - 1} - (t - 1) v_{t - 1})}{2}^2 d \rho
.
\end{align*}
We bound each of these two terms in succession.

\paragraph{\textbf{Bound on $I_1$:}}
In order to bound the term $I_1$, we apply Lemma~\ref{lemma:hessian-noise-mixed-bound} with the choice
\begin{align*}
\tilde{\theta}_{t - 1} = \rho \theta_{t - 1} + (1 - \rho)
\theta_{t - 2} \in \filtration_{t - 1},
\end{align*}
so as to obtain
\begin{align*}
\Exs\vecnorm{\left( H_t (\tilde{\theta}_{t - 1}) - H_t(\thetastar) \right) y_{t - 1}}{2}^2 
\leq 
\Contwo \sqrt{\Exs\vecnorm{\tilde{\theta}_t - \thetastar}{2}^2}
.
\end{align*}
Applying the Cauchy-Schwartz inequality yields
\begin{align*}
\Exs\vecnorm{ \left( \nabla^2 F (\tilde{\theta}_{t - 1}) - \nabla^2 F (\thetastar) \right) y_{t - 1}}{2}^2 
\leq 
\Exs\vecnorm{\left( H_t (\tilde{\theta}_{t - 1}) - H_t (\thetastar) \right) y_{t - 1}}{2}^2 
\leq 
\Contwo \sqrt{\Exs\vecnorm{\tilde{\theta}_t - \thetastar}{2}^2}
.
\end{align*}
Putting the two bounds together, we obtain:
\begin{align*}
&
\Exs\vecnorm{\left( \NoiseAt(\tilde{\theta}_{t - 1}) - \NoiseAt(\thetastar) \right) y_{t - 1}}{2}^2 
\\&\leq 
2 \Exs\vecnorm{\left( H_t (\tilde{\theta}_{t - 1}) - H_t (\thetastar) \right) y_{t - 1}}{2}^2 + 2 \Exs\vecnorm{ \left( \nabla^2 F (\tilde{\theta}_{t - 1}) - \nabla^2 F (\thetastar) \right) y_{t - 1}}{2}^2 
\notag \\&\leq 
4\Contwo \sqrt{\Exs\vecnorm{\tilde{\theta}_t - \thetastar}{2}^2}
.
\end{align*}
Thus, we find that
\begin{align*}
&\quad\,%%%%
\Exs\vecnorm{ \left(\NoiseAt (\thetastar) - \NoiseAt (\rho \theta_{t - 1} + (1 - \rho) \theta_{t - 2}) \right) y_{t - 1}}{2}^2 
\leq 
4 \Contwo \sqrt{\Exs\vecnorm{\rho \theta_{t - 1} + (1 - \rho) \theta_{t - 2} - \thetastar}{2}^2} 
\notag \\& \leq 
4 \Contwo \left( \sqrt{\Exs\vecnorm{\theta_{t - 1} - \thetastar}{2}^2} + \sqrt{\Exs\vecnorm{\theta_{t - 2} - \thetastar}{2}^2}\right) 
%\notag \\& 
\leq 
4 \Contwo \sqrt{\constA} \left(\frac{1}{\sqrt{t - 1}} + \frac{1}{\sqrt{t - 2}} \right) 
%\notag \\& 
\leq 
\frac{16 \Contwo \sqrt{\constA}}{\sqrt{t}}
,
\end{align*}
where in the last step we used the inequality~\eqref{eq:theta-bound-in-proof-of-clt}. 
Summing over $t$ from $\burnin+1$ to $\niters$ yields the bound
\begin{align*}
I_1
\leq
2 \stpsz^2 \sum_{t = \burnin+1}^T \frac{16 \Contwo \sqrt{\constA}}{\sqrt{t}}
\leq
64 \stpsz^2 \Contwo \sqrt{\constA T}
.
\end{align*}

%%%%%%%%%%%%%%%%%%%%%%%%%%%%%%%%%%%%%%%%%%%%%%%%%%%%%%%%%%%%%%%%%%%%%%%%%%%%%%%%%%%%%%%%%

\paragraph{\textbf{Bound on $I_2$:}}
Turning to the term $I_2$, by Assumption~\ref{assu_smoothnoise} and
Lemma~\ref{lemma-vt-coupling-estimate}, we note that:
\begin{align*}
    I_2 \leq 2 \stpsz^2 \sum_{t = \burnin+1}^T \sglip^2 \Exs
    \vecnorm{y_{t - 1} - (t - 1) v_{t - 1}}{2}^2 \leq 2 \stpsz^2
    \sglip^2 \sum_{t = \burnin+1}^T \frac{\kcon}{\sqrt{t}} \leq 4
    \stpsz^2 \sglip^2 \kcon \sqrt{T}.
\end{align*}
Putting these inequalities together, we conclude that:
\begin{align}
\label{EqnPsiUpBound}  
     \Exs\vecnorm{\Psi_T - \Upsilon_T}{2}^2 \leq (64 \stpsz^2 \Contwo
     \sqrt{\constA} + 4 \stpsz^2 \sglip^2 \kcon)\sqrt{T}.
\end{align}

Now we have the estimates for the quantities $\vecnorm{\Psi_T - \Upsilon_T}{2}$ and $\vecnorm{M_T - N_T}{2}$. In the following, we first prove the CLT for $N_T + \Upsilon_T$, and then use the error bounds to establish CLT for $M_T + \Psi_T$, which ultimately implies the desired limiting result for $\sqrt{T} (\theta_T - \thetastar)$ 

Define $\nu_t \mydefn \noise_t (\thetastar) + \stpsz \NoiseAt
(\thetastar) y_{t - 1}$. By definition, $N_T + \Upsilon_T = \sum_{t =
  \burnin}^T \nu_t$, and we have
\begin{multline*}
  \Exs (\nu_t \nu_t^\top) = \Exs (\noise_t (\thetastar) \noise_t
  (\thetastar)^\top ) + \Exs \left( \NoiseAt (\thetastar) y_{t - 1}
  y_{t - 1}^\top \NoiseAt (\thetastar)^\top \right) \\
  + \Exs \left(\noise_t (\thetastar) y_{t - 1}^\top \NoiseAt
  (\thetastar)^\top \right) + \Exs \left( (\NoiseAt (\thetastar) y_{t
    - 1} \noise_t (\thetastar)^\top \right).
\end{multline*}
For the first term, we have $ \Exs ( \noise_t (\thetastar) \noise_t
(\thetastar)^\top ) = \SigStar$ by definition.

For the second term, according to Lemma~\ref{lemma:yt-properties}, we note that the time-homogeneous Markov process $\{ y_t\}_{t\geq \burnin}$ converges asymptotically to a stationary distribution with covariance $Q_\stpsz$. 
Invoking the Birkhoff ergodic theorem, we have
\begin{align*}
\frac{1}{T} \sum_{t = \burnin+1}^T \Exs \left( \NoiseAt (\thetastar)
y_{t - 1} y_{t - 1}^\top \NoiseAt (\thetastar)^\top \mid
\filtration_{t - 1} \right) 
&= 
\Exs (\NoiseAplain (\thetastar) \otimes \NoiseAplain (\thetastar)) 
\left[ \frac{1}{T} \sum_{t = \burnin+1}^T y_{t - 1} y_{t - 1}^\top \right] 
\\&\xrightarrow{p} 
\Exs \left( \NoiseAplain (\thetastar) Q_\stpsz \NoiseAplain (\thetastar)^\top \right)
.
\end{align*}
For the cross terms, we note that:
\begin{align*}
     \Exs \left(\noise_t (\thetastar) y_{t - 1}^\top \NoiseAt
     (\thetastar)^\top \,|\, \filtration_{t - 1} \right) = \Exs (\noise
     (\thetastar) \otimes \NoiseAplain (\thetastar) ) [y_{t - 1}].
\end{align*}
Note that by Lemma~\ref{lemma:yt-properties}, we have $\Exs (y_t) = 0$
for any $t \geq \burnin$. By the weak law of large numbers, we have
$\frac{1}{T} \sum_{t = \burnin+1}^T y_t \xrightarrow{p} 0$. Putting
together these inequalities, we find that
\begin{multline*}
\frac{1}{T} \sum_{t = \burnin+1}^T \Exs \left( \nu_t \nu_t^\top \mid
\filtration_{t - 1} \right) = \frac{1}{T} \sum_{t = \burnin+1}^T \Big(
\SigStar + \stpsz^2 \Exs (\NoiseAplain (\thetastar) \otimes
\NoiseAplain (\thetastar)) [y_t y_t^\top] \\
+ \stpsz \Exs (\noise (\thetastar) \otimes \NoiseAplain
(\thetastar) ) [y_{t - 1}] + \stpsz \Exs ( \NoiseAplain (\thetastar)  \otimes \noise (\thetastar)) [y_{t - 1}]
\Big),
\end{multline*}
and hence the random matrix $\frac{1}{T} \sum_{t = \burnin+1}^T \Exs
\left( \nu_t \nu_t^\top \mid \filtration_{t - 1} \right)$ converges in
probability to the matrix
\begin{align*}
  \SigStar + \Exs \left( \NoiseAplain (\thetastar) \Lambda_\stpsz
  \NoiseAplain (\thetastar)^\top \right).
\end{align*}

To prove the limiting distribution result, we use standard martingale CLT (c.f.~Corollary 3.1 in~\cite{HALL-HEYDE}). It remains to verify the conditional Lindeberg condition. Indeed, for any $\varepsilon > 0$, a straightforward calculation yields:
\begin{multline*}
R_T (\varepsilon) 
\mydefn 
\sum_{t = \burnin+1}^T \Exs \left( \vecnorm{\frac{\nu_t}{\sqrt{T}}}{2}^2 \bm{1}_{\vecnorm{\frac{\nu_t}{\sqrt{T}}}{2} > \varepsilon} \,\bigg|\, \filtration_{t - 1} \right) 
\\\overset{(i)}{\leq} 
\frac{1}{T} \sum_{t = \burnin+1}^T \sqrt{\Exs \left( \vecnorm{\nu_t}{2}^4 \,|\, \filtration_{t - 1} \right)} \cdot \sqrt{ \Prob \left(\vecnorm{\nu_t}{2} > \varepsilon \sqrt{T} \,|\, \filtration_{t - 1}  \right)}
\overset{(ii)}{\leq} 
\frac{1}{T} \sum_{t = \burnin+1}^T \frac{1}{(\varepsilon \sqrt{T})^2 } \Exs \left( \vecnorm{\nu_t}{2}^4 \,|\, \filtration_{t - 1} \right)
.
\end{multline*}
In step $(i)$, we use the Cauchy-Schwartz inequality, and in step $(ii)$, we use the Markov inequality to bound the conditional probability.

Using the condition (CLT.B) and Young's inequality, we note that:
\begin{align*}
    \Exs \left( \vecnorm{\nu_t}{2}^4 \,|\, \filtration_{t - 1} \right) \leq 8 \Exs\vecnorm{\noise (\thetastar)}{2}^4 + 8 \sglipfourbracket^4 \vecnorm{y_{t - 1}}{2}^4.
\end{align*}
Plugging back to the upper bound for $R_T (\varepsilon)$, and applying Lemma~\ref{lemma:yt-properties}, as $\niters \rightarrow \infty$, we have
\begin{align*}
    \Exs [R_T (\varepsilon)] \leq \frac{8}{T \varepsilon^2} \Exs\vecnorm{\noise (\thetastar)}{2}^4 + \frac{8 \sglipfourbracket^4}{T^2 \varepsilon^2} \sum_{t = \burnin+1}^T \Exs\vecnorm{y_{t - 1}}{2}^4 \leq  \frac{8}{T \varepsilon^2} \Exs\vecnorm{\noise (\thetastar)}{2}^4 + \frac{8 \sglipfourbracket^4}{T \varepsilon^2} \aconprime \rightarrow 0.
\end{align*}

Note that $R_T (\varepsilon) \geq 0$ by definition. The limit statement implies that $R_T (\varepsilon) \xrightarrow{p} 0$, for any $\varepsilon > 0$. Therefore, the conditional Lindeberg condition holds true, and we have the CLT:
\begin{align*}
\frac{N_T + \Upsilon_T}{\sqrt{T}} 
\xrightarrow{d}
\mathcal{N} \left(0, \SigStar + \Exs\left[\NoiseAplain (\thetastar) \Lambda_\stpsz \NoiseAplain (\thetastar)\right] \right)
.
\end{align*}
By the second-moment estimates~\eqref{eq:mn-bound-in-proof-of-clt} and~\eqref{EqnPsiUpBound}, we have
\begin{align*}
\frac{\vecnorm{\Upsilon_T - \Psi_T}{2}}{\sqrt{T}}	\xrightarrow{p} 0
, \qquad
\frac{\vecnorm{M_T - N_T}{2}}{\sqrt{T}}		\xrightarrow{p} 0
.
\end{align*}
With the burn-in time $\burnin$ fixed, we also have $\frac{\burnin}{T} z_{\burnin} \xrightarrow{p} 0$. 
By Slutsky's theorem, we have
\begin{align*}
\sqrt{T} z_T	\xrightarrow{d}	\mathcal{N} \left(0, \SigStar + \Exs\left( \NoiseAplain (\thetastar) \Lambda_\stpsz \NoiseAplain (\thetastar)^\top \right) \right)
 .
\end{align*}
Note that $\nabla F (\theta_{t - 1}) = v_t - z_t$. 
By Lemma~\ref{lemma-vt-coupling-estimate} and Lemma~\ref{lemma:yt-properties}, we have
\begin{align*}
\Exs\vecnorm{v_t}{2}^2 
\leq 
\frac{2}{t^2} \Exs\vecnorm{t v_t -  y_t}{2}^2 + \frac{2}{t^2} \Exs\vecnorm{y_t}{2}^2 \leq \frac{2}{t^2}\left( \sqrt{\aconprime} + \frac{\kcon}{\sqrt{t}} \right)
,
\end{align*}
which implies that $\sqrt{t} v_t \xrightarrow{p} 0$. Recall that $z_t = v_t - \nabla F (\theta_{t - 1})$. By
Slutsky's theorem, we obtain:
\begin{align*}
\sqrt{T} \cdot \nabla F (\theta_T) \xrightarrow{d} \mathcal{N} \left(
0, \SigStar + \Exs\left[ \NoiseAplain (\thetastar) \Lambda_\stpsz \NoiseAplain (\thetastar) \right] \right)
 .
\end{align*}

Finally, we note that for $\theta \in \real^d$, we have
\begin{align*}
{\vecnorm{\nabla F (\theta) - \hessianstar (\theta - \thetastar)}{2}}
&= {\vecnorm{\int_0^1 \nabla^2 F (\thetastar + \rho (\theta -
    \thetastar)) (\theta - \thetastar) d \rho - \hessianstar (\theta
    - \thetastar)}{2}} 
    \\ &\leq 
    \int_0^1 \opnorm{\nabla^2 F
  (\thetastar + \rho (\theta - \thetastar)) - \hessianstar} \cdot
\vecnorm{\theta - \thetastar}{2} d \rho
\\ &\leq 
\vecnorm{\theta -
  \thetastar}{2} \cdot \sup_{\vecnorm{\theta' - \thetastar}{2} \leq
  \vecnorm{\theta - \thetastar}{2}} \opnorm{\nabla^2 F(\theta') -
  \hessianstar}.
\end{align*}
By Assumption (CLT.A), we have
\begin{align*}
\forall v \in \sphere^{d - 1},~\theta \in \real^d 
\quad 
\vecnorm{(\nabla^2 F(\theta) - \nabla^2 F(\thetastar))v}{2}^2 
\leq 
\Exs\vecnorm{(\nabla^2 f(\theta; \xi) - \nabla^2 f(\thetastar; \xi))v}{2}^2 
\leq 
\beta^2 \vecnorm{\theta - \thetastar}{2}^2
.
\end{align*}
Consequently, we have the bound:
\begin{align*}
\vecnorm{\nabla F(\theta) - \hessianstar (\theta - \thetastar)}{2} 
&\le
\vecnorm{\theta - \thetastar}{2} \cdot \sup_{\vecnorm{\theta' - \thetastar}{2}\le \vecnorm{\theta - \thetastar}{2}}~\sup_{ v \in \sphere^{d-1}} \vecnorm{ (\nabla^2 F(\theta') - \hessianstar) v}{2}
%\\&
\le
\beta \vecnorm{\theta - \thetastar}{2}^2
.
\end{align*}

By Eq.~\eqref{eq:theta-bound-in-proof-of-clt}, we
have $\sqrt{T} \vecnorm{\nabla F (\theta_T) - \hessianstar (\theta_T -
  \thetastar)}{2} \xrightarrow{p} 0$. Invoking Slutsky's theorem, this leads to $\sqrt{T}
\hessianstar (\theta_T - \thetastar)\xrightarrow{d} \mathcal{N}
\left(0, \SigStar + \Exs \left( \NoiseAplain (\thetastar)
\Lambda_\stpsz \NoiseAplain (\thetastar)^\top \right) \right)$, and
consequently,
\begin{align*}
\sqrt{T} (\theta_T - \thetastar) 
\xrightarrow{d} 
\mathcal{N} \left(0,
(\hessianstar)^{-1} \left( 
\SigStar + \Exs [ \NoiseAplain  (\thetastar) \Lambda_\stpsz \NoiseAplain (\thetastar)^\top ] 
\right) (\hessianstar)^{-1} \right)
,
\end{align*}
which finishes the proof.

\pb\section{Proofs of auxiliary lemmas in \S\ref{sec_proof,theo_rootHolder_single}, \S\ref{sec_proof,theo_rootHolder_multi} and \S\ref{subsec:proof-theo-other-bounds} } \label{sec:lem_holder}

\pb\subsection{Proof of Lemma~\ref{lemm_finitebddHolder}}\label{sec_proof,lemm_finitebddHolder}
Recall that we have the recursive update rule of $z_t$ as 
\begin{align*}
tz_t = (t - 1) z_{t - 1} + (t - 1) ( \noiset(\theta_{t - 1}) - \noiset(\theta_{t - 2})) + \noiset(\theta_{t - 1}) 
.
\end{align*}
Taking fourth moments on both sides, we have
\begin{align}\label{eq:tzt_decomp_fourth}
\Exs \normd{tz_t} 
	&=
\Exs \normd{(t - 1) z_{t - 1} + (t - 1) ( \noiset(\theta_{t - 1}) - \noiset(\theta_{t - 2})) + \noiset(\theta_{t - 1})}
	\notag \\&=
\Exs \normd{(t - 1) z_{t - 1}} + \Exs \normd{ (t - 1) ( \noiset(\theta_{t - 1}) - \noiset(\theta_{t - 2})) + \noiset(\theta_{t - 1})} 
	\notag \\&\quad\,%%%%
+4 \Exs \normc{ (t - 1) ( \noiset(\theta_{t - 1}) - \noiset(\theta_{t - 2})) + \noiset(\theta_{t - 1})} \norm{(t - 1) z_{t - 1}}
	\notag \\&\quad\,%%%%
+6 \Exs \normb{ (t - 1) ( \noiset(\theta_{t - 1}) - \noiset(\theta_{t - 2})) + \noiset(\theta_{t - 1})} \normb{(t - 1) z_{t - 1}}
,
\end{align} 
where one of the terms is zeroed out.
By H\"{o}lder's inequality and Young's inequality, we bound the third term and the fourth term of the RHS as 
\begin{align*}
& 
\Exs \normc{ (t - 1) ( \noiset(\theta_{t - 1}) - \noiset(\theta_{t - 2})) + \noiset(\theta_{t - 1})} \norm{(t - 1) z_{t - 1}}
	\notag \\&\leq 
\left( \Exs \normd{(t - 1) ( \noiset(\theta_{t - 1}) - \noiset(\theta_{t - 2})) + \noiset(\theta_{t - 1})}\right)^{3/4} \left( \Exs \normd{(t - 1) z_{t - 1}} \right)^{1/4} 
	\notag \\&\leq 
\frac12 \left( \Exs \normd{(t - 1) ( \noiset(\theta_{t - 1}) - \noiset(\theta_{t - 2})) + \noiset(\theta_{t - 1})}\right)^{1/2} \left( \Exs \normd{(t - 1) z_{t - 1}} \right)^{1/2} 
	\notag \\&~\quad + \frac12 \Exs \normd{ (t - 1) ( \noiset(\theta_{t - 1}) - \noiset(\theta_{t - 2})) + \noiset(\theta_{t - 1})}
,
\end{align*}
and 
\begin{align*}
& 
\Exs \normb{ (t - 1) ( \noiset(\theta_{t - 1}) - \noiset(\theta_{t - 2})) + \noiset(\theta_{t - 1})} \normb{(t - 1) z_{t - 1}}
	\notag \\&~\leq 
\left( \Exs \normd{(t - 1) ( \noiset(\theta_{t - 1}) - \noiset(\theta_{t - 2})) + \noiset(\theta_{t - 1})}\right)^{1/2} \left( \Exs \normd{(t - 1) z_{t - 1}} \right)
.
\end{align*}
Thus Eq.~\eqref{eq:tzt_decomp_fourth} continues as 
\begin{align*}
\Exs \normd{tz_t} 
	&\leq 
\Exs \normd{(t - 1) z_{t - 1} + (t - 1) ( \noiset(\theta_{t - 1}) - \noiset(\theta_{t - 2})) + \noiset(\theta_{t - 1})}
	\notag \\&=
\Exs \normd{(t - 1) z_{t - 1}} + 3\Exs \normd{ (t - 1) ( \noiset(\theta_{t - 1}) - \noiset(\theta_{t - 2})) + \noiset(\theta_{t - 1})}
	\notag \\&\quad\,%%%%
+8 \left( \Exs \normd{(t - 1) ( \noiset(\theta_{t - 1}) - \noiset(\theta_{t - 2})) + \noiset(\theta_{t - 1})}\right)^{1/2} \left( \Exs \normd{(t - 1) z_{t - 1}} \right)
	\notag \\&\leq 
\left( \sqrt{\Exs \normd{(t - 1) z_{t - 1}}} + 4\sqrt{\Exs \normd{ (t - 1) ( \noiset(\theta_{t - 1}) - \noiset(\theta_{t - 2})) + \noiset(\theta_{t - 1})}} \right)^2
,
\end{align*}
where 
\begin{align*}
& 
\Exs \normd{ (t - 1) ( \noiset(\theta_{t - 1}) - \noiset(\theta_{t - 2})) + \noiset(\theta_{t - 1})}
	\notag \\&\leq 
27 (t -1 )^4 \Exs \normd{\noiset(\theta_{t - 1}) - \noiset(\theta_{t - 2})} + 27 \Exs \normd{\noiset(\theta_{t - 1}) - \noiset(\thetastar)} + 27 \Exs \normd{\noiset(\thetastar)} 
	\notag \\&\leq 
27\sgliptild^4 \stpsz^4  (t -1 )^4 \Exs \normd{v_{t - 1}} + \frac{27 \sgliptild^4}{\strongconvex^4} \Exs \normd{\naF{\theta_{t - 1}}} + 27 \sigstartild^4
.
\end{align*}
Then 
\begin{align*}
t^2 \sqrt{\Exs \normd{z_t} }
	\leq
\sqrt{\Exs \normd{(t - 1) z_{t - 1}}} + 12 \sqrt{3} \sgliptild^2 \stpsz^2 \sqrt{\Exs \normd{(t - 1) v_{t - 1}}} + \frac{12 \sqrt{3} \sgliptild^2}{\strongconvex^2} \sqrt{\Exs \normd{\naF{\theta_{t - 1}}}} + 12 \sqrt{3} \sigstartild^2
.
\end{align*}

Combining this with Eq.~\eqref{eq_vtrecursive_holder} in Lemma~\ref{lemm_vtrecursive_holder} that
\begin{align*}
\sqrt{\Exs \normd{tv_t}}
	\le
\left( 1 - \frac{\stpsz\strongconvex}{2} \right) \sqrt{\Exs \normd{(t - 1) v_{t - 1}}} + \frac{5}{\stpsz \strongconvex} \sqrt{\Exs \normd{\naF{\theta_{t - 1}}}} + 14\tilde{\sigstar}^2 
.
\end{align*}
By the choice of $\stpsz$ satisfying $\stpsz \le \frac{1}{56\smoothness} \wedge \frac{\strongconvex}{64\sgliptild^2}$, we have $\frac{\sgliptild^2}{\strongconvex^2} \leq \frac{1}{64 \stpsz \strongconvex}$ and 
\begin{align*}
t^2 \sqrt{\Exs \normd{z_t}} + t^2 \sqrt{\Exs \normd{v_t}} 
	\leq 
\sqrt{\Exs \normd{(t - 1) z_{t - 1}}} + \sqrt{\Exs \normd{(t - 1) v_{t - 1}}} + \frac{6}{\stpsz\strongconvex} \sqrt{\Exs \normd{\naF{\theta_{t - 1}}}} + 35 \sigstartild^2 
.
\end{align*}
Recursively applying the above inequality and by observing that $\sqrt{\Exs \normd{\naF{\theta_{t - 1}}}} \leq 2\sqrt{\Exs \normd{z_t}} +  2\sqrt{\Exs \normd{v_t}} $, we have
\begin{align}\label{eq:F_bound_first}
&\quad
\niters^2 \sqrt{\Exs \normd{\naF{\theta_{\niters - 1}}}} 
\leq
2 \niters^2 \sqrt{\Exs \normd{z_\niters}} + 2 \niters^2 \sqrt{\Exs \normd{v_\niters}}
\notag \\&\leq 
2 \sqrt{\Exs \normd{\burnin z_{\burnin}}} + 2 \sqrt{\Exs \normd{\burnin v_{\burnin}}} 
+
\frac{12}{\stpsz \strongconvex} \sum_{t = \burnin + 1}^\niters \sqrt{\Exs \normd{\naF{\theta_{t - 1}}}}
+ 
70 (\niters - \burnin) \sigstartild^2
.
\end{align} 
Further for $\burnin z_{\burnin}$ and $\burnin v_{\burnin}$ we note that by applying Khintchine's inequality as well as Young's inequality we have

\begin{align}\label{eq:z_init_bound}
\Exs \normd{\burnin z_{\burnin}}
	&=
\Exs \normd{\sum_{t = 1}^{\burnin} \noiset(\theta_0)}
	\le
\Exs \left(\sum_{t = 1}^{\burnin} \normb{\noiset(\theta_0)} \right)^2
    \leq 
 \burnin \Exs \sum_{t = 1}^{\burnin} \normd{\noiset(\theta_0)}
    \leq 
8 \burnin^2 \left(\frac{\sgliptild^4}{\strongconvex^4} \Exs \normd{\naF{\theta_0}}
    +
\sigstartild^4
\right)
,
\end{align}
and 
\begin{align}\label{eq:v_init_bound}
\Exs \normd{\burnin v_{\burnin}} 
    &=
\Exs \normd{\burnin z_{\burnin}} + \Exs \normd{\burnin \naF{\theta_0}} + 4 \Exs \normc{\burnin z_{\burnin}} \norm{\burnin \naF{\theta_0}} + 6 \Exs \normb{ \burnin z_{\burnin}} \normb{\burnin \naF{\theta_0}}
    \notag \\&\leq 
7 \Exs \normd{\burnin v_{\burnin}} 
    +
5 \Exs \normd{\burnin \naF{\theta_0}}
    \leq 
56 \burnin^2 \left(\frac{\sgliptild^4}{\strongconvex^4} \Exs \normd{\naF{\theta_0}} +  \sigstartild^4 \right) + 5 \burnin^4 \Exs \normd{\naF{\theta_0}}
.
\end{align}
Taking squared root on Eq.~\eqref{eq:z_init_bound} and~\eqref{eq:v_init_bound} and recalling that $\stpsz \leq \frac{\strongconvex}{64\sgliptild^2}$, we have 
\begin{align}\label{eq:z_bound_holder}
\sqrt{\Exs \normd{\burnin z_{\burnin}}} 
    &\leq 2\sqrt{2} \burnin \left(\frac{\sgliptild^2}{\strongconvex^2} \sqrt{\Exs \normd{\naF{\theta_0}}} + \sigstartild^2 \right)
,
\end{align}
and
\begin{align}\label{eq:v_bound_holder}
\sqrt{\Exs \normd{\burnin v_{\burnin}}}
    &\leq 
(\sqrt{5} + 1/8) \burnin^2 \sqrt{\Exs \normd{\naF{\theta_0}}} + 8 \burnin \sigstartild^2
.
\end{align}
Bringing Eq.~\eqref{eq:z_bound_holder} and~\eqref{eq:v_bound_holder} into Eq.~\eqref{eq:F_bound_first}, we arrive at the following:
\begin{align}\label{eq:FT_recursive}
 T^2 \sqrt{\Exs \normd{\naF{\theta_{\niters - 1}}}} 
	&\leq 
4\sqrt{2} \burnin \left(\frac{\sgliptild^2}{\strongconvex^2} \sqrt{\Exs \normd{\naF{\theta_0}}} + \sigstartild^2 \right)
    +
(2\sqrt{5} + 1/4) \burnin^2 \sqrt{\Exs \normd{\naF{\theta_0}}} 
    \notag \\&~\quad 
    + 
16 \burnin \sigstartild^2
	+
\frac{12}{\stpsz \strongconvex} \sum_{t = \burnin + 1}^\niters \sqrt{\Exs \normd{\naF{\theta_{t - 1}}}} 
	+
 70 (\niters - \burnin) \sigstartild^2
 	\notag \\&\leq 
5 \burnin^2
\sqrt{\Exs \normd{\naF{\theta_0}}} 
		+
\frac{12}{\stpsz \strongconvex} \sum_{t = \burnin + 1}^\niters \sqrt{\Exs \normd{\naF{\theta_{t - 1}}}} 
	+
70 T \sigstartild^2
.
\end{align}
Dividing both sides by $\niters^2$, summing up Eq.~\eqref{eq:FT_recursive} from $\niters = \burnin + 1$ to $T^* \geq \burnin + 1$ and using the fact that $\stpsz \leq \frac{\strongconvex}{64\sgliptild^2}, \burnin \geq 2$, we have
\begin{align*}
\sum_{\niters = \burnin + 1}^{T^*} \sqrt{\Exs \normd{\naF{\theta_{\niters - 1}}}} 
	&\leq 
 5 \burnin 
\sqrt{\Exs \normd{\naF{\theta_0}}} 
		+
\frac{12}{\stpsz \strongconvex \burnin} \sum_{t = \burnin + 1}^{T^*} \sqrt{\Exs \normd{\naF{\theta_{t - 1}}}} 
	+
70 \sigstartild^2  \log \left(\frac{T^*}{\burnin}\right)
.
\end{align*}
Taking $\burnin = \left\lceil\frac{24}{\stpsz\strongconvex}\right\rceil$, we have 
\begin{align*}
 \sum_{\niters = \burnin + 1}^{T^*} \sqrt{\Exs \normd{\naF{\theta_{\niters - 1}}}}  
 	\leq 
10 \burnin \sqrt{\Exs\normd{\naF{\theta_0}}} 
	+
140 \sigstartild^2 \log \left(\frac{T^*}{\burnin}\right)
.
\end{align*}
Again by Eq.~\eqref{eq:FT_recursive}, we have
\begin{align*}
&~\quad T^2 \sqrt{\Exs \normd{\naF{\theta_{\niters - 1}}}} 
	\notag \\&\leq 
5 \burnin^2 \sqrt{\Exs \normd{\naF{\theta_0}}} 
	 +
\frac{12}{\stpsz\strongconvex} \left( 
10 \burnin \sqrt{\Exs\normd{\naF{\theta_0}}} 
	+
140 \sigstartild^2 \log \left(\frac{T}{\burnin}\right)
\right) 
	+
70 T \sigstartild^2
	\notag \\&\leq 
10 \burnin^2 \sqrt{\Exs \normd{\naF{\theta_0}}} 
	+
70 \burnin \sigstartild^2 \log \left(\frac{T}{\burnin}\right) + 70 T \sigstartild^2
.
\end{align*}
Dividing both sides by $\niters^2$ we conclude that 
\begin{align*}
\sqrt{\Exs\normd{\naF{\theta_{\niters - 1}}}} 
	&\leq 
 \frac{10 \burnin^2}{T^2} \sqrt{\Exs \normd{\naF{\theta_0}}} + 70 \left(1 + \frac{\burnin}{T} \log \left( \frac{T}{\burnin} \right)  \right)\frac{\sigstartild^2}{T}
    \notag \\&\leq 
 \frac{10 \burnin^2}{T^2} \sqrt{\Exs \normd{\naF{\theta_0}}} + 
 \frac{140\sigstartild^2}{T}
.
\end{align*}
which finishes our proof of Lemma~\ref{lemm_finitebddHolder}.

\pb\subsection{Proof of Lemma \ref{lemm_vtsharpbdd_holder}}\label{sec_proof,lemm_vtsharpbdd_holder}
Our main technical tools is the following lemma, which bound the fourth moment of the $v_t$ recursion.

\begin{lemma}\label{lemm_vtrecursive_holder}%Lemma A%lemm_vtrecursive
Under the setting of Proposition \ref{theo_rootHolder_single}, when $
	\stpsz \le \frac{1}{56\smoothness} \wedge \frac{\strongconvex}{64\sgliptild^2}
$,
 we have the following bound for $t\ge \burnin+1$
\begin{align}\label{eq_vtrecursive_holder}
\sqrt{\Exs \normd{tv_t}}
	\le
\left( 1 - \frac{\stpsz\strongconvex}{2} \right) \sqrt{\Exs \normd{(t - 1) v_{t - 1}}} + \frac{5}{\stpsz \strongconvex} \sqrt{\Exs \normd{\naF{\theta_{t - 1}}}} + 14\tilde{\sigstar}^2 
.
\end{align}
\end{lemma}
\noindent
The detailed proof is relegated to \S\ref{sec_proof,lemm_vtrecursive_holder}.
We are ready for the proof of Lemma \ref{lemm_vtsharpbdd_holder}.
Indeed, from \eqref{finitebdd_holder} and \eqref{eq_vtrecursive_holder}
\begin{align}\label{vtrec_holder}
t^2\sqrt{\Exs\vecnorm{v_t}{2}^4}
&\le%%%%
\left( 1 - \frac{\stpsz \strongconvex}{2} \right) (t-1)^2\sqrt{\Exs\vecnorm{v_{t - 1}}{2}^4}
+
\frac{5}{\stpsz\strongconvex} 
\left[
\frac{60 \; \vecnorm{\nabla F(\theta_0)}{2}^2}{\stpsz^2 \strongconvex^2 t^2} + \frac{140 \; \sigstartild}{t}
\right]
+
14 \sigstartild^2
%%%%%%%%
\notag \\&\le
\left( 1 - \frac{\stpsz \strongconvex}{2} \right) (t-1)^2\sqrt{\Exs\vecnorm{v_{t - 1}}{2}^4}
+
\frac{310\; \vecnorm{\nabla F(\theta_0)}{2}^2}{\stpsz^3 \strongconvex^3 t^2} 
+ 
714 \sigstartild^2
.
\end{align}
We have from \eqref{vtrec_holder}
\begin{align*}
t^4\sqrt{\Exs\vecnorm{v_t}{2}^4}
&\le%%%%
\left( 1 - \frac{\stpsz \strongconvex}{2} \right) t^2(t-1)^2\sqrt{\Exs\vecnorm{v_{t - 1}}{2}^4}
+
\frac{310 \; \vecnorm{\nabla F(\theta_0)}{2}^2}{\stpsz^3 \strongconvex^3}
+ 
714 \sigstartild^2 t^2
%%%%%%%%
\notag \\&\le
\left( 1 - \frac{\stpsz \strongconvex}{6} \right) (t-1)^4\sqrt{\Exs\vecnorm{v_{t - 1}}{2}^4}
+
\frac{310 \; \vecnorm{\nabla F(\theta_0)}{2}^2}{\stpsz^3 \strongconvex^3}
+ 
714 \sigstartild^2 t^2
,
\end{align*}
since the following holds $
\frac{ t^2}{ (t-1)^2}
\le
\frac{1 - \frac{\stpsz \strongconvex}{6}}{(1 - \frac{\stpsz \strongconvex}{6})^3}
\le
\frac{1 - \frac{\stpsz \strongconvex}{6}}{1 - \frac{\stpsz \strongconvex}{2}}
$
This gives, by solving the recursion,
\begin{align}\label{vt_approachone}
&%%%%
\niters^4\sqrt{\Exs\vecnorm{v_{\niters}}{2}^4}
\le
\left( 1 - \frac{\stpsz \strongconvex}{6} \right)^{\niters-\burnin} \burnin^4\sqrt{\Exs\vecnorm{v_{\burnin}}{2}^4}
+
\sum_{t=\burnin+1}^\niters
\left( 1 - \frac{\stpsz \strongconvex}{6} \right)^{\niters-t}\left[
\frac{310 \; \vecnorm{\nabla F(\theta_0)}{2}^2}{\stpsz^3 \strongconvex^3} 
+ 
714 \sigstartild^2 t^2
\right]
%%%%%%%%
\notag \\&\le
\left( 1 - \frac{\stpsz \strongconvex}{6} \right)^{\niters-\burnin} \burnin^4\sqrt{\Exs\vecnorm{v_{\burnin}}{2}^4}
+
\sum_{t=\burnin+1}^\niters
\left( 1 - \frac{\stpsz \strongconvex}{6} \right)^{\niters-t} \frac{310 \; \vecnorm{\nabla F(\theta_0)}{2}^2}{\stpsz^3 \strongconvex^3} 
+ 
\sum_{t=\burnin+1}^\niters
\left( 1 - \frac{\stpsz \strongconvex}{6} \right)^{\niters-t} 714 \sigstartild^2 t^2
%%%%%%%%
\notag \\&\le
\left( 1 - \frac{\stpsz \strongconvex}{6} \right)^{\niters-\burnin} \burnin^4\sqrt{\Exs\vecnorm{v_{\burnin}}{2}^4}
+
\frac{6}{\stpsz \strongconvex}
\cdot
\frac{310 \; \vecnorm{\nabla F(\theta_0)}{2}^2}{\stpsz^3 \strongconvex^3}
+ 
\frac{6}{\stpsz \strongconvex} \niters^2 \cdot 714 \sigstartild^2
%%%%%%%%
\notag \\&\le
\left( 1 - \frac{\stpsz \strongconvex}{6} \right)^{\niters-\burnin} \burnin^4\sqrt{\Exs\vecnorm{v_{\burnin}}{2}^4}
+
\frac{187500 \; \vecnorm{\nabla F(\theta_0)}{2}^2}{\stpsz^4 \strongconvex^4}
+
\frac{(714)(6)\sigstartild^2}{\stpsz \strongconvex} \niters^2
%%%%%%%%
.
\end{align}
where the summand is increasing so
\begin{align*}
\sum_{t=\burnin+1}^\niters \left( 1 - \frac{\stpsz \strongconvex}{6} \right)^{\niters-t} t^2
&\le
\frac{6}{\stpsz \strongconvex} \niters^2
.
\end{align*}

All in all, this concludes
\begin{align*}
&%%%%
\sqrt{\Exs\vecnorm{v_{\niters}}{2}^4}
\le
\left( 1 - \frac{\stpsz \strongconvex}{6} \right)^{\niters-\burnin} \frac{\burnin^4}{\niters^4}\sqrt{\Exs\vecnorm{v_{\burnin}}{2}^4}
+
\frac{187500 \; \vecnorm{\nabla F(\theta_0)}{2}^2}{\stpsz^4 \strongconvex^4\niters^4}
+
\frac{(714)(6)\sigstartild^2}{\stpsz \strongconvex\niters^2}
%%%%%%%%
\notag \\&\le
\left( 1 - \frac{\stpsz \strongconvex}{6} \right)^{\niters-\burnin} \frac{\burnin^4}{\niters^4}\sqrt{\Exs\vecnorm{v_{\burnin}}{2}^4}
+
\frac{187500 \; \vecnorm{\nabla F(\theta_0)}{2}^2}{\stpsz^4 \strongconvex^4\niters^4}
+
\frac{(714)(6)\sigstartild^2}{\stpsz \strongconvex\niters^2}
.
%%%%%%%%
\end{align*}
Bringing the burn-in upper bounds~\eqref{eq:v_bound_holder}, we arrive at our final result for bounding $\sqrt{\Exs \normd{v_\niters}}$:
\begin{align*}
\sqrt{\Exs \normd{v_\niters}} 
    &\leq 
\left( 
\frac{3 \burnin^4}{\niters^4} \sqrt{\Exs \normd{\naF{\theta_0}}} + \frac{8 \burnin^3}{\niters^4} \sigstartild^2
\right)
+
\frac{187500 \; \vecnorm{\nabla F(\theta_0)}{2}^2}{\stpsz^4 \strongconvex^4\niters^4}
+
\frac{(714)(6)\sigstartild^2}{\stpsz \strongconvex\niters^2}
    \notag \\&\leq 
\frac{1359375 \normb{\naF{\theta_0}}}{\stpsz^4 \strongconvex^4 \niters^4} + \frac{4484 \sigstartild^2}{\stpsz \strongconvex \niters^2}
.
\end{align*}

\pb\subsection{Proofs of recursive bounds on $v_t$}
In this section, we prove Lemmas~\ref{lemm_vtrecursive_holder} and~\ref{lemm_vtrecursive}, the two recursive bounds for $\{v_t\}_{t \geq \burnin}$ used in the proof of main theorems.

\pb\subsubsection{Proof of Lemma \ref{lemm_vtrecursive_holder}}\label{sec_proof,lemm_vtrecursive_holder}
By definition, we note that:
\begin{align*}
tv_t
=
\left(t-1\right) (v_{t-1} + \nabla f(\theta_{t-1};\xi_t) - \nabla f(\theta_{t - 2};\xi_t ) )
+
\nabla f(\theta_{t-1};\xi_t)
.
\end{align*}
Subtracting off a $\nabla F(\theta_{t-1})$ term from both sides we have
\begin{align*}
tv_t - \nabla F(\theta_{t-1})
=
\left(t-1\right) (v_{t-1} + \nabla f(\theta_{t-1};\xi_t) - \nabla f(\theta_{t - 2};\xi_t ) )
+
\underbrace{
\nabla f(\theta_{t-1};\xi_t) - \nabla F(\theta_{t-1})
}_{
=\noise_t(\theta_{t-1})
}
.
\end{align*}
Taking the fourth moments on both sides, we have
\begin{align}\label{vt_base_holder}
&%%%%
\Exs \norm{t v_t - \naF{\theta_{t - 1}}}^4
%%%%%%%%%
\notag \\&=
\Exs \norm{
(t-1)v_{t-1} + (t-1) (\nabla f(\theta_{t-1};\xi_t) - \nabla f(\theta_{t - 2};\xi_t ) ) 
+
\noise_t(\theta_{t-1})}^4	
%%%%%%%%
\notag \\&=
(t-1)^4\Exs \norm{v_{t - 1}}^4
	+
\underbrace{4 \Exs \left[
\norm{(t - 1) v_{t - 1}}^2 \binprod{ (t - 1) v_{t - 1}}{(t - 1)(\nabla f(\theta_{t - 1}; \xi_t) - \nabla f(\theta_{t - 2}; \xi_t)) + \noiset(\theta_{t - 1}) }
\right]}_{=: T_2}
	\notag \\&\quad\,%%%%
+\underbrace{6  \Exs \left[ 
\normb{(t - 1)v_{t - 1}} \normb{(t - 1)(\nabla f(\theta_{t - 1}; \xi_t) - \nabla f(\theta_{t - 2}; \xi_t)) + \noiset(\theta_{t - 1}) }
\right]}_{=: T_1}
	\notag \\&\quad\,%%%%
+\underbrace{4 \Exs \left[
\normc{(t - 1)(\nabla f(\theta_{t - 1}; \xi_t) - \nabla f(\theta_{t - 2}; \xi_t)) + \noiset(\theta_{t - 1}) } \norm{(t - 1) v_{t - 1}}
\right]}_{=: T_3}
	\notag \\&\quad\,%%%%
+\Exs \left[ 
\norm{(t - 1)(\nabla f(\theta_{t - 1}; \xi_t) - \nabla f(\theta_{t - 2}; \xi_t)) + \noiset(\theta_{t - 1}) }^4
\right]
.
\end{align}
To bound term $T_1$, we apply the H\"{o}lder's inequality and have
\begin{align}\label{eq:bound_I_holder}
T_1
\leq 
6 \left( \Exs \norm{(t - 1) v_{t - 1}}^4 \right)^{1/2} \left( \Exs \norm{(t - 1)(\nabla f(\theta_{t - 1}; \xi_t) - \nabla f(\theta_{t - 2}; \xi_t)) + \noiset(\theta_{t - 1}) }^4 \right)^{1/2}
.
\end{align}
To bound term $T_3$, we again apply the H\"{o}lder's inequality:
\begin{align}\label{eq:bound_III_holder}
T_3
&\leq 
4 \left(\Exs \normd{(t - 1) ( \nafxit{\theta_{t - 1}} - \nafxit{\theta_{t - 2}}) + \noiset(\theta_{t - 1})} \right)^{3/4} 
\left( \Exs \normd{(t - 1) v_{t - 1}} \right)^{1/4} 
	\notag \\&\leq 
2 \left(\Exs \normd{(t - 1) ( \nafxit{\theta_{t - 1}} - \nafxit{\theta_{t - 2}}) + \noiset(\theta_{t - 1})} \right)^{1/2} 
\left( \Exs \normd{(t - 1) v_{t - 1}} \right)^{1/2} 
	\notag \\&\quad\,%%%%
+2 \Exs \normd{ (t - 1) ( \nafxit{\theta_{t - 1}} - \nafxit{\theta_{t - 2}}) + \noiset(\theta_{t - 1})}
.
\end{align} 
To bound term $T_2$, we first take expectation with respect to $\xi_t$ and have 
\begin{align*}
T_2 &= 4 \Exs \left[ 
\normb{(t - 1)v_{t - 1}} \binprod{
(t - 1) v_{t - 1} }{
(t - 1) \left( \nabla F(\theta_{t - 1}) - \nabla F(\theta_{t - 2}) \right) } 
\right]
,
\end{align*}
where 
\begin{align*}
\binprod{v_{t - 1}}{\naF{\theta_{t - 1}} - \naF{\theta_{t - 2}}}
	\leq 
- \frac{1}{\stpsz L} \normb{\naF{\theta_{t - 1}} - \naF{\theta_{t - 2}}} 
\end{align*}
and 
\begin{align*}
\binprod{v_{t - 1}}{\naF{\theta_{t - 1}} - \naF{\theta_{t - 2}}}
	\leq 
- \frac{\strongconvex}{\stpsz} \normb{\theta_{t - 1} - \theta_{t - 2}} 
\end{align*}
holds true for any $\strongconvex$-strongly convex and $L$-smooth $F$. Then we have
\begin{align}\label{eq:bound_II_holder}
T_2 
&\leq  
- (t - 1)^4\Exs \left[ 
\normb{v_{t - 1}} \left(
\frac{1}{\stpsz L} \normb{\naF{\theta_{t - 1}} - \naF{\theta_{t - 2}}} 
	+
3 \stpsz \strongconvex \normb{v_{t - 1}} 
\right)
\right]
	\notag \\&= 
- 3 \stpsz \strongconvex (t - 1)^4 \Exs \normd{v_{t - 1}}  - \frac{(t - 1)^4}{\stpsz L} \Exs \normb{v_{t - 1}} \normb{\naF{\theta_{t - 1}} - \naF{\theta_{t - 2}}}
	\notag \\&\leq 
- 3 \stpsz \strongconvex (t - 1)^4 \Exs \normd{v_{t - 1}}  - \frac{(t - 1)^4}{\stpsz L} \left( \Exs \normd{v_{t - 1}} \right)^{1/2} \left( \Exs \normd{\naF{\theta_{t - 1}} - \naF{\theta_{t - 2}}}\right)^{1/2}
.
\end{align}
Combining Eqs.~\eqref{eq:bound_I_holder},~\eqref{eq:bound_III_holder} and~\eqref{eq:bound_II_holder} into Eq.~\eqref{vt_base_holder} we have 
\begin{align}\label{eq:vt_mid_holder} 
&~\quad \Exs \normd{tv_t - \naF{\theta_{t - 1}}} 
	\notag \\&\leq 
(t - 1)^4 \Exs \normd{v_{t - 1}} 
	+
3 \Exs \normd{(t - 1) ( \nafxit{\theta_{t - 1}} - \nafxit{\theta_{t - 2}}) + \noiset(\theta_{t - 1})}
	\notag \\&\quad\,%%%%
+8 \left( \Exs \normd{(t - 1) ( \nafxit{\theta_{t - 1}} - \nafxit{\theta_{t - 2}}) + \noiset(\theta_{t - 1})} \right)^{1/2} \left(\Exs \normd{(t - 1) v_{t - 1}} \right)^{1/2}
	\notag \\&~\quad 
- 3 \stpsz \strongconvex (t - 1)^4 \Exs \normd{v_{t - 1}}  - \frac{(t - 1)^4}{\stpsz L} \left( \Exs \normd{v_{t - 1}} \right)^{1/2} \left( \Exs \normd{\naF{\theta_{t - 1}} - \naF{\theta_{t - 2}}}\right)^{1/2}
.
\end{align}

We now turn to bound the term $\Exs \normd{(t - 1) (\nafxit{\theta_{t - 1}} - \nafxit{\theta_{t - 2}}) + \noiset(\theta_{t - 1})}$ by the following decomposition scheme:
\begin{align}\label{eq:vt_decomp_decomp_holder}
&~\quad \Exs \normd{(t - 1) (\nafxit{\theta_{t - 1}} - \nafxit{\theta_{t - 2}}) + \noiset(\theta_{t - 1})}
	\notag \\&\leq 
\Exs \normd{(t - 1) (\naF{\theta_{t - 1}} - \naF{\theta_{t - 2}}) + (t - 1) ( \noiset(\theta_{t - 1}) - \noiset(\theta_{t - 2})) + \noiset(\theta_{t - 1}) - \noiset(\theta^*) + \noiset(\theta^*)}
	\notag \\&\leq 
8 (t - 1)^4 \underbrace{\Exs \normd{\naF{\theta_{t - 1}} - \naF{\theta_{t - 2}} + \noiset(\theta_{t - 1}) - \noiset(\theta_{t - 2})}}_{=: I_1} 
	+
8 \underbrace{\Exs \normd{\noiset(\theta_{t - 1}) - \noiset(\theta^*) + \noiset(\theta^*)}}_{=: I_2} 
.
\end{align}
We claim that 
\begin{align}\label{eq:I_1}
I_1 \leq
5\Exs \normd{\naF{\theta_{t - 1}} - \naF{\theta_{t - 2}}} + 7\sglipfour^4 \stpsz^4 \Exs \normd{v_{t - 1}}
,
\end{align} 
and 
\begin{align}\label{eq:I_2}
I_2 \leq
	 8 \Exs \norm{\noiset(\theta_{t - 1}) - \noiset(\theta^*)}^4 + 8 \Exs \norm{\noiset(\theta^*)}^4 
	\leq \frac{8 \sgliptild^4}{\strongconvex^4} \Exs \norm{\nabla F(\theta_{t - 1})}^4 + 8 \tilde{\sigstar}^4
.
\end{align}

Combining Eqs.~\eqref{eq:vt_decomp_decomp_holder},~\eqref{eq:I_1} and~\eqref{eq:I_2} we have the bound 
\begin{align}\label{eq:bound_mid}
&~\quad \Exs \normd{(t - 1) (\nafxit{\theta_{t - 1}} - \nafxit{\theta_{t - 2}}) + \noiset(\theta_{t - 1})} 
	\notag \\&\leq 
40(t - 1)^4 \Exs \normd{\naF{\theta_{t - 1}} - \naF{\theta_{t - 2}}} + 56 \sglipfour^4 \stpsz^4 (t - 1)^4 \Exs \normd{v_{t - 1}}
	+
\frac{64 \sgliptild^4}{\strongconvex^4} \Exs \norm{\nabla F(\theta_{t - 1})}^4 + 64 \tilde{\sigstar}^4
.
\end{align}

Then, we bring Eq.~\eqref{eq:bound_mid} into Eq.~\eqref{eq:vt_mid_holder} and have

\begin{align*}
&~\quad \Exs \normd{tv_t - \naF{\theta_{t - 1}}} 
	\notag \\&\leq 
(t - 1)^4 \Exs \normd{v_{t - 1}} + 
120(t - 1)^4 \Exs \normd{\naF{\theta_{t - 1}} - \naF{\theta_{t - 2}}} 
	+ 168 \sglipfour^4 \stpsz^4 (t - 1)^4 \Exs \normd{v_{t - 1}}
	\notag \\&\quad\,%%%%
+\frac{192 \sgliptild^4}{\strongconvex^4} \Exs \norm{\nabla F(\theta_{t - 1})}^4 
+ 192 \tilde{\sigstar}^4
	+
8 \sqrt{40} (t - 1)^4 \left( \Exs \normd{\naF{\theta_{t - 1}} - \naF{\theta_{t - 2}}} \right)^{1/2} \left( \Exs \normd{v_{t - 1}}\right)^{1/2}
	\notag \\&\quad\,%%%%
+64 \sglipfour^2 \stpsz^2 (t - 1)^4 \Exs \normd{v_{t - 1}} 
	+
64\left( 
\frac{ \sgliptild^4}{\strongconvex^4} \Exs \norm{\nabla F(\theta_{t - 1})}^4 + \tilde{\sigstar}^4
\right)^{1/2} \left( \Exs \normd{v_{t - 1}}\right)^{1/2}
	\notag \\&~\quad 
- 3 \stpsz \strongconvex (t - 1)^4 \Exs \normd{v_{t - 1}}  - \frac{(t - 1)^4}{\stpsz L} \left( \Exs \normd{v_{t - 1}} \right)^{1/2} \left( \Exs \normd{\naF{\theta_{t - 1}} - \naF{\theta_{t - 2}}}\right)^{1/2}
	\notag \\&\leq 
\left( 
1 - 3 \stpsz \strongconvex + 64 \sglipfour^2 \stpsz^2 + 168 \sglipfour^4 \stpsz^4 
\right) \Exs \normd{(t - 1) v_{t - 1}} 
	\notag \\&\quad\,%%%%
+\left( 
8 \sqrt{40} - \frac{1}{\stpsz L} + 120 L^2 \stpsz^2 
\right)(t - 1)^4 \left( \Exs \normd{v_{t - 1}} \right)^{1/2} \left( \Exs \normd{\naF{\theta_{t - 1}} - \naF{\theta_{t - 2}}}\right)^{1/2} 
	\notag \\&\quad\,%%%%
+64\left( 
\frac{ \sgliptild^4}{\strongconvex^4} \Exs \norm{\nabla F(\theta_{t - 1})}^4 + \tilde{\sigstar}^4
\right)^{1/2} \left( \Exs \normd{v_{t - 1}}\right)^{1/2}
	+
\frac{192 \sgliptild^4}{\strongconvex^4} \Exs \norm{\nabla F(\theta_{t - 1})}^4 
+ 192 \tilde{\sigstar}^4
	\notag \\&\leq 
\left( 1 - \stpsz \strongconvex\right)^2 \Exs \normd{(t - 1) v_{t - 1}}+
64\left( 
\frac{ \sgliptild^4}{\strongconvex^4} \Exs \norm{\nabla F(\theta_{t - 1})}^4 + \tilde{\sigstar}^4
\right)^{1/2} \left( \Exs \normd{v_{t - 1}}\right)^{1/2}
	\notag \\&\quad\,%%%%
+\frac{192 \sgliptild^4}{\strongconvex^4} \Exs \norm{\nabla F(\theta_{t - 1})}^4 
+ 192 \tilde{\sigstar}^4 
,
\end{align*}
where the last inequality is due to the choice of $\stpsz \leq \frac{\strongconvex}{64 \sgliptild^2}$ and  $\stpsz \leq \frac{1}{56 L} $ such that 
\begin{align*}
168 \sgliptild^4 \stpsz^4 \leq \stpsz^2 \strongconvex^2, \sgliptild^2 \stpsz^2 \leq \stpsz \strongconvex
\qquad \text{and}\quad
8 \sqrt{40} - \frac{1}{\stpsz L} + 120 L^2 \stpsz^2 \leq 0
.
\end{align*}
Taking squared root on both sides, we have
\begin{align}\label{eq:bound_vtdiff_holder}
\sqrt{\Exs \normd{tv_t - \naF{\theta_{t - 1}}}} 
	\leq 
(1 - \stpsz \strongconvex) \sqrt{\Exs \normd{(t - 1) v_{t - 1}}} + 32\left(\frac{\sgliptild^2}{\strongconvex^2}\sqrt{ \Exs \normd{\naF{\theta_{t - 1}}}} + \sigstartild^2\right)
.
\end{align} 

Furthermore, Young's inequality gives%
\footnote{Here, a different coefficient from the analysis as in the proof of Theorem~\ref{theo_finitebdd_single} is adopted.}
\begin{align*}
&~\quad \Exs \normd{tv_t - \naF{\theta_{t - 1}}} 
	\notag \\&= 
t^4 \Exs \normd{v_t} + \Exs \normd{\naF{\theta_{t - 1}}} 
	+
6 \Exs \normb{tv_t} \normb{\naF{\theta_{t - 1}}} 
	 -
4 \Exs \normc{tv_t}\norm{\naF{\theta_{t - 1}}} 
	-
4 \Exs \norm{tv_t} \normc{\naF{\theta_{t - 1}}}
	\notag \\&\geq 
t^4 \Exs \normd{v_t} + \Exs \normd{\naF{\theta_{t - 1}}} 
	+
6 \Exs \normb{tv_t} \normb{\naF{\theta_{t - 1}}} 
	 -
2 \Exs \left[\frac{2}{\stpsz\strongconvex} \normb{tv_t} \normb{\naF{\theta_{t - 1}}} + \frac{\stpsz\strongconvex}{2} \normd{tv_t} \right]
	\notag \\&~\quad-
2 \Exs\left[ \frac{\stpsz\strongconvex}{2}\normb{tv_t} \normb{\naF{\theta_{t - 1}}} + \frac{2}{\stpsz\strongconvex}\normd{\naF{\theta_{t - 1}}}\right]
	\notag \\&\geq 
\left(1 - \stpsz\strongconvex \right) \Exs \normd{t v_t} + \left( 1 - \frac{4}{\stpsz \strongconvex}\right) \Exs \normd{\naF{\theta_{t - 1}}} 
	+
\left( 6 - \frac{4}{\stpsz \strongconvex} - \stpsz\strongconvex \right) \Exs \normb{tv_t} \normb{\naF{\theta_{t - 1}}}
	\notag \\&\geq 
\left(1 - \stpsz\strongconvex \right) \Exs \normd{t v_t} - \left( \frac{4}{\stpsz \strongconvex} - 1\right) \Exs \normd{\naF{\theta_{t - 1}}} 
	-
\left(1 - \stpsz\strongconvex\right) \left( \frac{4}{\stpsz\strongconvex} - 1 \right)  \Exs \normb{tv_t} \normb{\naF{\theta_{t - 1}}}
.
\end{align*}
Combining this we have 
\begin{align*}
&\quad\,%%%%
\left(1 - \stpsz\strongconvex \right) \Exs \normd{t v_t} - \left( \frac{4}{\stpsz \strongconvex} - 1\right) \Exs \normd{\naF{\theta_{t - 1}}} 
	-
\left(1 - \stpsz\strongconvex\right) \left( \frac{4}{\stpsz\strongconvex} - 1 \right)  \Exs \normb{tv_t} \normb{\naF{\theta_{t - 1}}}
%%%%%%%%
\notag \\&\le
\left( 1 - \stpsz \strongconvex\right)^2 \Exs \normd{(t - 1) v_{t - 1}}+
64\left( 
\frac{ \sgliptild^4}{\strongconvex^4} \Exs \norm{\nabla F(\theta_{t - 1})}^4 + \tilde{\sigstar}^4
\right)^{1/2} \left( \Exs \normd{v_{t - 1}}\right)^{1/2}
	\notag \\&\quad\,%%%%
+\frac{192 \sgliptild^4}{\strongconvex^4} \Exs \norm{\nabla F(\theta_{t - 1})}^4 
+ 192 \tilde{\sigstar}^4 
.
\end{align*}

Now we multiply both sides by $
	\left(1 - \stpsz \strongconvex\right)^{-1}
$ , noting that $(1 - \stpsz \strongconvex)^{-1} \leq (1 - \stpsz L)^{-1} \leq \frac{56}{55}$, rearranging, and have
\begin{align*}
 \Exs \normd{t v_t}
%%%%%%%%
&\le
\left( 1 - \stpsz \strongconvex\right) \Exs \normd{(t - 1) v_{t - 1}}+
\frac{(56)(64)}{(55)}\left( 
\frac{ \sgliptild^4}{\strongconvex^4} \Exs \norm{\nabla F(\theta_{t - 1})}^4 + \tilde{\sigstar}^4
\right)^{1/2} \left( \Exs \normd{v_{t - 1}}\right)^{1/2}
	\notag \\&\quad\,%%%%
+\frac{\frac{(56)(192)}{(55)} \sgliptild^4}{\strongconvex^4} \Exs \norm{\nabla F(\theta_{t - 1})}^4 
+ \frac{(56)(192)}{(55)} \tilde{\sigstar}^4 
	\notag \\&\quad\,%%%%
+\frac{\frac{(4)(56)}{(55)}}{\stpsz\strongconvex} \Exs \normd{\naF{\theta_{t - 1}}} + \frac{4}{\stpsz\strongconvex} \Exs \normb{tv_t} \normb{\naF{\theta_{t - 1}}}
	\notag \\&\leq 
\left(1 - \frac{\stpsz\strongconvex}{2} \right)^2 \Exs \normd{(t - 1)v_{t - 1}} 
	+
\frac{7}{55\stpsz^2 \strongconvex^2} \Exs \normd{\naF{\theta_{t - 1}}} + \frac{(56) (192)}{55} \tilde{\sigstar}^4 
	\notag \\&~\quad + 
\frac{56}{55}\left( 
\frac{ 1}{64 \stpsz^2\strongconvex^2} \Exs \norm{\nabla F(\theta_{t - 1})}^4 + 64\tilde{\sigstar}^4
\right)^{1/2} \left( \Exs \normd{v_{t - 1}}\right)^{1/2} 
	+
\frac{4}{\stpsz\strongconvex} \Exs \normb{tv_t} \normb{\naF{\theta_{t - 1}}}
.
\end{align*}
Rearranging and taking squared root on both sides we conclude that 
\begin{align*}
\sqrt{\Exs \normd{tv_t}} - \frac{2}{\stpsz\strongconvex} \sqrt{\Exs \normd{\naF{\theta_{t - 1}}}}
	&\leq 
\left( 1 - \frac{\stpsz\strongconvex}{2} \right) \sqrt{\Exs \normd{(t - 1) v_{t - 1}}} + \frac{3}{\stpsz\strongconvex} \sqrt{\Exs \normd{\naF{\theta_{t - 1}}}} + 14\sigstartild^2
.
\end{align*}
Further rearranging, we have 
\begin{align*}
\sqrt{\Exs \normd{tv_t}} \leq 
\left( 1 - \frac{\stpsz\strongconvex}{2} \right) \sqrt{\Exs \normd{(t - 1) v_{t - 1}}} + \frac{5}{\stpsz \strongconvex} \sqrt{\Exs \normd{\naF{\theta_{t - 1}}}} + 14\tilde{\sigstar}^2 
,
\end{align*}
which concludes our proof.

\paragraph{Proof of Eq.~\eqref{eq:I_1}:}
We use similar decomposition as in the decomposition in Eq.~\eqref{vt_base_holder} and have
\begin{align*}
I_1 
	&=
\Exs \normd{\naF{\theta_{t - 1}} - \naF{\theta_{t - 2}}} + \Exs \normd{\noiset(\theta_{t - 1}) - \noiset(\theta_{t - 2})} 
	\notag \\&\quad\,%%%%
+4 \Exs \normc{\noiset(\theta_{t - 1}) - \noiset(\theta_{t - 2})} \norm{\naF{\theta_{t - 1}} - \naF{\theta_{t - 2}}} 
	\notag \\&\quad\,%%%%
+6 \Exs \normb{\naF{\theta_{t - 1}} - \naF{\theta_{t - 2}}} \normb{\noiset(\theta_{t - 1}) - \noiset(\theta_{t - 2})} 
,
\end{align*}
where we note that we used the fact that one of the cross terms in the fourth moment decomposition $
\Exs\left[
\normb{\naF{\theta_{t - 1}} - \naF{\theta_{t - 2}}} \binprod{\naF{\theta_{t - 1}} - \naF{\theta_{t - 2}}}{ \noiset(\theta_{t - 1}) - \noiset(\theta_{t - 2})}
\right] = 0
$.
Further utilizing the H\"{o}lder's inequality, we have
\begin{align*}
I_1
%%%%%%%%
&\leq 
\Exs \normd{\naF{\theta_{t - 1}} - \naF{\theta_{t - 2}}} + \Exs \normd{\noiset(\theta_{t - 1}) - \noiset(\theta_{t - 2})} 
	\notag \\&\quad\,%%%%
+
4 \left(\Exs \normd{\noiset(\theta_{t - 1}) - \noiset(\theta_{t - 2})}\right)^{3/4} \left(\Exs \norm{\naF{\theta_{t - 1}}- \naF{\theta_{t - 2}}} \right)^{1/4}
	\notag \\&\quad\,%%%%
+
6 \left(\Exs \normd{\naF{\theta_{t - 1}} - \naF{\theta_{t - 2}}}\right)^{1/2} \left( \Exs \normb{\noiset(\theta_{t - 1}) - \noiset(\theta_{t - 2})} \right)^{1/2}
%%%%%%%%
\notag \\&\leq 
\Exs \normd{\naF{\theta_{t - 1}} - \naF{\theta_{t - 2}}} + 3\Exs \normd{\noiset(\theta_{t - 1}) - \noiset(\theta_{t - 2})} 
	\notag \\&\quad\,%%%%
+
8 \left(\Exs \normd{\naF{\theta_{t - 1}} - \naF{\theta_{t - 2}}}\right)^{1/2} \left( \Exs \normd{\noiset(\theta_{t - 1}) - \noiset(\theta_{t - 2})} \right)^{1/2}
%%%%%%%%
\notag \\&\leq 
\Exs \normd{\naF{\theta_{t - 1}} - \naF{\theta_{t - 2}}} + 3 \sglipfour^4 \stpsz^4 \Exs \normd{v_{t - 1}} 
	\notag \\&\quad\,%%%%
+
8 \sglipfour^2 \stpsz^2 \left(\Exs \normd{\naF{\theta_{t - 1}} - \naF{\theta_{t - 2}}}\right)^{1/2} \left( \Exs \normd{v_{t - 1}} \right)^{1/2}
%%%%%%%%
\notag \\&\leq 
5\Exs \normd{\naF{\theta_{t - 1}} - \naF{\theta_{t - 2}}} + 7 \sglipfour^4 \stpsz^4 \Exs \normd{v_{t - 1}}
.
\end{align*}
This completes the proof of Eq.~\eqref{eq:I_1}.

\pb\subsubsection{Proof of Lemma~\ref{lemm_vtrecursive}}\label{sec_proof,lemm_vtrecursive}

By definition, we note that:
\begin{align*}
  v_t = \left( 1 - \frac{1}{t} \right) (v_{t - 1} + \nabla f
  (\theta_{t - 1};\xi_t) - \nabla f (\theta_{t - 2};\xi_t ) ) +
  \frac{1}{t} \nabla f ( \theta_{t - 1};\xi_t).
\end{align*}
Taking the second moments for both sides, we have:
\begin{align*}
  \Exs\vecnorm{v_t}{2}^2 = &\left(1 - \frac{1}{t}\right)^2 \underbrace{\Exs
    \vecnorm{v_{t - 1} + \nabla f (\theta_{t - 1};\xi_t) - \nabla f
      (\theta_{t - 2};\xi_t)}{2}^2}_{I_1} + \frac{1}{t^2}
  \underbrace{\Exs\vecnorm{\nabla f (\theta_{t -
        1};\xi_t)}{2}^2}_{I_2}\\ &+ 2 \frac{t - 1}{t^2}
    \underbrace{\Exs \inprod{ v_{t - 1} + \nabla f (\theta_{t - 1};\xi_t) - \nabla
      f (\theta_{t - 2};\xi_t) }{ \nabla f (\theta_{t -
        1};\xi_t)}}_{I_3}.
\end{align*}
For the first term, using the fact that $ \theta_{t - 1} - \theta_{t - 2} = - \stpsz v_{t -
  1} $, we start with the following decomposition:
\begin{align*}
    &\Exs \left( \vecnorm{v_{t - 1} + \nabla f (\theta_{t - 1};\xi_t)
    - \nabla f (\theta_{t - 2};\xi_t) }{2}^2 \mid \filtration_{t - 1}
  \right)\\ &= \vecnorm{v_{t - 1}}{2}^2 + 2 \Exs \left( \inprod{v_{t -
      1}}{ \nabla f (\theta_{t - 1};\xi_t) - \nabla f (\theta_{t -
      2};\xi_t) } \mid \filtration_{t - 1} \right) +  \Exs \left(
  \vecnorm{ \nabla f (\theta_{t - 1};\xi_t) - \nabla f (\theta_{t -
      2};\xi_t)}{2}^2 \mid \filtration_{t - 1} \right)\\
      &=
  \vecnorm{v_{t - 1}}{2}^2 - \frac{2}{\stpsz}
  \inprod{\theta_{t - 1} - \theta_{t - 2}}{\nabla F (\theta_{t - 1}) -
    \nabla F (\theta_{t - 2})} + \Exs \left(
  \vecnorm{ \nabla f (\theta_{t - 1};\xi_t) - \nabla f (\theta_{t -
      2};\xi_t)}{2}^2 \mid \filtration_{t - 1} \right).
\end{align*}
Since $F$ is $\strongconvex$-strongly convex and $\smoothness$-smooth, we have the following standard inequality:
\begin{align*}
    \inprod{\theta_{t - 1} - \theta_{t - 2}}{\nabla F (\theta_{t - 1}) -
    \nabla F (\theta_{t - 2})} \geq  \frac{\vecnorm{\theta_{t - 1} - \theta_{t - 2}}{2}^2 \strongconvex
    \smoothness}{\strongconvex + \smoothness} + \frac{\vecnorm{\nabla
      F (\theta_{t - 1}) - \nabla F (\theta_{t - 2})}{2}^2
  }{\strongconvex + \smoothness}.
\end{align*}
Hence, when the step-size satisfies the bound $\stpsz_t \leq
\frac{1}{2 L} \wedge \frac{\strongconvex}{2 \sglip^2}$, there is the bound:
\begin{align*}
I_1  
&\leq 
\Exs
\vecnorm{v_{t - 1}}{2}^2 - \frac{2}{\stpsz} \Exs \left(
\frac{\vecnorm{\theta_{t - 1} - \theta_{t - 2}}{2}^2 \strongconvex \smoothness}{\strongconvex + \smoothness} 
+ 
\frac{\vecnorm{\nabla F (\theta_{t - 1}) - \nabla F (\theta_{t - 2})}{2}^2}{\strongconvex + \smoothness} \right) 
+ 
2 \Exs\vecnorm{\nabla F (\theta_{t - 1}) - \nabla F ( \theta_{t - 2} )}{2}^2  
\\& \quad\quad 
+  2 \Exs \left( \vecnorm{\noise (\theta_{t - 1}; \xi_t) - \noise ( \theta_{t - 2}; \xi_t )}{2}^2 \right)
\\&\leq 
(1 - \stpsz \strongconvex + 2 \stpsz^2 \sglip^2 ) \Exs\vecnorm{v_{t - 1}}{2}^2 + 2 \left(1 - \frac{1}{\stpsz (\strongconvex + \smoothness)} \right)  \Exs\vecnorm{\nabla F (\theta_{t - 1}) - \nabla F (\theta_{t - 2})}{2}^2
\\&\leq 
(1 - 3\stpsz \strongconvex / 4) \Exs\vecnorm{v_{t - 1}}{2}^2
.
\end{align*}

Now we study the second term $I_2$, note that
\begin{align*}
\Exs\vecnorm{\nabla f (\theta_{t - 1};\xi_t)}{2}^2 
&\leq  
\Exs\vecnorm{\nabla f (\theta_{t - 1};\xi_t) - \nabla f(\thetastar;\xi_t)}{2}^2 + 4 \Exs\vecnorm{\nabla f(\thetastar;\xi_t)}{2}^2
\\&\leq 
2 \Exs\vecnorm{\nabla F(\theta_{t - 1})}{2}^2 + 2 \Exs\vecnorm{\varepsilon (\theta_{t - 1}; \xi_t) - \varepsilon (\thetastar; \xi_t)}{2}^2 + 4 \Exs\vecnorm{\nabla f (\thetastar;\xi_t)}{2}^2
\\ &\leq 
2 \Exs\vecnorm{\nabla F (\theta_{t - 1})}{2}^2 + 2 \sglip^2 \Exs\vecnorm{\theta_{t - 1} - \thetastar}{2}^2 + 4 \sigstar^2
\\ &\leq 
2 \left( 1 +  \frac{\sglip^2}{\strongconvex^2}\right) \Exs\vecnorm{\nabla F (\theta_{t - 1})}{2}^2 + 4 \sigstar^2
.
\end{align*}
For the cross term $I_3$, we note that:
\begin{align*}
&\quad\,%
\Exs \left( \inprod{ v_{t - 1} + \nabla f (\theta_{t - 1};\xi_t) - \nabla f (\theta_{t - 2};\xi_t) }{ \nabla f (\theta_{t - 1};\xi_t)} \mid \filtration_{t - 1} \right)
\\ &= 
\Exs \left( \inprod{v_{t - 1}}{\nabla f (\theta_{t - 1}; \xi_t)} \mid \filtration_{t - 1} \right) + \Exs \left( \inprod{\nabla f (\theta_{t - 1}; \xi_t) - \nabla f (\theta_{ t- 2}; \xi_t)}{\nabla F (\theta_{t - 1})} \mid \filtration_{t - 1} \right)
\\ &\quad \quad 
+ \Exs \left( \inprod{\nabla f (\theta_{t - 1}; \xi_t) - \nabla f (\theta_{ t- 2}; \xi_t)}{\noise_t (\theta_{t - 1})} \mid \filtration_{t - 1} \right) 
\\ &= 
\underbrace{\inprod{v_{t - 1}}{ \nabla F (\theta_{t - 1})} + \inprod{\nabla F (\theta_{t - 1}) - \nabla F (\theta_{t - 2})}{\nabla F (\theta_{t - 1})} }_{\mydefn T_1} 
\\&\quad\,
+ \underbrace{\Exs \left( \inprod{\noise(\theta_{t - 1}, \xi_t) - \noise (\theta_{t - 2}; \xi_t)}{\noise (\theta_{t - 1}, \xi_t)} \mid \filtration_{t - 1} \right) }_{\mydefn T_2}
.
\end{align*}
For the term $T_1$, we note that:
\begin{align*}
    T_1 \leq \vecnorm{v_{t - 1}}{2} \cdot \vecnorm{\nabla F (\theta_{t
        - 1})}{2} + \vecnorm{\nabla F (\theta_{t - 1}) - \nabla F
      (\theta_{t - 2})}{2} \cdot \vecnorm{\nabla F (\theta_{t -
        1})}{2} \leq (1 + \stpsz \smoothness) \vecnorm{v_{t
        - 1}}{2} \cdot \vecnorm{\nabla F (\theta_{t - 1})}{2}.
\end{align*}
For the term $T_2$, we have:
\begin{align*}
     T_2 &\leq \Exs \left( \vecnorm{\noise(\theta_{t - 1}; \xi_t) -
       \noise (\theta_{t - 2}; \xi_t)}{2} \cdot \vecnorm{\noise
       (\theta_{t - 1}; \xi_{t})}{2} \mid \filtration_{t - 1} \right)
     \\ &\leq \sqrt{\Exs (\vecnorm{\noise(\theta_{t - 1}; \xi_t) -
         \noise (\theta_{t - 2}; \xi_t)}{2}^2 \mid \filtration_{t -
         1}) \cdot \Exs (\vecnorm{\noise (\theta_{t - 1},
         \xi_{t})}{2}^2 \mid \filtration_{t - 1})}\\ &\leq \sglip^2
     \stpsz \vecnorm{v_{t - 1}}{2} \cdot \vecnorm{\theta_{t
         - 1} - \thetastar}{2}\\ &\leq
     \frac{\sglip^2}{\strongconvex} \stpsz \vecnorm{v_{t
         - 1}}{2} \cdot \vecnorm{\nabla F (\theta_{t - 1} )}{2}
.
\end{align*}
So we have:
\begin{align*}
    I_3 &\leq \frac{3}{2} \Exs \left( \vecnorm{v_{t - 1}}{2} \cdot
    \vecnorm{\nabla F (\theta_{t - 1})}{2} \right) \leq \frac{3}{2} \sqrt{\Exs
      \vecnorm{v_{t - 1}}{2}^2 \cdot \Exs\vecnorm{\nabla F (\theta_{t
          - 1})}{2}^2}\\ &\leq \frac{t\stpsz
      \strongconvex}{8} \Exs\vecnorm{v_{ t- 1}}{2}^2 + \frac{9}{2t
      \strongconvex \stpsz} \Exs\vecnorm{\nabla F
      (\theta_{t - 1})}{2}^2.
\end{align*}
Putting above estimates together, we obtain:
\begin{align*}
    \Exs\vecnorm{v_t}{2}^2 &\leq \left(1 - \frac{1}{t} \right)^2 (1 -
    3\stpsz \strongconvex / 4) \Exs\vecnorm{v_{t - 1}}{2}^2
    + \frac{1}{t^2} \left( 4 \sigstar^2 + 2 \big( 1+ 
    \frac{\sglip^2}{\strongconvex^2}\big) \Exs\vecnorm{\nabla F
      (\theta_{t - 1} )}{2}^2 \right) \\ &\quad \quad + \frac{(t -
      1)\stpsz \strongconvex}{4 t} \Exs\vecnorm{v_{ t-
        1}}{2}^2 + \frac{9}{t^2 \strongconvex \stpsz} \Exs
    \vecnorm{\nabla F (\theta_{t - 1})}{2}^2\\ &\leq \left(1 -
    \frac{1}{t} \right)^2 \left( 1 - \frac{\stpsz
      \strongconvex}{2} \right) \Exs\vecnorm{v_{t - 1}}{2}^2 +
    \frac{10}{t^2 \strongconvex \stpsz} \Exs
    \vecnorm{\nabla F (\theta_{t - 1})}{2}^2 + \frac{4 \sigstar^2}{t^2}
,
\end{align*}
which completes the proof of this lemma.

\pb\subsection{Proof of Lemmas~\ref{lemm_ztsharpbdd_testmat} and~\ref{lemm_cross-term-bound-testmat}}\label{subsec:proof-test-mat-lemmas}

In this section, we present the proofs of lemma~\ref{lemm_ztsharpbdd_testmat} and~\ref{lemm_cross-term-bound-testmat}, the two technical lemmas involving a test matrix $\testMat \in \real^{d \times d}$.

\pb\subsubsection{Proof of Lemma~\ref{lemm_ztsharpbdd_testmat}}
The proof is similar to that of Lemma~\ref{lemm_ztsharpbdd}, and we follow the notation in such lemma throughout. Indeed, we note the following telescope result:
\begin{align*}
\niters^2 \Exs\vecnorm{\testMat z_\niters}{2}^2
-
\burnin^2 \Exs\vecnorm{\testMat z_{\burnin} }{2}^2
=
\sum_{t = \burnin + 1}^{\niters} \Exs\vecnorm{\testMat \noise_t(\thetastar)}{2}^2
+
\sum_{t = \burnin + 1}^{\niters} \Exs\vecnorm{\testMat \remnoise_t}{2}^2
+
2 \sum_{t = \burnin + 1}^{\niters} \Exs \inprod{\testMat \noise_t(\thetastar)}{ \testMat \remnoise_t}.
\end{align*}
Clearly, for each $t$, we have the following identity:
\begin{align*}
    \Exs\vecnorm{\testMat \noise_t(\thetastar)}{2}^2 = \mathrm{Tr} \big( \testMat \SigStar \testMat^\top \big).
\end{align*}
For the additional terms, we note that $\Exs\vecnorm{\testMat \remnoise_t}{2}^2 \leq \opnorm{\testMat}^2 \Exs\vecnorm{\remnoise_t}{2}^2$, and following the derivation in the proof of Lemma~\ref{lemm_ztsharpbdd}, we have the following identity:
\begin{align*}
&
    \sum_{t = \burnin + 1}^{\niters} \Exs
\inprod{\testMat \noise_t(\thetastar)}{\testMat \remnoise_t} 
\\&= 
\niters \cdot  \Exs \inprod{\testMat
\noise_\niters (\thetastar)}{\testMat \noise_\niters(\theta_{\niters - 1}) - \testMat\noise_\niters(\thetastar)} - \burnin \cdot \Exs \inprod{
\testMat \noise_{\burnin} (\thetastar)}{\testMat \noise_{\burnin}(\theta_{\burnin - 1}) - \testMat \noise_{\burnin}(\thetastar)}.
\end{align*}
Applying the Cauchy-Schwartz inequality, we obtain bounds similar to Eq.~\eqref{eq:zt-cross-term-bound-after-telescope}, for $t \in \{\burnin, \niters\}$:
\begin{align*}
\left|
t \cdot  \Exs \inprod{\testMat
\noise_t (\thetastar)}{\testMat \noise_t (\theta_{t - 1}) - \testMat\noise_t (\thetastar)}
\right|
&\leq
t \opnorm{\testMat}^2 \cdot \sqrt{\Exs\vecnorm{\noise_t(\thetastar)}{2}^2} \cdot \sqrt{\Exs\vecnorm{\noise_t(\theta_{t - 1}) - \noise_t(\thetastar)}{2}^2}
\\&\leq
\opnorm{\testMat}^2 \frac{t \sigstar \sglip}{\strongconvex} \sqrt{\Exs\vecnorm{\nabla F (\theta_{t - 1})}{2}^2}
.
\end{align*}
For the burn-in period, we have that:
\begin{align*}
\burnin^2 \Exs\vecnorm{\testMat z_{\burnin}}{2}^2 \leq 2 \burnin \Exs\vecnorm{\testMat \big( \noise_1 (\theta_0) - \noise_1 (\thetastar) \big)}{2}^2
+
2 \burnin \Exs\vecnorm{\testMat \noise_1 (\thetastar)}{2}^2
\leq
\frac{2 \burnin \sglip^2}{\strongconvex^2} \vecnorm{\nabla F(\theta_0)}{2}^2
+
2 \burnin \mathrm{Tr} \big( \testMat \SigStar \testMat \big)
.
\end{align*}
Putting them together, and following the derivation in Lemma~\ref{lemm_ztsharpbdd}, we obtain the conclusion of this lemma.

\pb\subsubsection{Proof of Lemma~\ref{lemm_cross-term-bound-testmat}}
The proof is similar to that of Lemma~\ref{lemm_cross-term-bound}. Following the notation in Lemma~\ref{lemm_cross-term-bound}, we have the decomposition:
\begin{align*}
\abss{\Exs \inprod{t \testMat z_t}{\testMat v_t}} 
\leq
(t - \Twindow) \abss{\Exs \inprod{\testMat z_{t - \Twindow}}{\testMat v_t}}
+
\abss{\Exs \inprod{\testMat ( t z_t -  (t - \Twindow) z_{t - \Twindow} )}{\testMat v_t}}
.
\end{align*}
Noting that
\begin{align*}
\abss{\Exs \inprod{\testMat ( t z_t -  (t - \Twindow) z_{t - \Twindow} )}{\testMat v_t}} 
\leq
\opnorm{\testMat}^2 \sqrt{\Exs\vecnorm{t z_t -  (t - \Twindow) z_{t - \Twindow}}{2}^2} \cdot \sqrt{\Exs\vecnorm{v_t}{2}^2}
,
\end{align*}
and that
\begin{align*}
\abss{\Exs \inprod{\testMat  z_{t - \Twindow}}{\testMat v_t}}
\leq
\opnorm{\testMat}^2 \sqrt{\Exs\vecnorm{z_{t - \Twindow}}{2}^2} \cdot \sqrt{\Exs\vecnorm{\Exs[ v_t \mid \filtration_{t - \Twindow}]}{2}^2}
.
\end{align*}
The rest of the proof simply follows that of Lemma~\ref{lemm_cross-term-bound}, with an additional factor of $\opnorm{\testMat}^2$ in each term.

\pb\section{Proofs of auxiliary lemmas in \S\ref{SecProofthm:asymptotic-const-step-size}}
\label{SecAsympAuxiliary}

In this section, we prove the three auxiliary lemmas used in the proof
of Proposition \ref{theo_asymptotic_iteration_const_step_size}. Note that the proofs
of the lemmas have inter-dependencies. In the following, we first prove
Lemma~\ref{lemma-vt-coupling-estimate} assuming
Lemma~\ref{lemma:hessian-noise-mixed-bound}, and then prove
Lemma~\ref{lemma:hessian-noise-mixed-bound} assuming
Lemma~\ref{lemma:yt-properties}. Finally, we give a
self-contained proof for Lemma~\ref{lemma:yt-properties}.

%%%%%%%%%%%%%%%%%%%%%%%%%%%%%%%%%%%%%%%%%%%%%%%%%%%%%%%%%%%%%%%%%%%%%%%%%%%%%%%%%%%%%%%%%%

\pb\subsection{Proof of Lemma~\ref{lemma-vt-coupling-estimate}}
\label{subsubsec:proof-coupling-estimate}

We begin by making note of the identities
\begin{align*}
t v_t 
&= 
(t - 1) ( v_{t - 1} + \nabla f(\theta_{t - 1}; \xi_t) - \nabla f(\theta_{t - 2}; \xi_t) ) + \nabla f(\theta_{t - 1}; \xi_t)
, \quad \mbox{and} \\
y_t 
&= 
y_{t - 1} - \stpsz \nabla^2 f (\thetastar; \xi_t) y_{t - 1} + \nabla f(\thetastar; \xi_t)
.
\end{align*}
Defining the quantity $e_t \mydefn t v_t - y_t$, we see that the two
identities above imply that
\begin{align*}
e_t & = e_{t - 1} + \left( (t - 1) (\nabla f(\theta_{t - 1}; \xi_t) -
\nabla f(\theta_{t - 2}; \xi_t) ) - \stpsz \nabla^2 f (\thetastar;
\xi_t) y_{t - 1} \right) + (\nabla f(\theta_{t - 1}; \xi_t) - \nabla
f (\thetastar; \xi_t) ) \\ & = \Quant_1(t) + \Quant_2(t) + \Quant_3(t),
\end{align*}
where we define
\begin{align*}
\Quant_1(t) & \mydefn e_{t - 1} - \stpsz \int_0^1 \nabla^2 f
(\rho \theta_{t - 1} + (1 - \rho) \theta_{t - 2}; \xi_t ) e_{t -
  1} d \rho, \quad \quad \Quant_2(t) \mydefn (\nabla f(\theta_{t -
  1}; \xi_t) - \nabla f(\thetastar; \xi_t) ), \\
\Quant_3(t) & \mydefn \stpsz \int_0^1 \left( \nabla^2 f (\rho
\theta_{t - 1} + (1 - \rho) \theta_{t - 2}; \xi_t ) - \nabla^2
f(\thetastar; \xi_t) \right) y_{t - 1} d \rho.
\end{align*}
By the triangle inequality, we have
\begin{align*}
  \Exs\vecnorm{e_t}{2}^2 & \leq \left( \sqrt{\Exs\vecnorm{\Quant_1
      (t)}{2}^2} + \sqrt{\Exs\vecnorm{\Quant_2 (t)}{2}^2} +
  \sqrt{\Exs\vecnorm{\Quant_3 (t)}{2}^2} \right)^2.
\end{align*}
In the following, we bound each term $\Exs\vecnorm{\Quant_i
  (t)}{2}^2$ in succession.

%%%%%%%%%%%%%%%%%%%%%%%%%%%%%%%%%%%%%%%%%%%%%%%%%%%%%%%%%%%%%%%%%%%%%%%%%%%%%%%%%%%%%
\paragraph{\textbf{Upper bound on $\Exs\vecnorm{\Quant_1(t)}{2}^2$:}}

Assumption~\ref{assu_StrcvxSmooth} and Assumption~\ref{assu_smoothnoise} together imply that
\begin{align*}
&\quad\,%%%%
\Exs\vecnorm{\Quant_1(t)}{2}^2 
%%%%%%%%
\\ &= 
\Exs\vecnorm{e_{t - 1}}{2}^2 - 2\stpsz \Exs \int_0^1 e_{t - 1}^\top \nabla^2 F (\rho \theta_{t - 1} + (1 - \rho) \theta_{t - 2} ) e_{t - 1} d \rho 
\\  &\quad \quad%%%%
+ \stpsz^2 \int_0^1 \Exs\vecnorm{\nabla^2 f (\rho \theta_{t - 1} + (1 - \rho) \theta_{t - 2}; \xi_t ) e_{t - 1}}{2}^2 d \rho 
%%%%%%%%
\notag \\&= 
\Exs\vecnorm{e_{t - 1}}{2}^2 - \Exs \int_0^1 e_{t - 1}^\top \left(2  \stpsz \nabla^2 F (\rho \theta_{t - 1} + (1 - \rho) \theta_{t - 2} ) - \stpsz^2 (\nabla^2 F (\rho \theta_{t - 1} + (1 - \rho) \theta_{t - 2} ) )^2 \right) e_{t - 1} d \rho
\notag \\&\quad \quad  %%%%
+ \stpsz^2 \int_0^1 \Exs\vecnorm{\NoiseAt (\rho \theta_{t - 1} + (1 - \rho) \theta_{t - 2} ) e_{t - 1}}{2}^2 d \rho
%%%%%%%%
\notag \\&\overset{(i)}{\leq} 
\Exs\vecnorm{e_{t - 1}}{2}^2 - (2 \stpsz - \stpsz^2 \smoothness) \int_0^1 e_{t - 1}^\top \nabla^2 F (\rho \theta_{t - 1} + (1 - \rho) \theta_{t - 2} ) e_{t - 1} d \rho + \stpsz^2 \sglip^2 \int_0^1 \vecnorm{e_{t - 1}}{2}^2 d \rho
%%%%%%%%
\notag \\&\overset{(ii)}{\leq} 
\Exs\vecnorm{e_{t - 1}}{2}^2 - \strongconvex \left(2 \stpsz - \stpsz^2 \smoothness\right) \Exs\vecnorm{e_{t - 1}}{2}^2 + \sglip^2 \stpsz^2 \Exs\vecnorm{e_{t - 1}}{2}^2
.
%%%%%%%%
\end{align*}
In step $(i)$, we are using the fact that $0 \preceq \nabla^2 F (\rho \theta_{t - 1} + (1 - \rho) \theta_{t - 2} ) \preceq \smoothness I_d$, and in step $(ii)$, we use the strong convexity of $F$.

For $\stpsz < \frac{1}{2 \smoothness} \wedge \frac{\strongconvex}{ 2\sglip^2}$,
we have $\Exs\vecnorm{\Quant_1 (t)}{2}^2 \leq (1 - \strongconvex
\stpsz ) \Exs\vecnorm{e_{t - 1}}{2}^2$.

%%%%%%%%%%%%%%%%%%%%%%%%%%%%%%%%%%%%%%%%%%%%%%%%%%%%%%%%%%%%%%%%%%%%%%%%%%%%%%%%%%%%%%%%

\paragraph{\textbf{Upper bound on $\Exs\vecnorm{\Quant_2(t)}{2}^2$:}}

By Assumption~\ref{assu_smoothnoise} and Eq.~\eqref{eq:theta-bound-in-proof-of-clt}, we have
\begin{align*}
\Exs\vecnorm{\Quant_2 (t)}{2}^2 \leq \sglip^2 \Exs\vecnorm{\theta_{t - 1}
  - \thetastar}{2}^2 \leq \frac{\constA \sglip^2}{t},
\end{align*}
where the last inequality follows from Theorem~\ref{theo_finitebdd_single}.

%%%%%%%%%%%%%%%%%%%%%%%%%%%%%%%%%%%%%%%%%%%%%%%%%%%%%%%%%%%%%%%%%%%%%%%%%%%%%%%%%%

\paragraph{\textbf{Upper bound on $\Exs\vecnorm{\Quant_3(t)}{2}^2$:}}

Applying Lemma~\ref{lemma:hessian-noise-mixed-bound} with
$\tilde{\theta}_{t - 1} \mydefn \rho \theta_{t - 1} + (1 - \rho)
\theta_{t - 2} \in \filtration_{t - 1}$, we have
\begin{multline*}
\Exs\vecnorm{\left( H_t (\rho \theta_{t - 1} + (1 - \rho) \theta_{t - 2}) - H_t (\thetastar) \right) y_{t - 1}}{2}^2 
\leq
\Contwo \sqrt{\Exs\vecnorm{\rho \theta_{t - 1} + (1 - \rho)
\theta_{t - 2}) - \thetastar }{2}^2 }
\\ \leq 
\Contwo \left(
\sqrt{\Exs\vecnorm{\theta_{t - 1} - \thetastar}{2}^2} +
\sqrt{\Exs\vecnorm{\theta_{t - 2} - \thetastar}{2}^2} \right)
\leq 
\Contwo \sqrt{\constA}\left( \frac{1}{\sqrt{t - 1}} + \frac{1}{\sqrt{t -
2}} \right) 
\leq 
\frac{16 \Contwo \sqrt{\constA}}{\sqrt{t}}
.
\end{multline*}
Putting the bounds for $(Q_1, Q_2, Q_3)$ together, we obtain:
\begin{align*}
\sqrt{\Exs\vecnorm{e_t}{2}^2} 
\leq 
\left( 1 - \frac{\stpsz \strongconvex}{2} \right) \sqrt{\Exs\vecnorm{e_{t - 1}}{2}^2} 
+ 
\frac{4 \Contwo^{1/2} \constA^{1/4}}{t^{1/4}} 
+ 
\frac{\sglip \sqrt{\constA}}{\sqrt{t}}
.
\end{align*}
Solving the recursion, we have
\begin{align*}
\sqrt{\Exs\vecnorm{e_T}{2}^2} 
\leq 
( 4 \Contwo^{1/2} \constA^{1/4} + \sglip \sqrt{\constA} ) \sum_{s = \burnin+1}^t s^{-
\frac{1}{4}} \exp \left( - \frac{\strongconvex \stpsz}{2} (T - s) \right) 
+ 
e^{- \frac{\strongconvex \stpsz (T - \burnin)}{2}} \sqrt{\Exs\vecnorm{e_{\burnin}}{2}^2}
.
\end{align*}
For the first term, we note that:
\begin{align*}
\sum_{s = \burnin+1}^T s^{- \frac{1}{4}} \exp \left( - \frac{\strongconvex \stpsz}{2} (T - s) \right) 
&\le
\sum_{s = 1}^{T / 2} \exp \left( - \frac{\strongconvex \stpsz}{2} T \right) 
+ 
\frac{1}{(T / 2)^{1/4}} \sum_{s = T / 2}^T e^{- \frac{\strongconvex \stpsz (T - s)}{2}}
\\ &\le 
\frac{T}{2} e^{- \frac{\strongconvex \stpsz T}{2}} + \frac{4}{\strongconvex \stpsz T^{1/4}}
.
\end{align*}

For $\niters$ large enough, the exponentially decaying term is dominated by the $\niters^{-1/4}$ term. 
So there exists a constant $\kcon > 0$, depending on the constants $(\constA, \Contwo, \aconprime, \stpsz, \strongconvex, \burnin)$ but independent of $t$, such that
\begin{align*}
\Exs\vecnorm{t v_t - y_t}{2}^2	\le	\frac{\kcon}{\sqrt{t}}
,
\end{align*}
which finishes the proof.

%%%%%%%%%%%%%%%%%%%%%%%%%%%%%%%%%%%%%%%%%%%%%%%%%%%%%%%%%%%%%%%%%%%%%%%%%%%%%%%%%%%%%%%%%%

\pb\subsection{Proof of Lemma~\ref{lemma:hessian-noise-mixed-bound}}\label{subsubsec:proof-hessian-noise-mixed-bound}

Observe that Assumption (CLT.A) guarantees that
\begin{align*}
\Exs \left( \vecnorm{ (H_t (\thetastar) - H_t (\tilde{\theta}_{t - 1}) ) y_{t - 1}}{2}^2 \,\bigg|\, \filtration_{t - 1} \right) 
\leq 
\beta^2 \vecnorm{\tilde{\theta}_{t - 1} - \thetastar}{2}^2 \cdot \vecnorm{y_{t - 1}}{2}^2
.
\end{align*}
On the other hand, by Assumption~\ref{assu_smoothnoise}, we have
\begin{align*}
\Exs \left( \vecnorm{ (H_t (\thetastar) - H_t (\tilde{\theta}_{t - 1}) ) y_{t - 1}}{2}^2 \,\bigg|\, \filtration_{t - 1} \right) 
\le
4\sglip^2 \vecnorm{y_{t - 1}}{2}^2
.
\end{align*}
Taking a geometric average and applying the tower law yields the bound
\begin{align*}
\Exs\vecnorm{ \left( H_t (\tilde{\theta}_{t - 1}) - H_t(\thetastar) \right) y_{t - 1}}{2}^2 
& \le 
2 \sglip \beta \Exs\left( \vecnorm{\tilde{\theta}_{t - 1} - \thetastar}{2} \cdot \vecnorm{y_{t - 1}}{2}^2 \right) 
\notag \\& \stackrel{(i)}{\leq} 
2 \sglip \beta \sqrt{\Exs\vecnorm{\tilde{\theta}_{t - 1} - \thetastar}{2}^2} \cdot \sqrt{ \Exs\vecnorm{y_{t - 1}}{2}^4}
,
\end{align*}
where step (i) follows from the Cauchy-Schwarz inequality. 
Applying Lemma~\ref{lemma:yt-properties}, we are guaranteed the existence of a constant $\aconprime > 0$ such that
\begin{align*}
\sup_{t \geq \burnin} \Exs\vecnorm{y_t}{2}^4 \leq \aconprime 
< 
\infty
.
\end{align*}
Setting $\Contwo = 2 \sglip \beta \sqrt{\aconprime}$ completes the proof of the claim.

%%%%%%%%%%%%%%%%%%%%%%%%%%%%%%%%%%%%%%%%%%%%%%%%%%%%%%%%%%%%%%%%%%%%%%%%%%%%%%%%%%%%%%%%%%%

\pb\subsection{Proof of Lemma~\ref{lemma:yt-properties}}\label{subsubsec:yt-properties}
Throughout this section, we adopt the shorthand notation $H_t \mydefn
H_t (\thetastar)$ and $\NoiseAt \mydefn \NoiseAt (\thetastar)$. We also use $\Xi$ to denote a generic random variable have the same law as $\Xi_1$.
Beginning with the proof of the first claim, we take expectations on
both sides of Eq.~\eqref{eq:yt-defn}, thereby finding that
\begin{align*}
 \Exs (y_t) = \Exs \left( y_{t - 1} - \stpsz H_t (\thetastar) y_{t
   - 1} + \noise_t (\thetastar) \right) = (I - \stpsz \hessianstar)
 \Exs (y_{t - 1}) = (I - \stpsz \hessianstar)^{t - \burnin} \Exs
 (y_{\burnin}) = 0.
\end{align*}
Our next step is to control the fourth moment. For $\stpsz \leq \frac{1}{2 \smoothness} < \frac{1}{2 \strongconvex}$, we observe that:
\begin{align*}
  \Exs\vecnorm{y_t}{2}^4 &= \Exs\vecnorm{ y_{t - 1} - \stpsz H_t
    (\thetastar) y_{t - 1} + \noise_t }{2}^4\\ &\leq \Exs\vecnorm{(I
    - \stpsz H_t) y_{t - 1}}{2}^4 + 4 \Exs (\vecnorm{(I - \stpsz
    H_t) y_{t - 1}}{2}^3 \cdot \vecnorm{\noise_t}{2}) + 6 \Exs (\vecnorm{(I - \stpsz H_t) y_{t -
      1}}{2}^2 \cdot \vecnorm{\noise_t}{2}^2) \\
  & \quad \quad + 4 \Exs (\vecnorm{\noise_t}{2}^3 \cdot \vecnorm{(I -
    \stpsz H_t) y_{t - 2}}{2}) + \Exs
  \vecnorm{\noise_t}{2}^4\\ 
  &\overset{(i)}{\leq} \left(1 + \frac{\stpsz
    \strongconvex}{2} \right) \Exs\vecnorm{(I - \stpsz H_t) y_{t -
      1}}{2}^4 + \frac{24}{(\stpsz \strongconvex)^3} \Exs
  \vecnorm{\noise_t}{2}^4 + \frac{216}{(\stpsz \strongconvex)^2}
  \Exs\vecnorm{\noise_t}{2}^4 + \frac{24}{(\stpsz \strongconvex)}
  \Exs\vecnorm{\noise_t}{2}^4 + \Exs\vecnorm{\noise_t}{2}^4\\ &\leq
  \left(1 + \frac{\stpsz \strongconvex}{2} \right) \Exs\vecnorm{(I
    - \stpsz H_t) y_{t - 1}}{2}^4 + \frac{157}{(\strongconvex
    \stpsz)^3} \Exs\vecnorm{\noise (\thetastar)}{2}^4,
\end{align*}
where in step $(i)$, we use Young's inequality for the last four terms.

Now we study the term $\Exs\vecnorm{(I - \stpsz H_t) y_{t -
    1}}{2}^4$. For $\stpsz < \frac{1}{\smoothness}$,
straightforward calculation yields:
\begin{align*}
    &\Exs \left( \vecnorm{(I - \stpsz H_t) y_{t - 1}}{2}^4 \mid
  \filtration_{t - 1} \right)\\
  &\leq \vecnorm{(I - \stpsz
    \hessianstar) y_{t - 1}}{2}^4 + 4 \Exs \left( \binprod{\stpsz    \NoiseAt y_{t - 1}}{(I - \stpsz \hessianstar) y_{t - 1}}
  \vecnorm{(I - \stpsz \hessianstar) y_{t - 1}}{2}^2 \mid
  \filtration_{t - 1}\right) + \Exs \left( \vecnorm{\stpsz \NoiseAt
    y_{t - 1}}{2}^4 \mid \filtration_{t - 1} \right)\\ 
    &\quad \quad +
  6 \Exs \left( \vecnorm{(I - \stpsz \hessianstar) y_{t - 1}}{2}^2
  \vecnorm{\stpsz \NoiseAt y_{t - 1}}{2}^2 \mid \filtration_{t - 1}
  \right) + 4 \Exs \left( \binprod{\stpsz \NoiseAt y_{t - 1}}{(I -    \stpsz \hessianstar) y_{t - 1}} \vecnorm{\stpsz \NoiseAt
    y_{t - 1}}{2}^2 \mid \filtration_{t - 1}\right)\\ 
    &\leq \vecnorm{(I - \stpsz
    \hessianstar) y_{t - 1}}{2}^4 + \stpsz^4 \Exs \left( \vecnorm{ \NoiseAt
    y_{t - 1}}{2}^4 \mid \filtration_{t - 1} \right) + 6 \stpsz^2 \sglip^2 \vecnorm{y_{t - 1}}{2}^4 \\
    & \quad \quad + 2 \Exs \left( \vecnorm{\stpsz \NoiseAt
    y_{t - 1}}{2}^4 \mid \filtration_{t - 1}\right) + 2 \Exs \left( \vecnorm{(I -
    \stpsz \hessianstar) y_{t - 1}}{2}^2 \cdot \vecnorm{\stpsz \NoiseAt
    y_{t - 1}}{2}^2 \mid \filtration_{t - 1}\right)  \\
    &\leq (1 - 3
  \stpsz \strongconvex) \vecnorm{y_{t - 1}}{2}^4 + 8 \stpsz^2
  \sglip^2 \vecnorm{y_{t - 1}}{2}^4 + 3 \stpsz^4 \sglip'^4
  \vecnorm{y_{t - 1}}{2}^4.
\end{align*}
For a step-size $\stpsz < \frac{1}{4 \smoothness} \wedge
\frac{\strongconvex}{16 \sglip^2} \wedge \frac{\strongconvex^{1/3}}{6
  \sglip'^{4/3}}$, we have $\Exs \left( \vecnorm{(I - \stpsz
  H_t) y_{t - 1}}{2}^4 \mid \filtration_{t - 1} \right) \leq (1 - 2
\strongconvex \stpsz) \vecnorm{y_{t - 1}}{2}^4$. Putting together
these bounds, we find that
\begin{align*}
    \Exs\vecnorm{y_t}{2}^4 \leq \left(1 - \strongconvex \stpsz
    \right) \Exs\vecnorm{y_{t - 1}}{2}^4 + \frac{157}{(\strongconvex
      \stpsz)^3} \Exs\vecnorm{\noise (\thetastar)}{2}^4,
\end{align*}
with the initial condition $\Exs\vecnorm{y_{\burnin}}{2}^4 =
0$. Solving this recursion leads to the bound
\begin{align*}
    \sup_{t \geq \burnin} \Exs\vecnorm{y_t}{2}^4 \leq \frac{157
    }{(\strongconvex \stpsz)^4} \Exs\vecnorm{\noise
      (\thetastar)}{2}^4.
\end{align*}
Let $\aconprime = \frac{157}{(\stpsz \strongconvex)^4}$, we prove the second claim.

Finally we study the stationary covariance of the process $\{y_t\}_{t \geq
  \burnin}$. The existence and uniqueness of the stationary
distribution was established in~\cite{mou2020linear}. Let
$\pi_\stpsz$ denote the stationary distribution of $(y_t)_{t \geq
  \burnin}$, and let $Q_\stpsz \mydefn \Exs_{Y \sim \pi_\stpsz}
(YY^\top)$. From the first part of this lemma, we can see that
$\Exs_{Y \sim \pi_\stpsz} (Y) = 0$. For $y_t \sim \pi_\stpsz$,
we have $y_{t + 1} \sim \pi_\stpsz$, and consequently,
\begin{align*}
  Q_\stpsz & = \Exs (y_{t + 1} y_{t + 1}^\top) \\
& = \Exs \left( (I - \stpsz H_{t + 1}) y_t y_t^\top (I - \stpsz
  H_{t + 1}^\top) + \noise_{t + 1} \noise_{t + 1}^\top \right) + \Exs
  \left( \noise_{t + 1} y_t^\top (I - \stpsz H_{t + 1}^\top) + (I -
  \stpsz H_{t + 1}) y_t \noise_{t + 1}^\top \right)\\
& = Q_\stpsz - \stpsz (\hessianstar Q_\stpsz + Q_\stpsz
  \hessianstar) + \stpsz^2 (\hessianstar Q_\stpsz \hessianstar +
  \Exs (\NoiseAplain Q_\stpsz \NoiseAplain)) + \SigStar.
\end{align*}
In the last equation, we use the fact that $\Exs (y_{t}) = 0$ and that $y_t$ is independent of $(H_{t + 1}, \noise_{t + 1})$, which leads to the following equation:
\begin{align*}
    \Exs \left( \noise_{t + 1} y_t^\top (I - \stpsz H_{t + 1}^\top) \right) = \Exs \left( \noise_{t + 1} (\thetastar) \otimes (I - \stpsz H_{t + 1} (\thetastar)) \right) \left[\Exs (y_{t}) \right] = 0. 
\end{align*}
Therefore, the matrix $Q_\stpsz$ satisfies the equation
\begin{align*}
\hessianstar Q_\stpsz + Q_\stpsz \hessianstar -
\stpsz(\hessianstar Q_\stpsz \hessianstar + \Exs (\NoiseAplain
Q_\stpsz \NoiseAplain)) = \frac{\SigStar}{\stpsz},
\end{align*}
which completes the proof of the last part of the lemma.

\end{document}